\pdfoutput=1
\documentclass[a4paper,10pt]{article}
\usepackage[margin=0.8in]{geometry}
\usepackage{graphicx,graphics}
\usepackage{amsfonts,amssymb}
\usepackage{amsmath,amsthm}
\usepackage[T1]{fontenc}
\usepackage[utf8]{inputenc}
\usepackage{algorithm,algorithmic}
\usepackage{bm}
\usepackage{subcaption}
\usepackage{placeins}
\usepackage{caption}
\usepackage{tikz}
\usepackage[hypertexnames=false]{hyperref}
\usepackage{sidecap}
\usepackage{soul}

\usepackage{color}

\newcommand{\iprod}[2]{\langle #1,#2\rangle}
\DeclareMathOperator{\Sym}{Sym}

\newcommand{\BDspace}{\mathvar{BD}}

\newcommand{\BVspace}{\mathvar{BV}}

\newcommand{\Meas}{\mathcal{M}}

\renewcommand{\L}{\mathcal{L}}

\renewcommand{\d}{\,d} %\text{d}}

\newcommand{\TGV}{\mathvar{TGV}}
\newcommand{\TV}{\mathvar{TV}}

\makeatletter
\def \weaktostar@sym{\setbox0=\hbox{$\rightharpoonup$}\rlap{\hbox 
        to\wd0{\hss\raise1ex\hbox{$\scriptscriptstyle{*\,}$}\hss}}\box0}
    \def \weaktostar    {\mathrel{\weaktostar@sym}}
\makeatother

\newcommand{\field}[1]{\mathbb{#1}}

\newcommand{\R}{\field{R}}

\newcommand{\norm}[1]{\|#1\|}

\newcommand{\abs}[1]{|#1|}

\newcommand{\grad}[1]{\nabla #1}
\newcommand{\freevar}{\,\boldsymbol\cdot\,}

\newcommand{\Union}\bigcup
\newcommand{\Isect}\bigcap
\newcommand{\union}\cup
\newcommand{\isect}\cap
\newcommand{\bigunion}\bigcup
\newcommand{\bigisect}\bigcap

\newcommand{\defeq}{:=}

\newcommand{\downto}{\searrow}
\newcommand{\upto}{\nearrow}

\DeclareMathOperator*{\argmin}{arg\,min}

\DeclareMathOperator{\Dom}{dom}

\makeatletter
\def \uminus@sym{\setbox0=\hbox{$\cup$}\rlap{\hbox 
        to\wd0{\hss\raise0.5ex\hbox{$\scriptscriptstyle{-}$}\hss}}\box0}
    \def \uminus    {\mathrel{\uminus@sym}}
\makeatother

\newcommand{\mathvar}[1]{\textup{#1}}

\theoremstyle{definition}
\newtheorem{assumptionx}{Assumption}

\DeclareMathOperator{\KRTV}{KRTV}
\DeclareMathOperator{\KR}{KR}

\def\DIMdomain{n}

\def\pixels{\ell}

\def\alphavec{\vec\alpha}

\def\costltwo{{\ensuremath{L_2^2}}}

\def\costhubertv{{\ensuremath{L_\eta^1\!\nabla}}}
\def\costsymbreg{{\ensuremath{B_sL_1\!\nabla_\eta}}}
\def\Nphi{M}%{M_\phi}
\def\NA{N}%{M_A}

\def\Smoother{H}
\newcommand{\JorigSYM}{J}%{^\infty}}

\newcommand{\JepsilonSYM}{J^{\gamma,\mu}}

\newcommand{\Jorig}[2][\lambda, \alpha]{\JorigSYM(#2; #1)}

\newcommand{\Jepsilon}[2][\lambda, \alpha]{\JepsilonSYM(#2; #1)}

\def\SolM{\mathcal{S}}

\def\nO{\mathfrak{N}}

\def\costf{F}

\newcommand{\huber}[2][\gamma]{\abs{#2}_{#1}}

\def\borel{\mathcal{B}}
\def\SPACEalpha{\mathcal{P}_\alpha^+}
\def\SPACEalphaPos{\SPACEalpha}

\def\SPACElambda{\mathcal{P}_\lambda^+}

\def\costK{}

\newcommand{\SPACEkuLp}[2][\DIMkurange]{L^{#2}(\Omega; \R^{#1})}

\newcommand{\SPACEkuLtwo}[1][\DIMkurange]{\SPACEkuLp[#1]{2}}

\newcommand{\ICTV}{\mathvar{ICTV}}
\newcommand{\POLYTV}{\mathvar{POLYTV}}

\newlength{\colw}

\def\iftgv#1#2#3{{\def\iftgvpar{#1}\def\iftgvtgv{tgv2}\ifx\iftgvpar\iftgvtgv#2\else#3\fi}}
\def\ifkrtv#1#2#3{{\def\ifkrtvpar{#1}\def\ifkrtvkrtv{krtv2}\ifx\ifkrtvpar\ifkrtvkrtv#2\else#3\fi}}
\def\ifictv#1#2#3{{\def\ifictvpar{#1}\def\ifictvictv{ictv}\ifx\ifictvpar\ifictvictv#2\else#3\fi}}
\def\ifpolytv#1#2#3{{\def\ifpolytvpar{#1}\def\ifpolytvpolytv{polytv}\ifx\ifpolytvpar\ifpolytvpolytv#2\else#3\fi}}
\def\iftwoparam#1#2#3{\iftgv{#1}{#2}{\ifictv{#1}{#2}{#3}}}
\def\ifhuber#1#2#3{{\def\ifhpar{#1}\def\ifhh{huberonly}\ifx\ifhpar\ifhh#2\else#3\fi}}
\def\ifsymbreg#1#2#3{{\def\ifhpar{#1}\def\ifhh{symbregman}\ifx\ifhpar\ifhh#2\else#3\fi}}
\def\ifssn#1#2#3{{\def\ifhpar{#1}\def\ifhh{ssn}\ifx\ifhpar\ifhh#2\else#3\fi}}
\def\mathname#1{\iftgv{#1}{$\TGV^2$}{\ifkrtv{#1}{$\KRTV$}{\ifictv{#1}{$\ICTV$}{\ifpolytv{#1}{$\POLYTV_1$}{$\TV$}}}}}
\def\costname#1{\ifhuber{#1}{\costhubertv}{\ifsymbreg{#1}{\costsymbreg}{\costltwo}}}

\def\TVINIT{(\alpha_{\TV}^*/\pixels, \alpha_{\TV}^*)}
\def\TVINITKRTV{(\alpha_{\TV}^*/\pixels^{1.5}, \alpha_{\TV}^*)}

\def\XXint#1#2#3{{\setbox0=\hbox{$#1{#2#3}{\int}$ }
\vcenter{\hbox{$#2#3$ }}\kern-.6\wd0}}

\newcommand{\resplotxx}[2][]{
    \setlength{\colw}{\textwidth}
    \begin{tikzpicture}
        \pgftext[at=\pgfpoint{0}{0},left,bottom]{%
            \includegraphics[width=\colw]{{#2}.png}
        }
        #1
    \end{tikzpicture}%
}

\newlength{\imw}
\newcommand{\resplot}[5][]{
    \begin{subfigure}[t]{\imw}%
    \resplotxx[#1]{tgv-learn-resimg/#2}
    \caption{\mathname{#3} denoising, \costname{#4} cost}
    \ifdefined\subfigprefix\label{\subfigprefix:#3-#4}\else\relax\fi 
    \end{subfigure}
    }
\newcommand{\resplottexture}[5][]{
    \begin{subfigure}[t]{\imw}%
    \resplotxx[#1]{tgv-learn-resimg/#2}
    \caption{Texture component for \mathname{#3}}
    \ifdefined\subfigprefix\label{\subfigprefix:#3-#4}\else\relax\fi 
    \end{subfigure}
    }
\newcommand{\inplot}[3][]{
    \begin{subfigure}[t]{\imw}%
    \resplotxx[#1]{tgv-learn-img/#2}
    \caption{#3}
    \ifdefined\subfigprefix\label{\subfigprefix:#2}\else\relax\fi
    \end{subfigure}
    }

\def\SQl{48}
\def\SQb{160}
\def\SQr{112}
\def\SQt{224}
\newcommand{\igraph}[2]{\includegraphics[#1,bb=48 160 112 224,clip]{{tgv-learn-resimg/#2}.png}}
\newlength{\scf}

\newcommand{\drawzoomarea}{\draw[line width=1.5,dashed,color=red] (\SQl\scf, \SQb\scf) rectangle (\SQr\scf, \SQt\scf);}

\newcommand{\resplotzxx}[2][]{
    \setlength{\colw}{\textwidth}
    \setlength{\scf}{0.0039\colw} % \colw/256
    \begin{tikzpicture}
        \pgftext[at=\pgfpoint{0}{0},left,bottom]{%
            \includegraphics[width=\colw]{{tgv-learn-resimg/#2}.png}
        }
        \pgftext[at=\pgfpoint{0.55\colw}{0.05\colw},left,bottom]{%
            \igraph{width=0.4\colw}{#2}
        }
        \draw[line width=1.5, color=red] (0.55\colw, 0.05\colw) rectangle (0.95\colw, 0.45\colw);
        #1
    \end{tikzpicture}%
}

\newcommand{\restabll}[6][]{%
    \input{tgv-learn-resimg/#2-vals.tex}%
    \mathname{#3} &  % Denoising
    \costname{#4} & % Cost
    #5 & % Initial $\alphavec$
    \ifpolytv{#3}{%
        $(\RESalphaSCone,\RESalphaSCtwo,\RESalphaSCthree,\RESalphaSCfour)/\pixels$%
    }{%
        \iftwoparam{#3}%
              {$(\RESalphaSCtwo/\pixels^2, \RESalphaSCone/\pixels)$}%
              {$\RESalphaSCone/\pixels$}% Result $\alphavec$
    } &
   \RESdist & % Value
   \RESssim & % SSIM
   \RESpsnr &
   \RESiters &
   \def\arg{#6}
   \def\nofig{}
   \ifx\arg\nofig\else\ref{#6}(\subref{#6:#3-#4})\fi \\%
}

\newtheorem{theorem}{Theorem}[section]

\newtheorem{definition}[theorem]{Definition}

\graphicspath{ {./images/}{./tgv-learn-resimg/}{./tgv-learn-img/} }

\begin{document}

\begin{center}
{\Large \bf Bilevel approaches for learning of variational imaging models} \\[1cm]
\end{center}

\begin{center}
     \small{            {\sc Luca Calatroni}\\
 Cambridge Centre for Analysis, University of Cambridge\\
 Wilberforce Road, CB3 0WA, Cambridge, UK \\
 (lc524@cam.ac.uk) \\[0.4cm]
 {\sc Cao Chung}\\
 Research Center on Mathematical Modelling (MODEMAT), \\
 Escuela Polit\'ecnica Nacional, Quito, Ecuador \\
 (cao.vanchung@epn.edu.ec) \\[0.4cm]
 {\sc Juan Carlos De Los Reyes}\\
 Research Center on Mathematical Modelling (MODEMAT), \\
 Escuela Polit\'ecnica Nacional, Quito, Ecuador \\
 (juan.delosreyes@epn.edu.ec) \\[0.4cm]
 {\sc Carola-Bibiane Sch\"{o}nlieb} \\
 Department of Applied Mathematics and Theoretical Physics (DAMTP)\\
 University of Cambridge, Wilberforce Road, CB3 0WA, Cambridge, UK \\
 (cbs31@cam.ac.uk)\\ [0.4cm]
  {\sc Tuomo Valkonen} \\
 Department of Applied Mathematics and Theoretical Physics (DAMTP)\\
 University of Cambridge, Wilberforce Road, CB3 0WA, Cambridge, UK \\
 (tuomo.valkonen@iki.fi)}\\ [0.4cm]
       \end{center}

\begin{abstract}
We review some recent learning approaches in variational imaging, based on bilevel optimisation, and emphasize the importance of their treatment in function space. The paper covers both analytical and numerical techniques. Analytically, we include results on the existence and structure of minimisers, as well as optimality conditions for their characterisation. Based on this information, Newton type methods are studied for the solution of the problems at hand, combining them with sampling techniques in case of large databases. The computational verification of the developed techniques is extensively documented, covering instances with different type of regularisers, several noise models, spatially dependent weights and large image databases.
\end{abstract}


\vspace{-0.3cm}

\paragraph{Classification:} 49J40, 49J21, 49K20, 68U10, 68T05, 90C53, 65K10

%\keywords{Image denoising, variational methods, bilevel optimisation, supervised learning}

%\classification{49J40, 49J21, 49K20, 68U10, 68T05, 90C53, 65K10}

%-------------------------------------------------------------------------%
% To include a table of chapters, write                                   %
%   \tableofcontents
%-------------------------------------------------------------------------%

\section{Overview of learning in variational imaging}
A myriad of different imaging models and reconstruction methods exist in the literature and their analysis and application is mostly being developed in parallel in different disciplines. The task of image reconstruction from noisy and under-sampled measurements, for instance, has been attempted in engineering and statistics (in particular signal processing) using filters \cite{natterer2001mathematical,scharf1991statistical,cichocki2002adaptive} and multi scale analysis \cite{unser1995texture,kingsbury2001complex,unser2013unifying}, in statistical inverse problems using Bayesian inversion and machine learning \cite{evans2002inverse} and in mathematical analysis using variational calculus, PDEs and numerical optimisation \cite{rudin1992nonlinear}. Each one of these methodologies has its advantages and disadvantages, as well as multiple different levels of interpretation and formalism. In this paper we focus on the formalism of variational reconstruction approaches. 

A variational image reconstruction model can be formalised as follows. Given data $f$ which is related to an image (or to certain image information, e.g. a segmented or edge detected image) $u$ through a generic forward operator (or function) $K$ the task is to retrieve $u$ from $f$. In most realistic situations this retrieval is complicated by the ill-posedness of $K$ as well as random noise in $f$. A widely accepted method that approximates this ill-posed problem by a well-posed one and counteracts the noise is the method of Tikhonov regularisation. That is, an approximation to the true image is computed as a minimiser of 
\begin{equation}\label{reg1:eq}
\alpha~ R(u) + d(K(u),f),
\end{equation}
where $R$ is a regularising energy that models a-priori knowledge about the image $u$, $d(\cdot , \cdot)$ is a suitable distance function that models the relation of the data $f$ to the unknown $u$, and $\alpha>0$ is a parameter that balances our trust in the forward model against the need of regularisation. The parameter $\alpha$ in particular, depends on the amount of ill-posedness in the operator $K$ and the amount (amplitude) of the noise present in $f$. A key issue in imaging inverse problems is the correct choice of $\alpha$, image priors (regularisation functionals $R$), fidelity terms $d$ and (if applicable) the choice of what to measure (the linear or nonlinear operator $K$). Depending on this choice, different reconstruction results are obtained. 

Several strategies for conceiving optimization problems have been considered. One approach is the a-priori modelling of image priors, forward operator $K$ and distance function $d$. Total variation regularisation, for instance, has been introduced as an image prior in \cite{rudin1992nonlinear} due to its edge-preserving properties. Its reconstruction qualities have subsequently been thoroughly analysed in works of the variational calculus and partial differential equations community, e.g. \cite{ambcosdal96,caselles2007discontinuity,allard2008total,benning2012ground,bredies2013properties,benning2011higher,papafitsoros2013study,tuomov-jumpset} only to name a few. The forward operator in magnetic resonance imaging (MRI), for instance, can be derived by formalising the physics behind MRI which roughly results in $K=\mathcal S\mathcal F$ a sampling operator applied to the Fourier transform. Appropriate data fidelity distances $d$ are mostly driven by statistical considerations that model our knowledge of the data distribution \cite{idier2013bayesian,kaipio2006statistical}. Poisson distributed data, as it appears in photography \cite{costantini2004virtual} and emission tomography applications \cite{vardi1985statistical}, is modelled by the Kullback-Leibler divergence \cite{sawatzky2009total}, while a normal data distribution, as for Gaussian noise, results in a least squares fit model. In the context of data driven learning approaches we mention statistically grounded methods for optimal model design \cite{haber2003learning} and marginalization \cite{bui2008model,kolehmainen2011marginalization}, adaptive and multiscale regularization \cite{tadmor2004multiscale,dong2011automated,frick2012statistical}, learning in the context of sparse coding and dictionary learning \cite{OF96,MBPSZ08,mairal2009online,YSM10,peyre2011learning}, learning image priors using Markov Random fields \cite{roth2005fields,tappen2007utilizing,domke2012generic}, deriving optimal priors and regularised inversion matrices by analysing their spectrum \cite{chung2011designing,gilboa2014total}, and many recent approaches that -- based on a more or less generic model setup such as \eqref{reg1:eq} -- aim to optimise operators (i.e., matrices and expansion) and functions (i.e. distance functions $d$) in a functional variational regularisation approach by bilevel learning from `examples'  \cite{huang2012optimal,de2013image,kunisch2013bilevel,baus2014fully,schmidt2014shrinkage,chen_cvpr2015,reyes2015a,reyes2015b}, among others. All these approaches vary in their philosophy and mathematics. The main dividing line is between model-based derivation of \eqref{reg1:eq} that uses a-priori knowledge on the data and the image, and data-based derivation of \eqref{reg1:eq} that learns the setup of the model from the data itself.

While functional modelling constitutes a mathematically rigorous and physical way of setting up the reconstruction of an image -- providing reconstruction guarantees in terms of error and stability estimates -- it is limited with respect to its adaptivity for real data. On the other hand, data-based modelling of reconstruction approaches is set up to produce results which are optimal with respect to the given data. However, in general it neither offers insights into the structural properties of the model nor provides comprehensible reconstruction guarantees. Indeed, we believe that for the development of reliable, comprehensible and at the same time effective models \eqref{reg1:eq} it is essential to aim for a unified approach that seeks tailor-made regularisation and data models by combining model- and data-based approaches.

To do so we focus on a bilevel optimisation strategy for finding an optimal setup of variational regularisation models \eqref{reg1:eq}. That is, given a set of training images we find a setup of \eqref{reg1:eq} which minimises an a-priori determined cost functional $F$ measuring the performance of \eqref{reg1:eq} with respect to the training set, compare Section \ref{sec:analysis} for details. The setup of \eqref{reg1:eq} can be optimised for the choice of regularisation $R$ as will be discussed in Section \ref{sec:image}, for the data fitting distance $d$ as in Section \ref{sec:data}, or for an appropriate forward operator $K$ as in blind image deconvolution \cite{hintermuller2014bilevel} for example. In the present article, rather than working on the discrete problem, as is done in standard parameter learning and model optimisation methods, we discuss the optimisation of variational regularisation models in infinite dimensional function space. We will explain this approach in more detail in the next section. Before, let us give an account to the state of the art of bilevel optimisation for model learning. In machine learning bilevel optimisation is well established. It is a semi-supervised learning method that optimally adapts itself to a given dataset of measurements and desirable solutions. In \cite{roth2005fields,tappen2007utilizing,domke2012generic,chen2014}, for instance the authors consider bilevel optimization for finite dimensional Markov random field models. In inverse problems the optimal inversion and experimental acquisition setup is discussed in the context of optimal model design in works by Haber, Horesh and Tenorio \cite{haber2003learning,haber2010numerical}, as well as Ghattas et al. \cite{bui2008model,biegler2011large}.  Recently parameter learning in the context of functional variational regularisation models \eqref{reg1:eq} also entered the image processing community with works by the authors \cite{de2013image,calatronidynamic,reyes2015a,reyes2015b,lucainfimal,chungdelosreyes}, Kunisch, Pock and co-workers \cite{kunisch2013bilevel,Chen2012,chen_cvpr2015}, Chung et al. \cite{chung2014optimal}, and others \cite{baus2014fully,schmidt2014shrinkage}. Interesting recent works also include bilevel learning approaches for image segmentation \cite{Ranftl_GCPR2014} and learning of support vector machines \cite{Klatzer2015}.

\begin{figure}[t]
\centering
\includegraphics[width=0.7\textwidth]{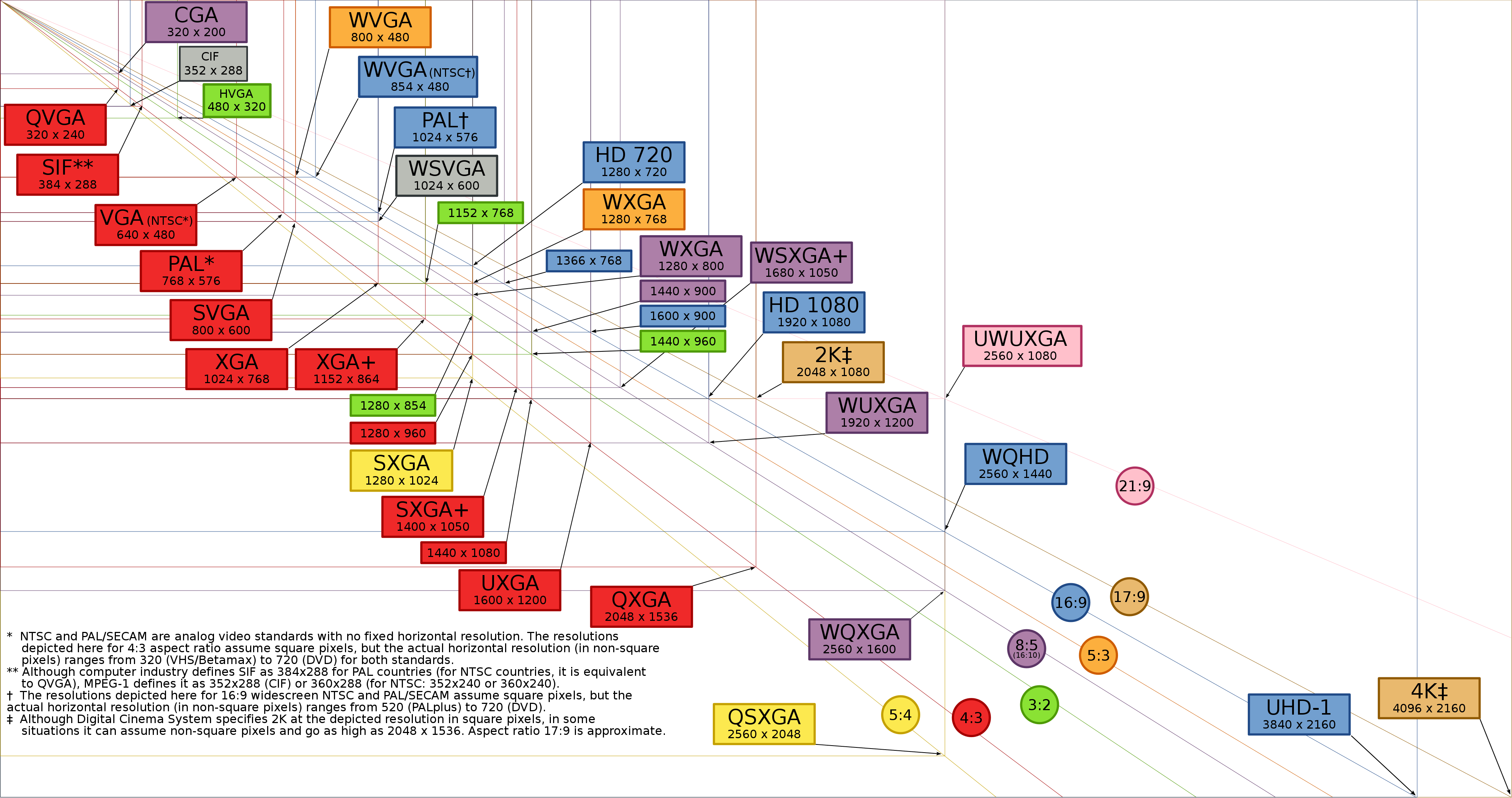}
\caption{Camera technology tending towards continuum images? Most image processing and analysis algorithms are designed for a finite number of pixels. But camera technology allows to capture images of higher and higher resolution and therefore the number of pixels in images changes constantly. Functional analysis, partial differential equations and continuous optimisation allow us to design image processing models in the continuum.}
\label{fig:resolutionchart}
\end{figure}

Apart from the work of the authors \cite{de2013image,calatronidynamic,reyes2015a,reyes2015b,chungdelosreyes,lucainfimal}, all approaches so far are formulated and optimised in the discrete setting. In what follows, we review modelling, analysis and optimisation of bilevel learning approaches in function space rather than on a discretisation of \eqref{reg1:eq}. While digitally acquired image data is of course discrete, the aim of high resolution image reconstruction and processing is always to compute an image that is close to the real (analogue, infinite dimensional) world. HD photography produces larger and larger images with a frequently increasing number of megapixels, compare Figure \ref{fig:resolutionchart}. Hence, it makes sense to seek images which have certain properties in an infinite dimensional function space. That is, we aim for a processing method that accentuates and preserves qualitative properties in images independent of the resolution of the image itself \cite{viola2012unifying}. Moreover, optimisation methods conceived in function space potentially result in numerical iterative schemes which are resolution and mesh-independent upon discretisation \cite{hintermuller2006infeasible}. 

\paragraph{Outline of the paper} In what follows we focus on bilevel learning of an optimal variational regularisation model in function space. We give an account on the analysis for a generic learning approach in infinite dimensional function space presented in \cite{reyes2015a} in Section \ref{sec:analysis}. In particular, we will discuss under which conditions on the learning approach, in particular the regularity of the variational model and the cost functional, we can indeed prove existence of optimal parameters in the interior of the domain (guaranteeing compactness), and derive an optimal system exemplarily for parameter learning for total variation denoising. Section \ref{sec:algorithms} discusses the numerical solution of bilevel learning approaches. Here, we focus on the second-order iterative optimisation methods such as quasi and semismooth Newton approaches \cite{de2015numerical}, which are combined with stochastic (dynamic) sampling strategies for efficiently solving the learning problem even in presence of a large training set \cite{calatronidynamic}. In Sections \ref{sec:image} and \ref{sec:data} we discuss the application of the generic learning model from Section \ref{sec:analysis} to conceiving optimal regularisation functionals (in the simplest setting this means computing optimal regularisation parameters; in the most complex setting this means computing spatially dependent and vector valued regularisation parameters) \cite{reyes2015b,chungdelosreyes} and optimal data fidelity functions in presence of different noise distributions \cite{de2013image,lucainfimal}.

\section{The learning model and its analysis in function space}\label{sec:analysis}
%!TEX root = article.tex

\subsection{The abstract model}
\label{subsec:abstract model}

Our image domain will be an open bounded set $\Omega \subset \R^\DIMdomain$ with Lipschitz boundary. Our data $f$ lies in $Y=L^2(\Omega; \R^m)$.
We look for positive parameters $\lambda=(\lambda_1,\ldots,\lambda_\Nphi)$ and $\alpha=(\alpha_1,\ldots,\alpha_\NA)$ in abstract parameters sets $\SPACElambda$ and $\SPACEalpha$ 
%\TODO{These are not spaces!!!!!! With positivity constraints are not spaces. That's why the difference}. \TODO{---Naah. Any subset of a metric space or a topological space is still a space of the same kind. We're not specifying here that we are in vector spaces or anything. In practise, of course the formulation of $J$ below forces these spaces/sets to be subsets of $L^2$. But that would not have to be so with a more abstract formulation, $J(u)=\Phi(Ku; \lambda)+R(u; \alpha)$. Which we could very well do here, to remove excess detail from the general model} 
They are intended to solve for some convex, proper, weak* lower semicontinuous cost functional $F: X \to \R$ the problem
\begin{equation}
    \label{eq:learn}
    \tag{$\mathrm{P}$}
    \min_{\alpha \in \SPACEalpha,\, \lambda \in \SPACElambda} \costf(\costK u_{\alpha,\lambda})
    \quad
    \text{s.t.}
    \quad
    u_{\alpha,\lambda} \in \argmin_{u \in X} \Jorig{u},
\end{equation}
for
\[
    \Jorig{u} \defeq
        \sum_{i=1}^{\Nphi} \int_\Omega \lambda_i(x) \phi_i(x, [Ku](x)) \d x
        +
        \sum_{j=1}^{\NA} \int_\Omega \alpha_j(x) \d \abs{A_j u}(x).
\]
Our solution $u$ lies in an abstract space $X$, mapped by the linear operator $K$ to $Y$.
Several further technical assumptions discussed in detail in \cite{reyes2015a} cover $A$, $K$, and the $\phi_i$. In Section \ref{sec:l2theory} of this review we concentrate on specific examples of the framework. %Before that we discuss some modifications required for the practical numerical realisation of the bi-level framework.

%\subsection{Numerical considerations}
For the approximation of problem \eqref{eq:learn} we consider various smoothing steps.
%The numerical solution of the denoising sub-problem, using an infeasible semi-smooth quasi-Newton approach \cite{hintermuller2006infeasible} depends on various smoothing steps. 
For one, we require Huber regularisation of the Radon norms.
Secondly, we take a convex, proper, and weak* lower-semicontinous smoothing functional $\Smoother: X \to [0,\infty]$.
The typical choice that we concentrate on is $H(u)=\frac{1}{2}\norm{\grad u}^2$.

For parameters $\mu \ge 0$ and $\gamma \in (0, \infty]$, we then consider the problem
\begin{equation}
    \label{eq:learn-numerical-single}
    \tag{$\mathrm{P}^{\gamma,\mu}$}
    \min_{\alpha \in \SPACEalpha,\, \lambda \in \SPACElambda} \costf(\costK u_{\alpha,\lambda,\gamma,\mu})
    \quad
    \text{s.t.} 
    \quad
    u_{\alpha,\lambda,\gamma,\mu} \in \argmin_{u \in X \isect \Dom \mu\Smoother } \Jepsilon{u}
\end{equation}
for
\begin{equation}
    %\label{eq:j-hilbert}
    %\tag{$J$}
    \notag
    \Jepsilon{u} :=
        \mu \Smoother(u)
        +
        \sum_{i=1}^{\Nphi} \int_\Omega \lambda_i(x) \phi_i(x, [Ku](x)) \d x
        +
        \sum_{j=1}^{\NA} \int_\Omega \alpha_j(x) \d \abs{A_j u}_\gamma(x).
\end{equation}
Here we denote by $\abs{A_j u}_\gamma$ the Huberised total variation measure per the following definition.

\begin{definition}\label{def:huber}
Given $\gamma \in (0, \infty]$, we  define for the norm $\norm{\freevar}_2$ on $\R^\DIMdomain$, 
the Huber regularisation
\[
    \huber[\gamma]{g} = 
    \begin{cases}
        \norm{g}_2 - \frac{1}{2\gamma}, & \norm{g}_2 \ge 1/\gamma,
        \\
        \frac{\gamma}{2}\norm{g}_2^2, & \norm{g}_2 < 1/\gamma.
    \end{cases}
\]
% We observe that this can equivalently be written using convex conjugates as
% \begin{equation}
%     \label{eq:huber-twonorm}
%     \alpha \huber[\gamma]{g}=
%     \sup
%     \Bigl\{
%         \iprod{q}{g} - \frac{1}{2\gamma\alpha} \norm{q}_2^2
%         \Bigm|
%         \norm{q}_2 \le \alpha
%     \Bigr\}.
% \end{equation}
Then if $\nu = f \L^\DIMdomain + \nu^s$ is the Lebesgue decomposition of $\nu \in \Meas(\Omega; \R^\DIMdomain)$ into the absolutely continuous part $f\L^\DIMdomain$ and the singular part $\nu^s$, we set
\[
    \huber[\gamma]{\nu}(V) \defeq \int_V \huber[\gamma]{f(x)} \d x + \abs{\nu^s}(V),
    \quad
    (V \in \borel(\Omega)).
\]
The measure $\huber[\gamma]{\nu}$ is the Huber-regularisation of the total variation measure $\abs{\nu}$. 
\end{definition}
In all of these, we interpret the choice $\gamma=\infty$ to give back the standard unregularised total variation measure or norm.

\subsection{Existence and structure: $L^2$-squared cost and fidelity}
\label{sec:l2theory}

We now choose
\begin{equation}
    \label{eq:l2setup}
    \costf(u)=\frac{1}{2}\norm{K u - f_0}_Y^2,
    \quad
    \text{and}
    \quad
    \phi_1(x, v)=\frac{1}{2}\abs{f(x)-v}^2,
\end{equation}
with $\Nphi=1$.
We also take $\SPACElambda=\{1\}$, i.e., we do not look for the fidelity weights.
%It is easily verified that the latter fidelity satisfies Assumption \ref{ass:phi}.
Our next results state for specific regularisers with discrete parameters $\alpha=(\alpha_1,\ldots,\alpha_\NA) \in \SPACEalpha = [0,\infty]^\NA$, conditions for the optimal parameters to satisfy $\alpha>0$. Observe how we allow infinite parameters, which can in some cases distinguish between different regularisers.

We note that these results are not a mere existence results; they are structural results as well. If we had an additional lower bound $0 < c \le \alpha$ in \eqref{eq:learn}, we could without the conditions \eqref{eq:l2cost-interior-condition-tv} for $\TV$ and \eqref{eq:l2cost-interior-condition-tgv2} for $\TGV^2$ \cite{bredies2011tgv} denoising, show the existence of an optimal parameter $\alpha$. Also with fixed numerical regularisation ($\gamma<\infty$ and $\mu>0$), it is not difficult to show the existence of an optimal parameter without the lower bound. What our very natural conditions provide is existence of optimal interior solution $\alpha>0$ to \eqref{eq:learn} without any additional box constraints or the numerical regularisation. Moreover, the conditions \eqref{eq:l2cost-interior-condition-tv} and \eqref{eq:l2cost-interior-condition-tgv2} guarantee convergence of optimal parameters of the numerically regularised $H^1$ problems \eqref{eq:learn-numerical-single} to a solution of the original $\BVspace(\Omega)$ problem \eqref{eq:learn}.

% \begin{theorem}[\cite{reyes2015a}]
%     \label{thm:l2cost-main}
%     Let $Y$ and $Z$ be Hilbert spaces, $\costf(u)=\frac{1}{2}\norm{K_0 u - f_0}_Z^2$, and $\Phi(v)=\frac{1}{2}\norm{f-v}_Y^2$ for some $f \in \range{K}$, $f_0 \in Z$, and a bounded linear operator $K_0: X \to Z$ satisfying
%     \begin{equation}
%         \label{eq:k0-condition}
%     \norm{K_0 u}_Z \le C_0 \norm{Ku}_Y, \quad \text{for all } u \in X
%         \quad \text{for some constant $C_0>0$}.
%     \end{equation}
%     Suppose Assumption \ref{ass:a-k} and \ref{ass:aj-zeroterm-approx} hold.
%     If for some $\bar\alpha \in \SPACEalphaInt$ and $t \in (0, 1/C_0]$ holds
%     \begin{equation}
%         \label{eq:l2cost-interior-condition-k0}
%         \regfMARG{f} > \regfMARG{f-t(K_0\invK)^*(K_0 \invK f-f_0)},
%     \end{equation}
%     then the problem  \eqref{eq:learn}, admits a solution $\hat \alpha \in \SPACEalphaCompactPosInt$.
% \end{theorem}

% \begin{remark}
%     \label{remark:l2cost-interior-condition}
%     Let $Z=Y$ and $K_0=K$ in Theorem \ref{thm:l2cost-main}. Then \eqref{eq:l2cost-interior-condition-k0}
%     reduces into
%     \begin{equation}
%         \label{eq:l2cost-interior-condition}
%         \regfMARG{f} > \regfMARG{f_0}.
%     \end{equation}
%     Also observe that our result requires, through \eqref{eq:k0-condition}, $\Phi \circ K$ to measure all the data that $F$ measures; otherwise for example an oscillating solution $u_\alpha$ for $\alpha \in \BD \SPACEalphaCompactPos$ could have low cost.
% \end{remark}

\begin{theorem}[Total variation Gaussian denoising \cite{reyes2015a}]
    \label{theorem:tv-l2cost-interior}
    %In case of the total variation denoising of Example \ref{example:tv-1}, 
    Suppose $f, f_0 \in \BVspace(\Omega) \isect L^2(\Omega)$, and
    \begin{equation}
        \label{eq:l2cost-interior-condition-tv}
        \TV(f)
        >
        \TV(f_0).
    \end{equation}
    Then there exist $\bar \mu, \bar\gamma > 0$ such that any optimal solution $\alpha_{\gamma,\mu} \in [0, \infty]$
    to the problem
    \[
        \min_{\alpha \in [0, \infty]} \frac{1}{2}\norm{f_0-u_\alpha}_{L^2(\Omega)}^2
    \]
    with
    \[
        %\begin{split}
        u_\alpha \in \argmin_{u \in \BVspace(\Omega)} 
            \Bigl(
            %&
            \frac{1}{2}\norm{f-u}_{L^2(\Omega)}^2
            + \alpha\abs{Du}_\gamma(\Omega)
            + \frac{\mu}{2}\norm{\grad v}_{L^2(\Omega; \R^\DIMdomain)}^2
            \Bigr)
        %\end{split}
    \]
    satisfies $\alpha_{\gamma,\mu} > 0$ whenever $\mu \in [0, \bar \mu]$, $\gamma \in [\bar\gamma, \infty]$.
\end{theorem}

This says that for the optimal parameter to be strictly positive, the noisy image $f$ should, in terms of the total variation, oscillate more than the noise-free image $f_0$ -- exactly what we would naturally expect!

\begin{proof}[First steps of proof: modelling in the abstract framework]
    The modelling of total variation is is based on the choice of $K$ as the embedding of $X=\BVspace(\Omega) \isect L^2(\Omega)$ into $Y=L^2(\Omega)$, and $A_1=D$. For the smoothing term we take $H(u)=\frac{1}{2}\norm{\grad v}_{L^2(\Omega; \R^\DIMdomain)}^2$.
    For the rest of the proof we refer to \cite{reyes2015a}.
\end{proof}

% \begin{proof}
%     Assumption \ref{ass:a-k}, \ref{ass:aj-zeroterm-approx}, and \ref{ass:sol-smoothness} we
%     have already verified in Example \ref{example:tv-1}.
%     We then observe that $K_0=K$, so we are in the setting of Remark \ref{remark:l2cost-interior-condition}.
%     Following the mapping of the TV problem to the general framework using the construction in Example \ref{example:tv-1}, we have $K=I$ and $\invK=I$ embeddings with $Y = L^2(\Omega)$.
%     $\invK$ is bounded on $\range{K}=L^2(\Omega) \isect \BVspace(\Omega)$.
%     Moreover, by Example \ref{example:tv-marginal}, $\regfMARG{v}=\bar\alpha\TV(v)$.
%     Thus \eqref{eq:l2cost-interior-condition} reduces with the choice $t=1$ into \eqref{eq:l2cost-interior-condition-tv}.
% \end{proof}

\begin{theorem}[Second-order total generalised variation Gaussian denoising \cite{reyes2015a}]
    \label{theorem:tgv2-l2cost-interior}
    %In case of second order total generalised varition denoising of Example \ref{example:tgv2-1},
    Suppose that the data $f, f_0 \in L^2(\Omega) \isect \BVspace(\Omega)$ satisfies for some $\alpha_2>0$ the condition
    \begin{equation}
        \label{eq:l2cost-interior-condition-tgv2}
        \TGV^2_{(\alpha_2, 1)}(f)
        >
        \TGV^2_{(\alpha_2, 1)}(f_0).
    \end{equation}
    Then there exists $\bar \mu,\bar\gamma > 0$ such any optimal solution $\alpha_{\gamma,\mu}=((\alpha_{\gamma,\mu})_1, (\alpha_{\gamma,\mu})_2)$ to the problem
    \[
        \min_{\alpha \in [0, \infty]^2} \frac{1}{2}\norm{f_0-v_\alpha}_{L^2(\Omega)}^2
    \]
    with
    \[
        \begin{split}
        (v_\alpha, w_\alpha) \in \argmin_{\substack{v \in \BVspace(\Omega)\\ w \in \BDspace(\Omega)}} 
            \Bigl(
            &
            \frac{1}{2}\norm{f-v}_{L^2(\Omega)}^2
            + \alpha_1\abs{Dv-w}_\gamma(\Omega)
            + \alpha_2\abs{Ew}_\gamma(\Omega)
            \\
            &
            + \frac{\mu}{2}\norm{(\grad v, \grad w)}_{L^2(\Omega; \R^\DIMdomain \times \R^{\DIMdomain \times \DIMdomain})}^2
            \Bigr)
        \end{split}
    \]
    satisfies $(\alpha_{\gamma,\mu})_1, (\alpha_{\gamma,\mu})_2 > 0$ whenever $\mu \in [0, \bar \mu]$, $\gamma \in [\bar\gamma, \infty]$.
\end{theorem}

Here we recall that $\BDspace(\Omega)$ is the space of vector fields of bounded deformation \cite{temam1985mpp}.
Again, the noisy data has to oscillate more in terms of $\TGV^2$ than the ground-truth does, for the existence of an interior optimal solution to \eqref{eq:learn}. This of course allows us to avoid constraints on $\alpha$.

Observe that we allow for infinite parameters $\alpha$. We do not seek to restrict them to be finite, as this will allow us to decide between $\TGV^2$, $\TV$, and $\TV^2$ regularisation.

\begin{proof}[First steps of proof: modelling in the abstract framework]
    To present $\TGV^2$ in the abstract framework, we take take $X=(\BVspace(\Omega) \isect L^2(\Omega))\times \BDspace(\Omega)$, and $Y=L^2(\Omega)$. We denote $u=(v, w)$, and set
    \[
         K(v, w)=v, \quad A_1u = Dv-w,\quad \text{ and } \quad A_2u =Ew
    \]
    for $E$ the symmetrised differential. 
    For the smoothing term we take
     \[
        H(u)=\frac{1}{2}\norm{(\grad v, \grad w)}_{L^2(\Omega; \R^\DIMdomain \times \R^{\DIMdomain \times \DIMdomain})}^2.
    \]
    For all the gory details we again point the reader to \cite{reyes2015a}.
\end{proof}

We also have a result on the approximation properties of the numerical models as $\gamma \upto \infty$ and $\mu \downto 0$.
%Outer semicontinuity of a set-valued map $S: \R^k \setto \R^\DIMdomain$ means \cite{rockafellar-wets-va} that for any convergent sequence $x^k \to x$ and $S(x^k) \ni y^k \to y$, we have $y \in S(x)$. 
Roughly, the the outer semicontinuity \cite{rockafellar-wets-va} of the solution map $\SolM$ in the next theorem means that as the numerical regularisation vanishes, any optimal parameters for the regularised models \eqref{eq:learn-numerical-single} tend to some optimal parameters of the original model \eqref{eq:learn}.

\begin{theorem}[\cite{reyes2015a}]
   In the setting of Theorem \ref{theorem:tv-l2cost-interior} and Theorem \ref{theorem:tgv2-l2cost-interior}, there exist $\bar\gamma \in (0, \infty)$ and $\bar \mu \in (0, \infty)$ such that the solution map 
    \[
        %\SolM(\gamma, \mu) \defeq \argmin_{\alpha \in \SPACEalphaCompactPos}~\costf(u_{\alpha,\gamma,\mu})
        (\gamma, \mu) \mapsto \alpha_{\gamma,\mu}
    \]
    is outer semicontinuous within $[\bar\gamma,\infty] \times [0, \bar\mu]$.
\end{theorem}

We refer to \cite{reyes2015a} for further, more general results of the type in this section. These include analogous of the above ones for a novel Huberised total variation cost functional.

\subsection{Optimality conditions}
In order to compute optimal solutions to the learning problems, a proper characterization of them is required. Since \eqref{eq:learn-numerical-single} constitute PDE-constrained optimisation problems, suitable techniques from this field may be utilized. For the limit cases, an additional asymptotic analysis needs to be performed in order to get a sharp characterization of the solutions as $\gamma \to \infty$ or $\mu \to 0$, or both.

Several instances of the abstract problem \eqref{eq:learn-numerical-single} have been considered in previous contributions. The case with Total Variation regularization was considered in \cite{de2013image} in presence of several noise models. There the G\^ateaux differentiability of the solution operator was proved, which lead to the derivation of an optimality system. Thereafter an asymptotic analysis with respect to $\gamma \to \infty$ was carried out (with $\mu>0$), obtaining an optimality system for the corresponding problem. In that case the optimisation problem corresponds to one with variational inequality constraints and the characterization concerns C-stationary points.

Differentiability properties of higher order regularisation solution operators were also investigated in \cite{reyes2015b}. A stronger Fr\'echet differentiability result was proved for the $\TGV^2$ case, which also holds for TV. These stronger results open the door, in particular, to further necessary and sufficient optimality conditions. %and to convergence proofs of infinite dimensional solution algorithms.

For the general problem \eqref{eq:learn-numerical-single}, using the Lagrangian formalism the following optimality system is obtained:
\begin{multline} \label{eq: OS state}
\mu \int_\Omega \iprod{\grad u}{\grad v} \d x + \sum_{i=1}^M \int_\Omega \lambda_i \, \phi_i'(Ku)Kv \d x \\ + \sum_{j=1}^N \int_\Omega \alpha_j \iprod{h_\gamma(A_j u)}{A_j v} \d x =0, \quad \forall v \in V,
\end{multline}
\begin{multline} \label{eq: adjoint equation}
\mu \int_\Omega \iprod{\grad p}{\grad v} \d x + \sum_{i=1}^M \int_\Omega \iprod{\lambda_i\phi_i''(Ku)Kp}{Kv} \d x \\ + \sum_{j=1}^N \int_\Omega \alpha_j \iprod{h_\gamma'^*(A_j u) A_j p}{A_j v} \d x  =-F'(u)v, \quad \forall v \in V,
\end{multline}
\begin{equation} \label{eq: OS VI1}
\int_\Omega \phi_i(Ku)Kp (\zeta-\lambda_i) \d x \geq 0, \quad \forall \zeta \geq 0, \ i=1, \dots, M,
\end{equation}
\begin{equation} \label{eq: OS VI2}
\int_\Omega h_\gamma(A_j u)A_j p (\eta-\alpha_j) \d x \geq 0, \quad \forall \eta \geq 0, \ j=1, \dots, N,
\end{equation}
where $V$ stands for the Sobolev space where the regularised image lives (typically a subspace of $H^1(\Omega; \mathbb R^m)$ with suitable homogeneous boundary conditions), $p \in V$ stands for the adjoint state and $h_\gamma$ is a regularized version of the TV subdifferential, for instance,
\begin{equation}\label{eq:local reg. of q}
h_{\gamma}(z):=
\begin{cases}
\frac{z}{|z|} &\text{ if }~\gamma |z|-1 \geq \frac{1}{2\gamma}\\
 \frac{z}{|z|} (1- \frac{\gamma}{2} (1- \gamma |z|+\frac{1}{2\gamma})^2) &\text{ if }~\gamma |z|-1 \in (-\frac{1}{2\gamma}, \frac{1}{2\gamma})\\
\gamma z &\text{ if }~\gamma |z|-1 \leq -\frac{1}{2\gamma}.
\end{cases}
\end{equation}
This optimality system is stated here formally. Its rigorous derivation has to be justified for each specific combination of spaces, regularisers, noise models and cost functionals.

With help of the adjoint equation \eqref{eq: adjoint equation} also gradient formulas for the reduced cost functional $\mathcal F (\lambda,\alpha):= F(u_{\alpha, \lambda}, \lambda, \alpha)$ are derived:
\begin{equation} \label{eq: gradient formulas}
(\nabla_\lambda \mathcal F)_i= \int_\Omega \phi_i(Ku)Kp \d x, \quad  \quad (\nabla_\alpha \mathcal F)_j= \int_\Omega h_\gamma(A_j u) A_j p \d x,
\end{equation}
for $i=1, \dots, M$ and $j=1, \dots, N$, respectively. The gradient information is of numerical importance in the design of solution algorithms. In the case of finite dimensional parameters, thanks to the structure of the minimisers reviewed in Section 2, the corresponding variational inequalities \eqref{eq: OS VI1} and \eqref{eq: OS VI2} turn into equalities. This has important numerical consequences, since in such cases the gradient formulas \eqref{eq: gradient formulas} may be used without additional projection steps. This will be commented in detail in the next section.

\section{Numerical optimisation of the learning problem}\label{sec:algorithms}
%!TEX root = article.tex

\subsection{Adjoint based methods}    \label{sec:adjoint}
The derivative information provided through the adjoint equation \eqref{eq: adjoint equation}
%\begin{multline}
%\mu \int_\Omega D p \ D v + \sum_{i=1}^M \int_\Omega \lambda_i \, \Phi_i''(Ku)Kv\\ + \sum_{j=1}^N \int_\Omega \alpha_j \, h_\gamma'^*(A_j u) A_j p \ A_j v =-F'(u)v, \quad \forall v \in V,
%\end{multline}
may be used in the design of efficient second-order algorithms for solving the bilevel problems under consideration. Two main directions may be considered in this context: Solving the original problem via optimisation methods \cite{calatronidynamic,reyes2015b,ochs_ssvm2015}, and solving the full optimality system of equations \cite{kunisch2013bilevel,chungdelosreyes}. The main advantage of the first one consists in the reduction of the computational cost when a large image database is considered (this issue will be treated in detail below). When that occurs, the optimality system becomes extremely large, making it difficult to solve it in a manageable amount of time. The advantage of the second approach, on the other hand, consists in the possibility of using efficient (possibly generalized) Newton solvers, which have been intensively developed in the last years. 

Let us first describe the quasi-Newton methodology considered in \cite{calatronidynamic,reyes2015b} and further developed in \cite{reyes2015b}. For the design of a quasi-Newton algorithm for the bilevel problem with, e.g., one noise model ($\lambda_1=1$), the cost functional has to be considered in reduced form as $\mathcal F (\alpha):= F(u_\alpha, \alpha),$ where $u_\alpha$ is implicitly determined by solving the denoising problem 
\begin{equation} \label{alg eq: denoising}
u_\alpha = \arg \min_{u \in V} \ \frac{\mu}{2} \int_\Omega \norm{\grad u}^2 \d x +\sum_{j=1}^N \int_\Omega \alpha_j \d |A_j u|_\gamma +\int_\Omega \phi(u) \d x, \quad \mu >0.
\end{equation}
%where $X^H:=X \isect \Dom \mu\Smoother$. 
Using the gradient formula for $\mathcal F$, 
\begin{equation}
(\nabla \mathcal F (\alpha^{(k)}))_j=\int_\Omega h_\gamma(A_j u) A_j p \d x, \quad j=1, \dots,N,
\end{equation}
the BFGS matrix may be updated with the classical scheme
\begin{equation}  \label{alg def:BFGS update}
B_{k+1}=B_k- \frac{B_k s_k \otimes B_k s_k}{(B_k s_k,s_k)}+\frac{z_k \otimes z_k}{(z_k,s_k)},\end{equation}
where $s_k=\alpha^{(k+1)}-\alpha^{(k)}$, $z_k= \nabla \mathcal F (\alpha^{(k+1)})-\nabla \mathcal F (\alpha^{(k)})$ and $(w \otimes v)\varphi:=(v,\varphi)w$.
For the line search strategy, a backtracking rule may be considered,  with the classical Armijo criteria 
\begin{equation}   \label{alg def:Armijo}
\mathcal F (\alpha^{(k)}+ t_k d^{(k)} )- \mathcal F (\alpha^{(k)}) \leq t_k \beta \nabla  \mathcal F (\alpha^{(k)})^T d^{(k)}, \quad \beta \in (0,1],
\end{equation}
where $d^{(k)}$ stands for the quasi-Newton descent direction and $t_k$ the length of the quasi-Newton step.
We consider, in addition, a cyclic update based on curvature verification, i.e., we update the quasi-Newton matrix only if the curvature condition $(z_k,s_k)>0$ is satisfied. The positivity of the parameter values is usually preserved along the iterations, making a projection step superfluous in practice. In more involved problems, like the ones with $\TGV^2$ or $\text{ICTV}$ denoising, an extra criteria may be added to the Armijo rule, guaranteeing the positivity of the parameters in each iteration. Experiments with other line search rules (like Wolfe) have also been performed. Although these line search strategies automatically guarantee the satisfaction of the curvature condition (see, e.g., \cite{nocedal2006numerical}), the interval where the parameter $t_k$ has to be chosen appears to be quite small and is typically missing.

The denoising problems \eqref{alg eq: denoising} may be solved either by efficient first- or second-order methods. In previous works we considered primal-dual Newton type
algorithms (either classical or semismooth) for this purpose. Specifically, by introducing the dual variables $q_i, ~i=1, \dots, N$, a necessary and sufficient condition for the lower level is given by 
\begin{align} 
&\mu \int_\Omega \iprod{\grad u}{\grad v} \d x +\sum_{i=1}^N  \int_\Omega \iprod{q_i}{A_i v} \d x +\int_\Omega \iprod{\phi'(u)}{v} \d x =0,  && \forall v \in V, \label{eq:primal-dual1}\\
&q_i=\alpha_i \, h_\gamma(A_i u) \quad \text{a.e. in }\Omega, ~ i=1, \dots,N, \label{eq:primal-dual2}
\end{align}
where $h_\gamma(z):= \frac{z}{\max(1/\gamma, |z|)}$ is a regularized version of the TV subdifferential, and the generalized Newton step has the following Jacobi matrix
\begin{equation}
    \begin{pmatrix}
        L+\phi''(u) & A_1^* & \ldots & A_\NA^* \\
        - \alpha_1 \left[ \nO(A_1 u) -\chi_1 \frac{A_1 u \otimes A_1 u}{|A_1 u|^3} \right] A_1 & I & 0 & 0 \\
        \vdots & 0 & \ddots & 0 \\
        - \alpha_\NA \left[ \nO(A_\NA u) -\chi_\NA \frac{A_\NA u \otimes A_\NA u}{|A_\NA u|^3} \right] A_\NA & 0 & 0 & I
    \end{pmatrix}
\end{equation}
where $L$ is an elliptic operator, $\chi_i(x)$ is the indicator function of the set $\{x: \gamma |A_i u| > 1 \}$ and $\nO(A_i u):= \frac{\min(1,\gamma |A_i u|)}{|A_i u|}$, for $i=1, \dots, \NA$.
In practice, the convergence neighbourhood of the classical method is too small and some sort of globalization is required. Following \cite{hintermuller2006infeasible} a modification of the matrix was systematically considered, where the term $\frac{A_i u \otimes A_i u}{|A_i u|^3}$ is replaced by $\frac{q_i}{\max(|q_i|,\alpha_i)} \otimes \frac{A_i u}{|A_i u|^2}$. The resulting algorithm exhibits both a global and a local superlinear convergent behaviour.

For the coupled BFGS algorithm a warm start of the denoising Newton methods was considered, using the image computed in the previous quasi-Newton iteration as initialization for the lower level problem algorithm. The adjoint equations, used for the evaluation of the gradient of the reduced cost functional, are solved by means of sparse linear solvers.

Alternatively, as mentioned previously, the optimality system may be solved at once using nonlinear solvers. In this case the solution is only a stationary point, which has to be verified a-posteriori to be a minimum of the cost functional. This approach has been considered in \cite{kunisch2013bilevel} and \cite{chungdelosreyes} for the finite- and infinite-dimensional cases, respectively.
The solution of the optimality system also presents some challenges due to the nonsmoothness of the regularisers and the positivity constraints.

For simplicity, consider the bilevel learning problem with the TV-seminorm, a single Gaussian noise model and a scalar weight $\alpha$. The optimality system for the problems reads as follows
\begin{subequations} \label{alg eq: optimality system reg}
 \begin{gather} \label{alg eq:optimality condition state equation}
  \mu \int_\Omega \iprod{\grad u}{\grad v} \d x + \int_\Omega \alpha h_\gamma(\grad u) \grad v \d x + \int_\Omega (u-f)v \d x =0, \forall v \in V,
\\
\label{alg eq:optimality condition adjoint equation}
\begin{split}
\mu \int_\Omega \iprod{\grad p}{\grad v} \d x + \int_\Omega \alpha \iprod{h_\gamma'^*(\grad u) \grad p}{\grad v} \d x 
+ & \int_\Omega  p\,v \d x \\ & =-F'(u)v, 
 \quad \forall v \in V,
 \end{split}
\\
 \sigma=  \int_\Omega \iprod{h_\gamma(\grad u)}{\grad p} \d x.
\\
\label{alg eq:optimality condition complementarity system}
 \sigma \geq 0, ~ \alpha \geq 0, ~\sigma \cdot \alpha =0.
\end{gather}
\end{subequations}
where $h_\gamma$ is given by, e.g., equation \eqref{eq:local reg. of q}.
%a $C^{1,1}-$regularization of the TV subdifferential, for instance,
%\begin{equation}\label{eq:local reg. of q}
%h_{\gamma}(z):=
%\begin{cases}
%\frac{z}{|z|} &\text{ if }~\gamma |z|-1 \geq \frac{1}{2\gamma}\\
% \frac{z}{|z|} (1- \frac{\gamma}{2} (1- \gamma |z|+\frac{1}{2\gamma})^2) &\text{ if }~\gamma |z|-1 \in (-\frac{1}{2\gamma}, \frac{1}{2\gamma})\\
%\gamma z &\text{ if }~\gamma |z|-1 \leq -\frac{1}{2\gamma}.
%\end{cases}
%\end{equation}
The Newton iteration matrix for this coupled system has the following form:
    \begin{equation*}
    \begin{pmatrix}
    L + \,\grad^* \alpha^{(k)} h'_\gamma(\grad u^k)\grad &0 &  \grad^* h_\gamma(\grad u^k)\\
    \grad^* \alpha^{(k)} h_\gamma''(\grad u^k) \, \grad p \,\grad+F''(u^k) &L + \grad^* \alpha^{(k)} h'_\gamma(\grad u^k)\grad &\grad^*h'_\gamma(\grad u^k) \grad p\\
    h'_\gamma(\grad u^k)\grad p \grad &h_\gamma(\grad u^k) \grad &0
    \end{pmatrix}.
    \end{equation*}
The structure of this matrix leads to similar difficulties as for the denoising Newton iterations described above. To fix this and get good convergence properties, Kunisch and Pock \cite{kunisch2013bilevel} proposed an additional feasibility step, where the iterates are projected on the nonlinear constraining manifold. In \cite{chungdelosreyes}, similarly as for the lower level problem treatment, modified Jacobi matrices are built by replacing the terms $h'_\gamma(u_k)$ in the diagonal, using projections of the dual multipliers. Both approaches lead to globally convergent algorithm with locally superlinear convergence rates. Also domain decomposition techniques were tested in \cite{chungdelosreyes} for the efficient numerical solution of the problem. 

By using this optimize-then-discretize framework, resolution independent solution algorithms may be obtained. Once the iteration steps are well specified, both strategies outlined above use a suitable discretization of the image. Typically a finite differences scheme with mesh size step $h>0$ is used for this purpose. The minimum possible value of $h$ is related to the resolution of the image. For the discretization of the Laplace operator the usual five point stencil is used, while forward and backward finite differences are considered for the discretization of the divergence and gradient operators, respectively. Alternative discretization methods (finite elements, finite volumes, etc) may be considered as well, with the corresponding operators. %The resulting discretized denoising problem is then given by ..., where... stands for the discrete Laplace, gradient and divergence operators, respectively. 

\subsection{Dynamic sampling}

For a robust and realistic learning of the optimal parameters, ideally, a rich database of $K$ images, $K\gg 1$ should be considered (like, for instance, MRI applications, compare Section \ref{sec:numerics single}). Numerically, this consists in solving a large set of nonsmooth PDE-constraints of the form \eqref{eq:primal-dual1}- \eqref{eq:primal-dual2} in each iteration of the BFGS optimisation algorithm \eqref{alg def:BFGS update}.

In \cite{calatronidynamic} we extended to our imaging framework a \emph{dynamic sample size} stochastic approximation method proposed by Byrd et al. \cite{byrd2012sample}. 
%The main idea of this method is to consider an initial, small, training sample of the dictionary to start the algorithm with and \emph{dynamically} increasing its size, if needed, throughout the different steps of the optimisation process. 
The algorithm starts by selecting from the whole dataset a sample $S$ whose size $|S|$ is small compared to the original size $K$. In the following iterations, if the approximation of the optimal parameters computed produces an improvement in the cost functional, then the sample size is kept unchanged and the optimisation process continues selecting in the next iteration a new sample of the same size. Otherwise, if the approximation computed is not a good one, a new, larger, sample size is selected and a new sample $S$ of this new size is used to compute the new step. The key point in this procedure is clearly the rule that checks throughout the progression of the algorithm, whether the approximation we are performing is good enough, i.e. the sample size is big enough, or has to be increased. Because of this systematic check, such sampling strategy is called \emph{dynamic}. Denoting by $u_\alpha^k$ the solution of \eqref{eq:primal-dual1}-\eqref{eq:primal-dual2} and by $f_0^k$ the ground-truth images for every $k=1,\dots,K$,  we consider now the reduced cost functional $\mathcal{F}(\alpha)$ in correspondence of the whole database
\begin{equation*}   
\mathcal{F}(\alpha)=\frac{1}{2K}\sum_{k=1}^{K} \| u^k_\alpha-f_0^k \|^2_{L^2},
\end{equation*}
we consider, for every sample $S\subset\left\{1,\ldots,K\right\}$, the batch objective function:
\begin{equation*}
\mathcal{F}_S(\alpha):=\frac{1}{2|S|}\sum_{k\in S}\|u^k_\alpha-f_0^k\|^2_{L^2}.
\end{equation*}

As in \cite{byrd2012sample}, we formulate in \cite{calatronidynamic} a condition on the batch gradient $\nabla \mathcal{F}_S$ which imposes in every stage of the optimisation that the direction $-\nabla \mathcal{F}_S$ is a descent direction for $\mathcal{F}$ at $\alpha$ if the following condition holds:
\begin{equation} \label{descentcond}
\|\nabla \mathcal{F}_S(\alpha)-\nabla \mathcal{F}(\alpha)\|_{L^2}
\leq \theta\|\nabla \mathcal{F}_S(\alpha)\|_{L^2},\quad \theta\in[0,1).
\end{equation}

The computation of $\nabla \mathcal{F}$ may be very expensive for applications involving large databases and nonlinear constraints, so we rewrite \eqref{descentcond} as an estimate of the variance of the random vector $\nabla \mathcal{F}_S(\alpha)$. We do not report here the details of the derivation of such estimate, but we refer the interested reader to \cite[Section 2]{calatronidynamic}. Here, we just underline that through such a condition on the variance one can control in each iteration of the BFGS optimisation whether the sampling approximation is accurate enough and, if this is not the case, a new larger sample size may be determined in order to reach the desired level of accuracy, depending on the parameter $\theta$ in \eqref{descentcond}. 

\begin{algorithm}[t!]
\caption{Dynamic Sampling BFGS}
\label{alg:dynbfgs}
\label{alg:dynsamp}
\begin{algorithmic}[1]
\STATE\small Initialize: $\alpha_0$, sample $\mathcal{S}_0$ with $|S_0|\ll K$ and model parameter $\theta$, $k=0$.
\STATE \textbf{while} {BFGS not converging, $k\geq0$}
\STATE{\quad  sample $|S_k|$ PDE constraints to solve};
\STATE{\quad update the BFGS matrix};
\STATE{\quad compute direction $d_k$ by BFGS and steplength $t_k$ by Armijo cond. \eqref{alg def:Armijo};}
\STATE{\quad define new iterate: $\alpha_{k+1}=\alpha_{k}+t_k d_k$;}
\STATE{\quad \textbf{if variance condition} is satisfied then}
\STATE{\quad\quad maintain the sample size: $|S_{k+1}|=|S_k|$;}
\STATE{\quad\textbf{else} augment $S_k$ such that condition \textbf{variance condition} is verified.}
\STATE{\textbf{end}}
\end{algorithmic}
\end{algorithm}

In order to improve upon the traditional slow convergence drawback of steepest descent methods, we combined the Dynamic Sampling strategy described above with BFGS method \eqref{alg def:BFGS update}, as described in Algorithm \ref{alg:dynbfgs}.

\section{Learning the image model}\label{sec:image}
%!TEX root = article.tex

One of the main aspects of discussion in the modelling of variational image reconstruction is the type and strength of regularisation that should be imposed on the image. That is, what is the correct choice of regularity that should be imposed on an image and how much smoothing is needed in order to counteract imperfections in the data such as noise, blur or undersampling. In our variational reconstruction approach \eqref{reg1:eq} this boils down to the question of choosing the regulariser $R(u)$ for the image function $u$ and the regularisation parameter $\alpha$. In this section we will demonstrate how functional modelling and data learning can be combined to derive optimal regularisation models. To do so, we focus on {\bf T}otal {\bf V}ariation (TV) type regularisation approaches and their optimal setup. The following discussion constitutes the essence of our derivations in \cite{reyes2015b}, including an extended numerical discussion with an interesting application of our approach to cartoon-texture decomposition.

\subsection{Total variation type regularisation}

\begin{figure}[t]  \label{fig:effect}
\centering
% Do not use subfigure, use the more modern subcaption!!!
\newcommand{\SUBFIG}[3]{\begin{subfigure}[t]{#2}\includegraphics[width=\textwidth]{#3}\caption{#1}\end{subfigure}}
\SUBFIG{original}{0.24\textwidth}{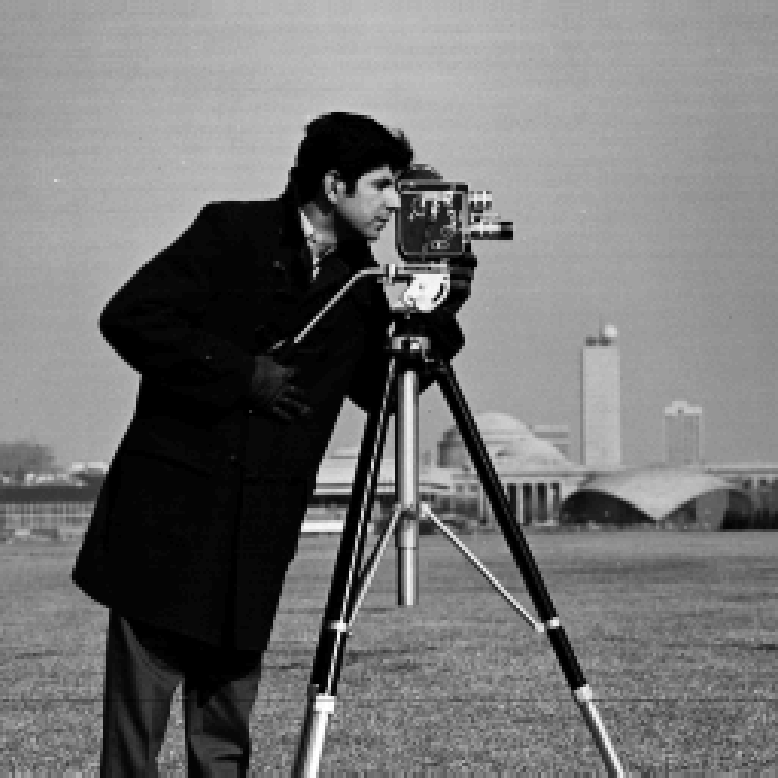}\,%
\SUBFIG{noisy}{0.24\textwidth}{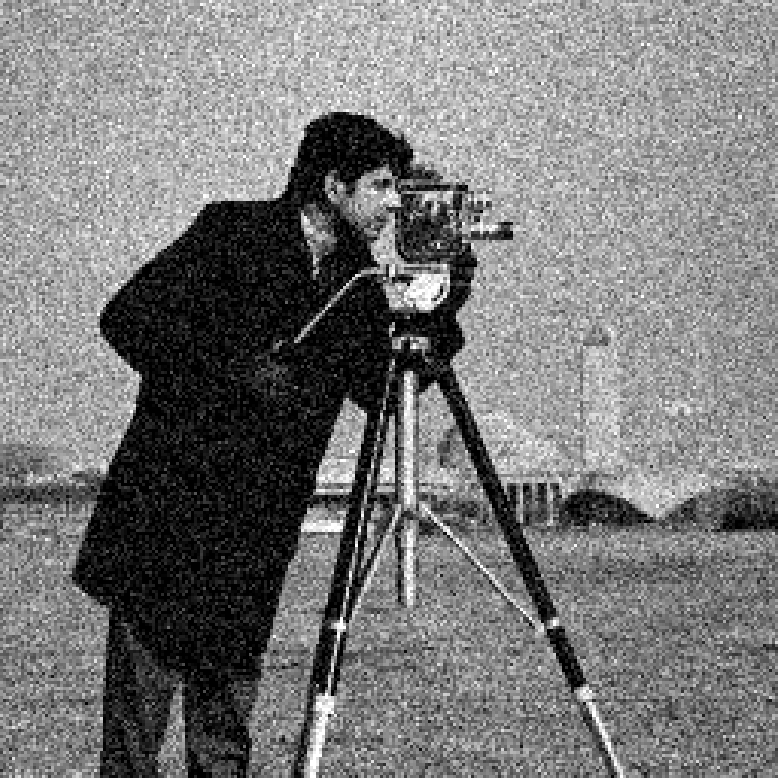}\,%
\SUBFIG{$R(u)=\|\nabla u\|_2^2$}{0.24\textwidth}{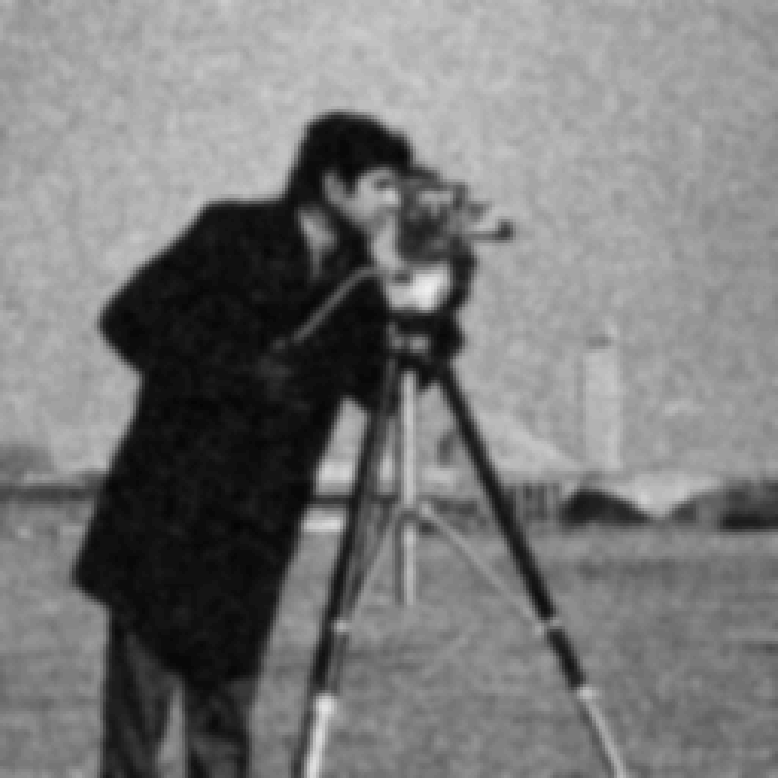}\,%
\SUBFIG{$R(u)=|Du|(\Omega)$}{0.24\textwidth}{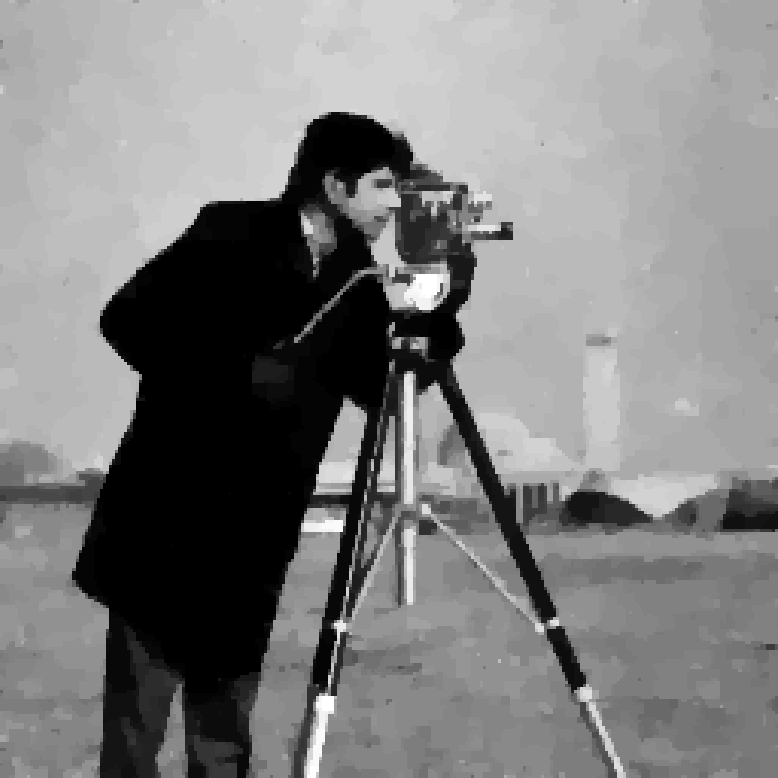}
\caption{The effect of the choice of regularisation in \eqref{reg1:eq}: Choosing the $L^2$ norm squared of the gradient of $u$ as a regulariser imposes isotropic smoothing on the image and smoothes the noise equally as blurring the edges. Choosing the total variation (TV) as a regulariser we are able to eliminate the noise while preserving the main edges in the image.}
\end{figure}

The TV is the total variation measure of the distributional derivative of $u$ \cite{ambrosio2000fbv}, that is for $u$ defined on $\Omega$
\begin{equation}
TV(u)=|Du|(\Omega) = \int_\Omega \d \abs{D u}.
\label{eq:TVreg}
\end{equation}
As the seminal work of Rudin, Osher and Fatemi \cite{rudin1992nonlinear} and many more contributions in the image processing community have proven, a non-smooth first-order regularisation procedure as TV results in a nonlinear smoothing of the image, smoothing more in homogeneous areas of the image domain and preserving characteristic structures such as edges, compare Figure \ref{fig:effect}. More precisely, when TV is chosen as a regulariser in \eqref{reg1:eq} the reconstructed image is a function in $BV$ the space of functions of bounded variation, allowing the image to be discontinuous as its derivative is defined in the distributional sense only. Since edges are discontinuities in the image function they can be represented by a BV regular image. In particular, the TV regulariser is tuned towards the preservation of edges and performs very well if the reconstructed image is piecewise constant. 

\begin{figure}[t]
\begin{center}
\includegraphics[height=2.5cm]{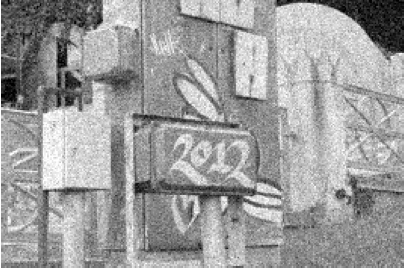}\includegraphics[height=2.5cm]{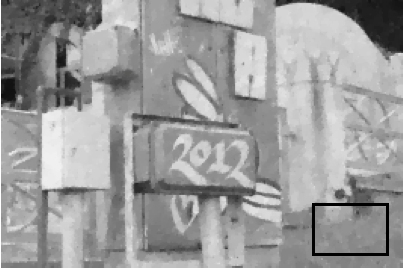} \includegraphics[height=2.5cm]{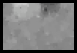}
%\begin{picture}(100,0)
%\thicklines\put(105,35){\vector(1,1){20}}
%\end{picture}
\caption{TV image denoising and the staircasing effect: (l.) noisy image, (m.) denoised image, (r.) detail of the bottom right hand corner of the denoised image to visualise the staircasing effect (the creation of blocky-like patterns due to the first-order regulariser).}
\label{fig:staircasing}
\end{center}
\end{figure}

Because one of the main characteristics of images are edges as they define divisions between objects in a scene, the preservation of edges seems like a very good idea and a favourable feature of TV regularisation. The drawback of such a regularisation procedure becomes apparent as soon as images or signals (in 1D) are considered which do not only consist of constant regions and jumps, but also possess more complicated, higher-order structures, e.g. piecewise linear parts. The artefact introduced by TV regularisation in this case is called staircasing \cite{ring2000structural}, compare Figure \ref{fig:staircasing}. 

One possibility to counteract such artefacts is the introduction of higher-order derivatives in the image regularisation. Here, we mainly concentrate on two second-order total variation models: the recently proposed {\bf T}otal {\bf G}eneralized {\bf V}ariation (TGV)  \cite{bredies2011tgv} and the {\bf I}nfimal-{\bf C}onvolution {\bf T}otal {\bf V}ariation (ICTV) model of Chambolle and Lions \cite{chambolle97image}. We focus on second-order TV regularisation only since this is the one which seems to be most relevant in imaging applications \cite{knoll2010second,bredies2012total}. For $\Omega\subset\mathbb R^2$ open and bounded, the ICTV regulariser reads
\begin{equation}\label{eq:ictv}
\ICTV_{\alpha,\beta}(u) \defeq \min_{v\in W^{1,1}(\Omega),~ \nabla v\in BV(\Omega)} \alpha \norm{Du-\nabla v}_{\Meas(\Omega; \R^2)} + \beta \norm{D\nabla v}_{\Meas(\Omega; \R^{2\times 2})}.
\end{equation}
On the other hand, second-order TGV \cite{sampta2011tgv,bredies2013properties} reads
\begin{equation}
    \label{eq:tgv}
    \TGV^2_{\alpha,\beta}(u) 
    \defeq
    \min_{w \in BD(\Omega)}
    \alpha \norm{Du-w}_{\Meas(\Omega; \R^2)}
    +
    \beta \norm{Ew}_{\Meas(\Omega; \Sym^2(\R^2))}.
\end{equation}
Here $
    \BDspace(\Omega)
    \defeq
    \{ w \in L^1(\Omega; \R^\DIMdomain)
        \mid
        \norm{Ew}_{\Meas(\Omega; \R^{\DIMdomain \times \DIMdomain})}
     < \infty \}
$
is the space of vector fields of bounded deformation on $\Omega$, $E$ denotes the \emph{symmetrised gradient} and $\mathrm{Sym}^2(\mathbb{R}^2)$ the space of symmetric tensors of order $2$ with arguments in $\mathbb{R}^2$. The parameters $\alpha,\beta$ are fixed positive parameters. The main difference between \eqref{eq:ictv} and \eqref{eq:tgv} is that we do not generally have that $w=\nabla v$ for any function $v$. That results in some qualitative differences of ICTV and TGV regularisation, compare for instance \cite{benning2011higher}. Substituting $\alpha R(u)$ in \eqref{reg1:eq} by $\alpha TV(u)$, $\TGV^2_{\alpha,\beta}(u)$ or $\ICTV_{\alpha,\beta}(u)$ gives the TV image reconstruction model, TGV image reconstruction model and the ICTV image reconstruction model, respectively.

\subsection{Optimal parameter choice for TV type regularisation}

\begingroup
    \def\SQl{128}
    \def\SQb{160}
    \def\SQr{192}
    \def\SQt{224}
    \renewcommand{\igraph}[2]{\includegraphics[#1,bb=128 160 192 224,clip]{{../tgv-learn-resimg/#2}.png}}

    \begin{figure}[t]
        \centering
        \begin{subfigure}[t]{0.3\textwidth}%
            \resplotzxx{intro-low}
            \caption{Too low $\beta$ / High oscillation}
        \end{subfigure}
        \begin{subfigure}[t]{0.3\textwidth}%
            \resplotzxx[\drawzoomarea]{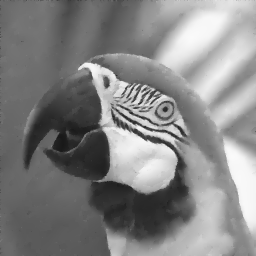}
            \caption{Optimal $\beta$}
        \end{subfigure}
        \begin{subfigure}[t]{0.3\textwidth}%
            \resplotzxx{intro-high}
            \caption{Too high $\beta$ / almost $\TV$}
        \end{subfigure}
        \caption{Effect of $\beta$ on $\TGV^2$ denoising with optimal $\alpha$}
        \label{fig:intro}
    \end{figure}

    \begin{figure}[t]
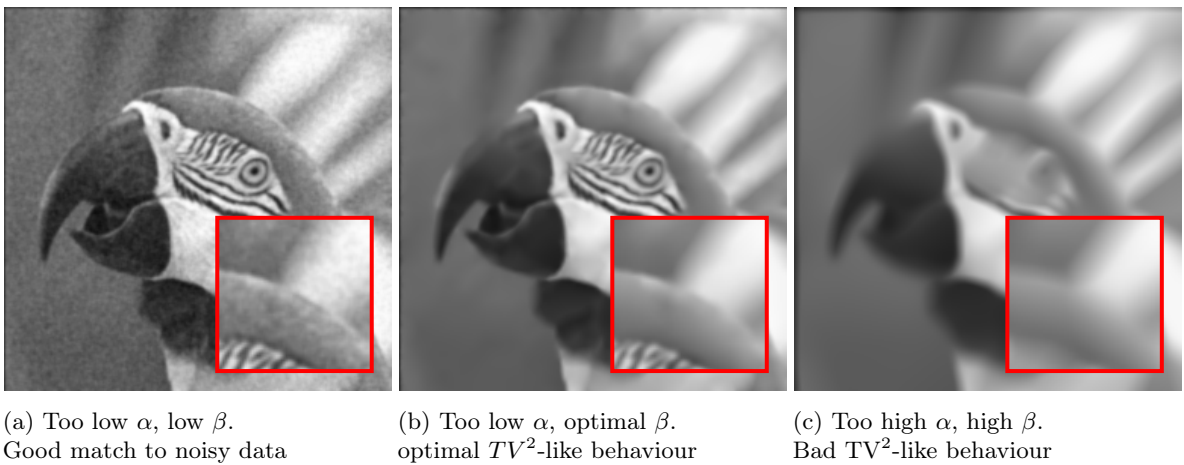

        \centering
        \begin{subfigure}[t]{0.3\textwidth}%
            \resplotzxx{intro-alphahigh-betalow}
            \caption{Too low $\alpha$, low $\beta$. \\Good match to noisy data}
        \end{subfigure}
        \begin{subfigure}[t]{0.3\textwidth}%
            \resplotzxx{intro-alphahigh-betaopt}
            \caption{Too low $\alpha$, optimal $\beta$. \\ optimal $TV^2$-like behaviour}
        \end{subfigure}
        \begin{subfigure}[t]{0.3\textwidth}%
            \resplotzxx{intro-alphahigh-betahigh}
            \caption{Too high $\alpha$, high $\beta$.\\ Bad $\TV^2$-like behaviour}
        \end{subfigure}
        \caption{Effect of choosing $\alpha$ too large in $\TGV^2$ denoising}
        \label{fig:intro2}
    \end{figure}
\endgroup

The regularisation effect of TV and second-order TV approaches as discussed above heavily depends on the choice of the regularisation parameters $\alpha$ (i.e. $(\alpha,\beta)$ for second-order TV approaches). In Figures \ref{fig:intro} and \ref{fig:intro2} we show the effect of different choices of $\alpha$ and $\beta$ in $\TGV^2$ denoising. In what follows we show some results from \cite{reyes2015b} applying the learning approach \eqref{eq:learn-numerical-single} to find optimal parameters in TV type reconstruction models, as well as a new application of bilevel learning to optimal cartoon-texture decomposition. 

\begin{SCfigure}
\centering
        \includegraphics[width=0.6\textwidth]{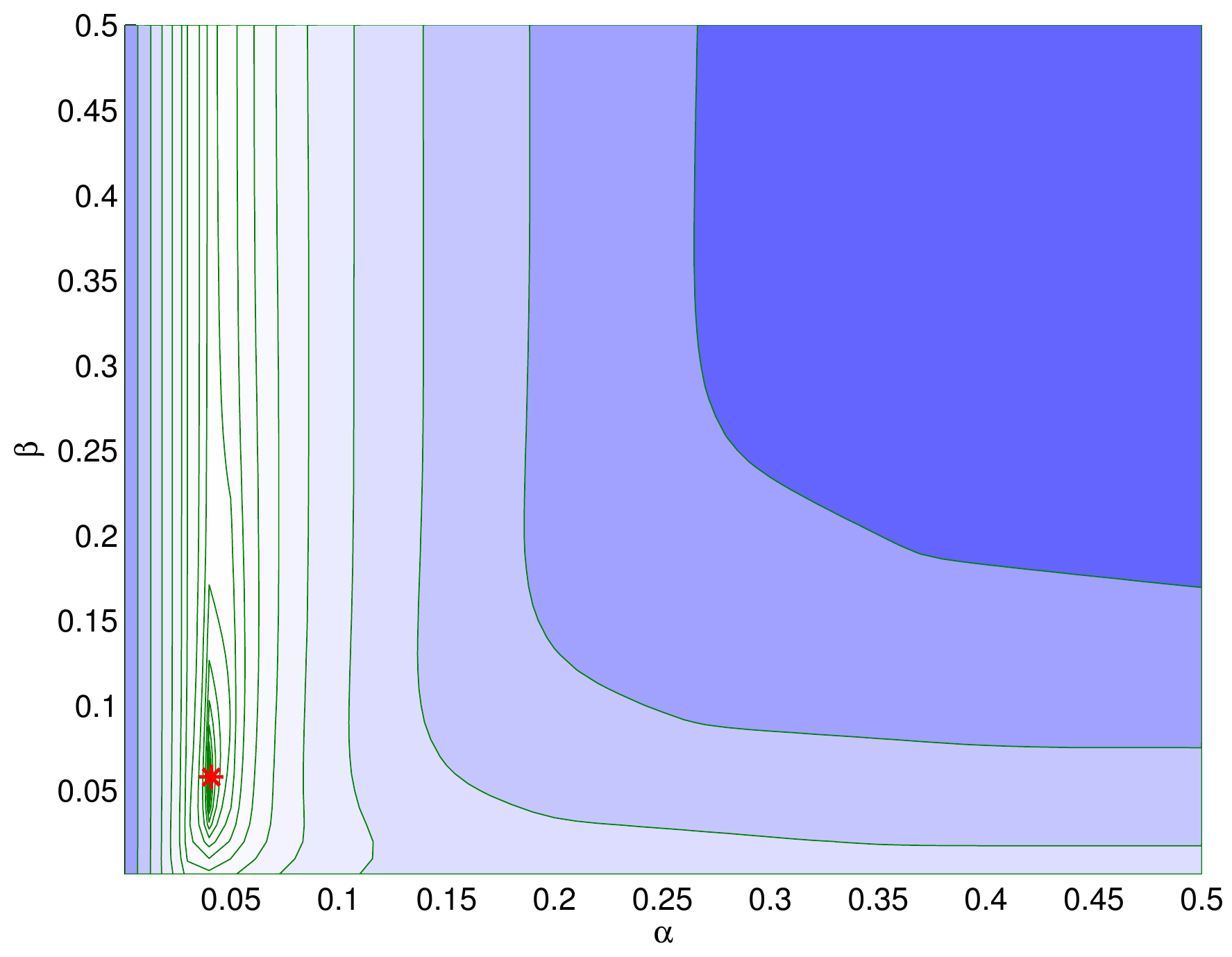}
               \caption{Cost functional value for the {\costltwo} cost functional plotted versus $(\alpha, \beta)$ for $\TGV^2$ denoising.
        The illustration is a contour plot of function value versus $(\alpha, \beta)$.
        }
    \label{fig:landscape-tgv2}
\end{SCfigure}

\paragraph{Optimal TV, $\TGV^2$ and $ICTV$ denoising}
We focus on the special case of $K=Id$ and $L^2$-squared cost $F$ and fidelity term $\Phi$ as introduced in Section \ref{sec:l2theory}. In \cite{reyes2015a,reyes2015b} we also discuss the analysis and the effect of Huber regularised $L^1$ costs, but this is beyond the scope of this paper and we refer the reader to the respective papers. We consider the problem for finding optimal parameters $(\alpha,\beta)$ for the variational regularisation model  
$$
u_{(\alpha,\beta)} \in \argmin_{u \in X} R_{(\alpha,\beta)}(u) + \|u-f\|_{L^2(\Omega)}^2,
$$
where $f$ is the noisy image, $R_{(\alpha,\beta)}$ is either TV in \eqref{eq:TVreg} multiplied by $\alpha$ (then $\beta$ is obsolete), $\TGV^2_{(\alpha,\beta)}$ in \eqref{eq:tgv} or $ICTV_{(\alpha,\beta)}$ in \eqref{eq:ictv}. We employ the framework of \eqref{eq:learn-numerical-single} with a training pair $(f_0,f)$ of original image $f_0$ and noisy image $f$, using $L^2$-squared cost $\costf_{\costltwo}(v) \defeq \frac{1}{2}\norm{f_0 - v}_{\SPACEkuLtwo}^2$. As a first example we consider a photograph of a parrot to which we add Gaussian noise such that the PSNR of the parrot image is $24.72$. In Figure \ref{fig:landscape-tgv2}, we plot by the red star the
discovered regularisation parameter $(\alpha^*, \beta^*)$ reported in
Figure \ref{fig:res-dataset2}. Studying the location of the red star, we
may conclude that the algorithm managed
to find a nearly optimal parameter in very few BFGS 
iterations, compare Table \ref{table:res-dataset2}.

\begin{figure}[t]
    \setlength{\imw}{0.32\textwidth}
    \centering
    \def\subfigprefix{fig:res-dataset2}
%    \inplot{kodim23gray-crop}
%            {Original image}
    \inplot{kodim23gray-crop-noisy}
            {Noisy image}%
    \resplot{tgv2-l2-dataset2-linesearch1-bfgs1-beta0-uinitreg_0.100000_}
            {tgv2}{l2}{$\TVINIT$}%\\
    \resplot{ictv-l2-dataset2-linesearch1-bfgs1-beta0-uinitreg_0.100000_}
            {ictv}{l2}{$\TVINIT$}\\
    \hspace*{.07\textwidth}\parbox{0.5\textwidth}{
            \caption{Optimal denoising results for initial guess $\alphavec=\TVINIT$ 
                     for $\TGV^2$ and $ICTV$, and $\alphavec=0.1/\pixels$ for $\TV$}
            \vspace{6ex}
    \label{\subfigprefix}
    }\hspace*{.0875\textwidth}\ %
    \resplot{tv-l2-dataset2-linesearch1-bfgs1-beta0-uinitreg_0.100000_}
            {tv}{l2}{$0.1/\pixels$}
\end{figure}

\begin{table}
    \caption{Quantified results for the parrot image ($\pixels=256=\text{image width/height in pixels}$)}   
    \label{table:res-dataset2}
    \centering\footnotesize
    %\begin{adjustbox}{max width=\textwidth}
    \begin{tabular}{lll|lllll|l}
    Denoise & 
    Cost & 
    Initial $(\alpha, \beta)$ & 
    Result $(\alpha^*, \beta^*)$ & 
    Cost & 
    SSIM & 
    PSNR &
    Its. &
    Fig. \\
    \hline
    \restabll{tgv2-l2-dataset2-linesearch1-bfgs1-beta0-uinitreg_0.100000_}
            {tgv2}{l2}{$\TVINIT$}{fig:res-dataset2}
    \restabll{ictv-l2-dataset2-linesearch1-bfgs1-beta0-uinitreg_0.100000_}
            {ictv}{l2}{$\TVINIT$}{fig:res-dataset2}
    %\restabll{tgv2-symbregman-dataset2-linesearch1-bfgs1-beta0-uinitreg(0.100000)}
    %        {tgv2}{symbregman}{$\TVINIT$}{}%fig:res-dataset2}
     \restabll{tv-l2-dataset2-linesearch1-bfgs1-beta0-uinitreg_0.100000_}
            {tv}{l2}{$0.1/\pixels$}{fig:res-dataset2}
    %\restabll{tv-symbregman-dataset2-linesearch1-bfgs1-beta0-uinitreg(0.100000)}
    %        {tv}{symbregman}{$0.1/\pixels$}{}%fig:res-dataset2}
    \end{tabular}
    %\end{adjustbox}
\end{table}

\paragraph{Optimizing cartoon-texture decomposition using a sketch}
It is not possible to distinguish noise from texture by the $G$-norm and related approaches \cite{meyer2002oscillating}. Therefore, learning an optimal cartoon-texture decomposition based on a noise image and a ground-truth image is not feasible. What we did instead, is to make a hand-drawn sketch as our expected ``cartoon'' $f_0$, and then use the bi-level framework to find the true ``cartoon'' and ``texture'' as split by the model
\[
    J(u, v; \alpha) = \frac{1}{2}\norm{f-u-v}^2 + \alpha_1 \norm{v}_{\KR} + \alpha_2\TV(u)
\]
for the Kantorovich-Rubinstein norm of \cite{krtv}.
For comparison we also include basic TV regularisation results, where we define $v=f-u$. The results for two different iages are in Figure \ref{fig:res-dataset8-cartoon} and Table \ref{table:res-dataset8-cartoon}, and Figure \ref{fig:res-dataset9-cartoon} and Table \ref{table:res-dataset9-cartoon}, respectively.

%\graphicspath{{../../tgv_learn/resimg}}
\begin{figure}
    \centering
    \def\subfigprefix{fig:res-dataset8-cartoon}
    \inplot{kodim23gray-crop}
            {Original image}
    \inplot{kodim23gray-crop-sketch}
            {Cartoon sketch}
    \resplot{krtv2sketch-l2-dataset8-linesearch1-bfgs1-beta0-uinitreg_0.100000_}
            {krtv2}{l2}{$\TVINITKRTV$}
    \resplot{tv-l2-dataset8-linesearch1-bfgs1-beta0-uinitreg_0.100000_}
            {tv}{l2}{$0.1/\pixels$}
    \resplottexture{krtv2sketch-l2-dataset8-linesearch1-bfgs1-beta0-uinitreg_0.100000_-diff}
            {krtv2}{l2}{$\TVINITKRTV$}
    \resplottexture{tv-l2-dataset8-linesearch1-bfgs1-beta0-uinitreg_0.100000_-diff}
            {tv}{l2}{$0.1/\pixels$}
    \caption{Optimal sketch-based cartoonification for initial guess $\alphavec=\TVINITKRTV$ 
             for $\KRTV$ and $\alphavec=0.1/\pixels$ for $\TV$}
    \label{\subfigprefix}
\end{figure}

\begin{table}
    \caption{Quantified results for cartoon-texture decomposition of the parrot image ($\pixels=256=\text{image width/height in pixels}$)}
    
    \label{table:res-dataset8-cartoon}
    \centering\small
     %\begin{adjustbox}{max width=\textwidth}
    \begin{tabular}{lll|lllll|l}
    Denoise & 
    Cost & 
    Initial $\alphavec$ & 
    Result $\alphavec^*$ & 
    Value & 
    SSIM & 
    PSNR &
    Its. &
    Fig. \\
    \hline
    \restabll{krtv2sketch-l2-dataset8-linesearch1-bfgs1-beta0-uinitreg_0.100000_}
            {krtv2}{l2}{$\TVINITKRTV$}{fig:res-dataset8-cartoon}
    \restabll{tv-l2-dataset8-linesearch1-bfgs1-beta0-uinitreg_0.100000_}
            {tv}{l2}{$0.1/\pixels$}{fig:res-dataset8-cartoon}
    \end{tabular}
    %\end{adjustbox}
\end{table}

\begin{figure}
    \centering
    \def\subfigprefix{fig:res-dataset9-cartoon}
    \inplot{barbara}
            {Original image}
    \inplot{barbara-sketch}
            {Cartoon sketch}
    \resplot{krtv2sketch-l2-dataset9-linesearch1-bfgs1-beta0-uinitreg_0.100000_}
            {krtv2}{l2}{$\TVINIT$}
    \resplot{tv-l2-dataset9-linesearch1-bfgs1-beta0-uinitreg_0.100000_}
            {tv}{l2}{$0.1/\pixels$}
    \resplottexture{krtv2sketch-l2-dataset9-linesearch1-bfgs1-beta0-uinitreg_0.100000_-diff}
            {krtv2}{l2}{$\TVINIT$}
    \resplottexture{tv-l2-dataset9-linesearch1-bfgs1-beta0-uinitreg_0.100000_-diff}
            {tv}{l2}{$0.1/\pixels$}
    \caption{Optimal sketch-based cartoonification for initial guess $\alphavec=\TVINIT$ 
             for $\KRTV$ and $\alphavec=0.1/\pixels$ for $\TV$}
    \label{\subfigprefix}
\end{figure}

\begin{table}
    \caption{Quantified results for cartoon-texture decomposition of the Barbara image ($\pixels=256=\text{image width/height in pixels}$)}
    
    \label{table:res-dataset9-cartoon}
    \centering\small
     %\begin{adjustbox}{max width=\textwidth}
    \begin{tabular}{lll|lllll|l}\small
    Denoise & 
    Cost & 
    Initial $\alphavec$ & 
    Result $\alphavec^*$ & 
    Value & 
    SSIM & 
    PSNR &
    Its. &
    Fig. \\
    \hline
    \restabll{krtv2sketch-l2-dataset9-linesearch1-bfgs1-beta0-uinitreg_0.100000_}
            {krtv2}{l2}{$\TVINIT$}{fig:res-dataset9-cartoon}
    \restabll{tv-l2-dataset9-linesearch1-bfgs1-beta0-uinitreg_0.100000_}
            {tv}{l2}{$0.1/\pixels$}{fig:res-dataset9-cartoon}
    \end{tabular}
    %\end{adjustbox}
\end{table}

%\FloatBarrier

\section{Learning the data model}\label{sec:data}
%!TEX root = article.tex

The correct mathematical modelling of the data fidelity terms $\phi_i, i=1,\ldots,M$ in \eqref{eq:learn} is crucial for the design of a realisitc denoising model. Their choice corresponds to physical and statistical properties of the noise distribution corrupting the ground-truth $f_0$ and varies significantly depending on applications. Typically, the noise is assumed to be additive, Gaussian-distributed with $0$ mean and variance $\sigma^2$ determining the noise intensity. This assumption is reasonable in most of the applications because of the Central Limit Theorem. However, there are cases where this modelling assumption does not correspond to the actual statistical properties characterising the physics of the application considered. For instance, when considering astronomical images, different physical properties corresponding to the quantised (discrete) nature of light and to the independence of photons detection lead to consider a \emph{Poisson} noise distribution, which is signal dependent. Impulse noise seems to be more appropriate for modelling transmission errors affecting only some of the pixels in the image. For those pixels, the intensity value of the signal is switched to either the maximum/minimum value of the dynamic range of the image intensity or to a random value, with positive probability. 

For what follows, we will focus on these three noise distributions and on their possible combination. Other  distributions can be considered as well: in general, they suit specific applications (like radar or medical ultrasound images) where intrinsically the noise corrupting the image cannot be considered signal-independent. 

\medskip

From a mathematical point of view, variational models reflecting the statistical properties of the noise have been derived for the design of consistent denoising models. Starting from the pioneering work of Rudin, Osher and Fatemi \cite{rudin1992nonlinear}, in the case of Gaussian noise a $L^2$-type data fidelity $\phi$ is typically considered. In the case of impulse noise, a variational model based on the use of the $L^1$ norm has been considered in \cite{nikolova2004variational}: statistically, this corresponds to consider a Laplace distribution. Poisson noise-based models have been considered in several papers by approximating such distribution with a weighted-Gaussian distribution through variance-stabilising techniques \cite{starckmurtagh1994, papoutsellis2014}. In \cite{sawatzky2009total} a statistically-consistent analytical modelling for Poisson noise distributions has been derived: this results in a Kullback-Leibler-type fidelity.

As a result of different physical factors, very often in applications the presence of different noise distributions has to be considered as well. In \cite{hintermuller2013subspace} a combined $L^1$-$L^2$ TV-based model is considered for impulse and Gaussian noise removal. A two-phase approach is considered in \cite{impulsegauss2008} where the selection of the $L^1/L^2$ term is performed depending on the intensity of the noise. In general, though, the literature on these combined noise models is rather scarse. Gaussian-Poisson noise mixture has been considered in several papers from different point of views: in \cite{poissongauss2013} the exact log-likelihood estimator of the model is derived and then computed via a primal-dual splitting, while in other works (see, e.g.,  \cite{Foigausspoiss}) the discrete-continuous nature of the model (due to the Poisson-Gaussian component, respectively) is approximated by neglecting or modifying one of the two noise models, typically by means of variance-stabilising techniques or a weighted-$L^2$ approximation. 

\smallskip

We now proceed differently from Section \ref{sec:l2theory} and focus on the modelling of the optimal fidelity terms $\phi_i$ best fitting the acquired data, providing some examples for the single and multiple noise estimation case. In particular, we focus on the estimation of the optimal fidelity weights $\lambda_i, i=1,\ldots, M$ appearing in \eqref{eq:learn} and \eqref{eq:learn-numerical-single}, focusing on the Total-Variation regularisation \eqref{eq:TVreg} only applied to denoising problems. Compared to Section \ref{subsec:abstract model}, this corresponds to fix $\SPACEalphaPos=\left\{1\right\}$ and $K=Id$.  We base our presentation on \cite{de2013image,calatronidynamic}, where a careful analysis in term of well-posedness of the problem and derivation of the optimality system in this framework is carried out. 

\paragraph{Shorthand notation} 
In order not to make the notation too heavy, we warn the reader that we will use a shorthand notation for the quantities appearing in the regularised problem \eqref{eq:learn-numerical-single}, that is we will write $\Phi_i(v)$ for the data fidelities $\phi_i(x,v), i=1\ldots,M$ and $u$ for $u_{\lambda,\gamma,\mu}$, the minimiser of $J^{\gamma,\mu}(\cdot;\lambda)$.
\subsection{Single noise estimation}
\label{sec:numerics single}

In this section we consider the one-noise distribution case ($M=1$)  where we aim to determine the constant optimal fidelity weight $\lambda$ by solving the following optimisation problem:
\begin{subequations} \label{eq:optimization problem}
\begin{equation}
\min_{\lambda \geq 0} ~ \frac{1}{2}\| f_0-u\|^2_{L^2}
\end{equation}
subject to (compare \eqref{alg eq: denoising})
\begin{multline}  \label{eq:optimization problem eq 2}
\mu \iprod{\grad u}{\grad(v-u)}_{L^2} +\lambda ~ \int_{\Omega}
\Phi'(u)(v-u) \d x\\+ \int_{\Omega} \norm{\grad v} \d x-\int_{\Omega}
\norm{\grad u} \d x \geq 0 \, \, \text{ for all } v \in H_0^1(\Omega),
\end{multline}
\end{subequations}
where the fidelity term $\Phi$ will change according to the different noise distributions considered and the pair $(f_0,f)$ is the training pair composed by a noise-free and noisy version of the same image, respectively. 
%We then look for the optimal $\lambda^*$ such that the Huber-regularised TV-reconstruction of $f$ computed in correspondence of $\lambda^*$ matches at best the corresponding  ground truth $f_0$ in the $L^2$-sense. 

\smallskip

Note that in the case the noise level is known there are classical techniques in inverse problems for choosing an optimal parameter $\lambda^*$ in a variational regularisation approach, e.g. the discrepancy principle or the L-curve approach \cite{engl1996regularization}. In our discussion we do not use any knowledge of the noise level but rather extract this information indirectly from our training set and translate it to the optimal choice of $\lambda$. As we will see later such an approach is also naturally extendable to multiple noise models as well as inhomogeneous noise.

\paragraph{Gaussian noise} 

We start by considering \eqref{eq:optimization problem} for determining the regularisation parameter in the standard TV denoising model assuming that the noise in the image is normally distributed. In this case the fidelity term reads $\Phi(u)=|u-f|^2$. The optimisation problem  \ref{eq:optimization problem} takes the following form:
\begin{subequations} \label{eq:gaussian noise problem}
\begin{equation}
\min_{\lambda \geq 0}  ~\frac{1}{2} \|f_0 - u\|^2_{L^2} 
\end{equation}
subject to:
\begin{multline}
\mu \iprod{\grad u}{\grad(v-u)}_{L^2} +\int_{\Omega} \lambda (u-f) (v-u) ~dx\\+ \int_{\Omega} \norm{\grad v} \d x-\int_{\Omega}  \norm{\grad u} \d x \geq 0, \forall v \in H_0^1(\Omega).   \label{eq:gaussian noise problem lower}
\end{multline}
\end{subequations}

%The problem consists therefore in the optimal choice of the TV regularization parameter, if the original image is known in advance. This is a toy example for proof of concept only. In practice this image would be replaced by a training set of images as motivated in the introduction.

For the numerical solution of the regularised variational inequality we use a 
primal-dual algorithm  presented in \cite{hintermuller2006infeasible}.

\begin{figure}[h!] 
\begin{center}
\hspace{0.01cm} \includegraphics[height=4cm,width=0.43\textwidth]{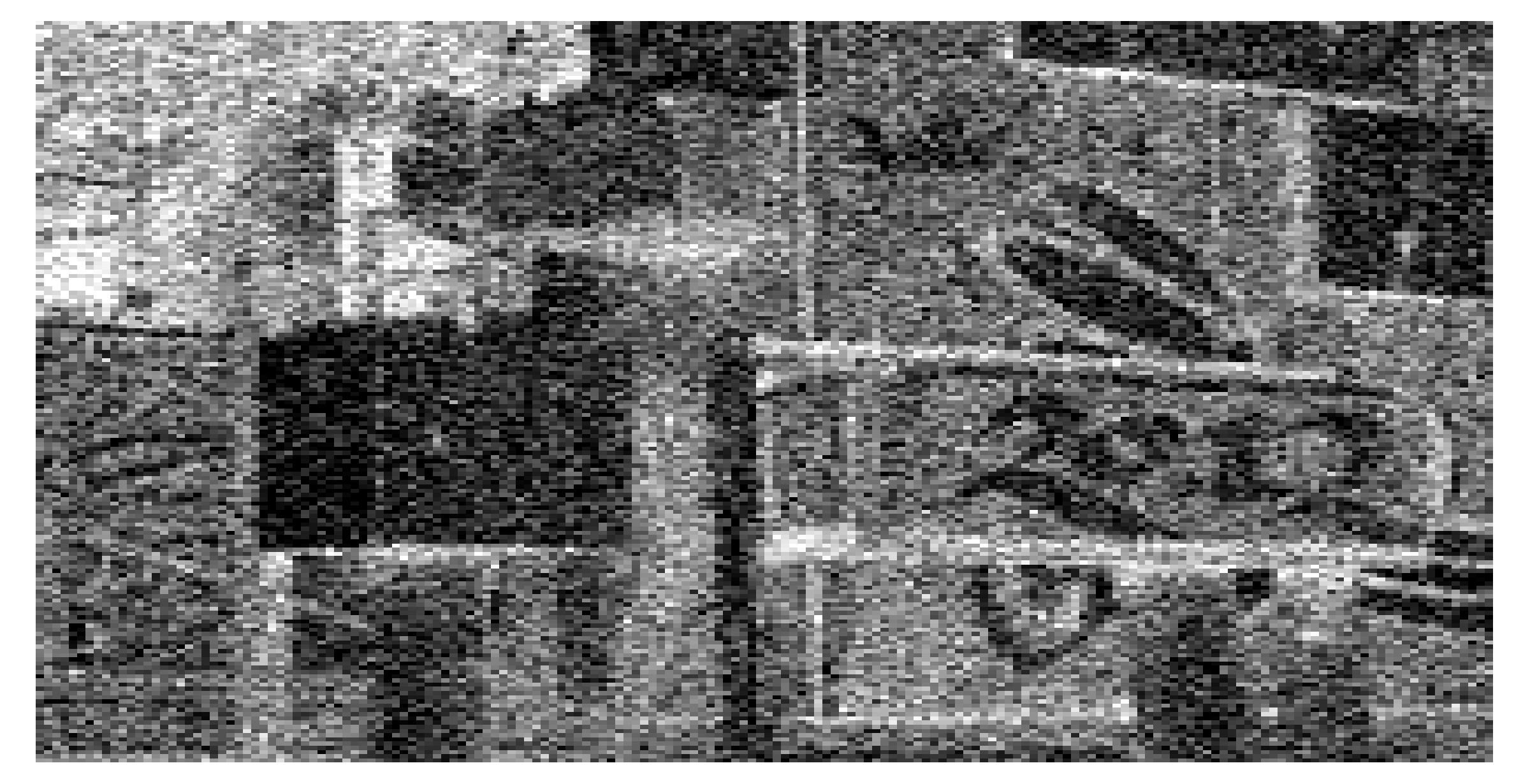} \hspace{0.15cm} \includegraphics[height=4cm,width=0.43\textwidth]{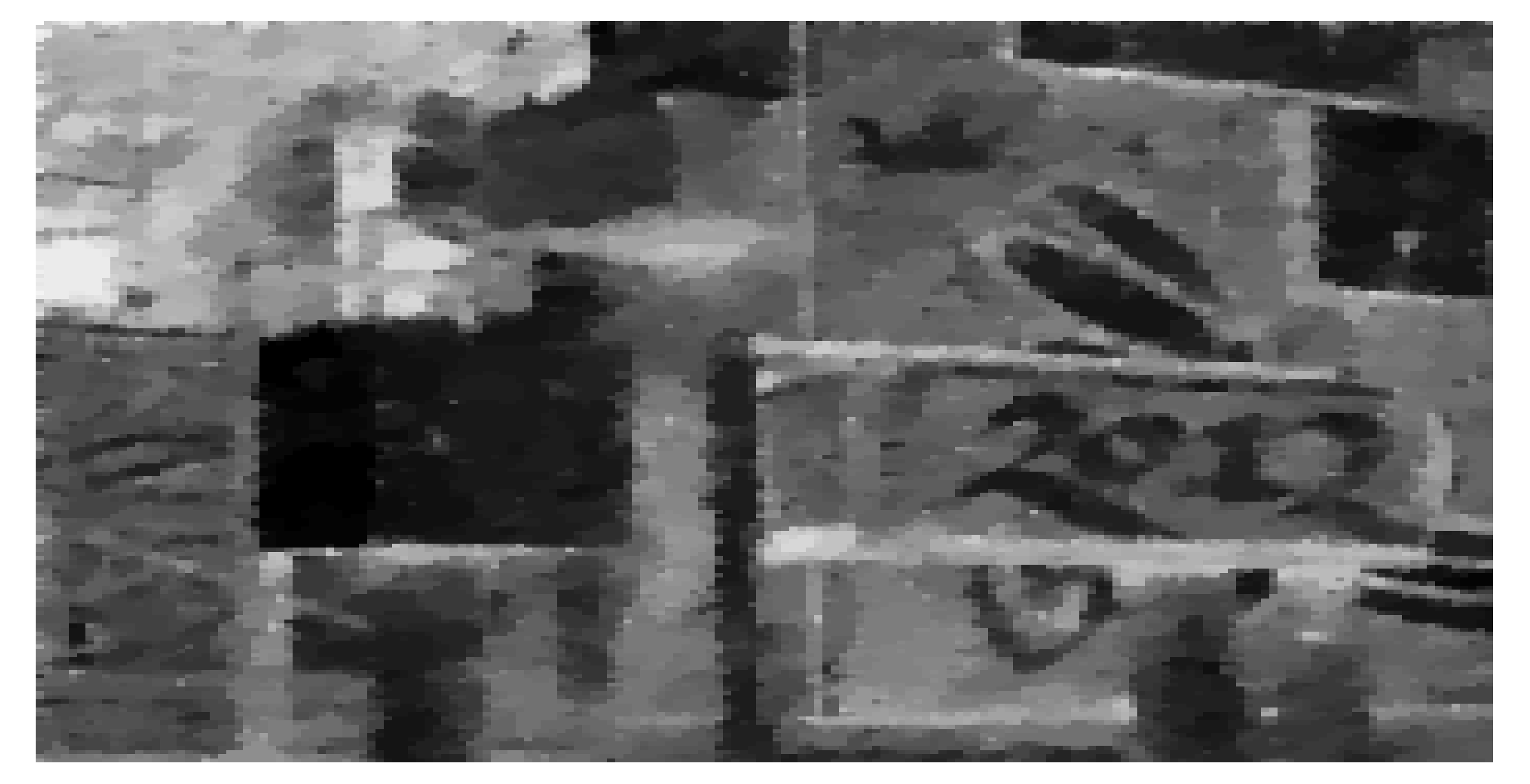}
\end{center}
\caption{Noisy (left) and optimal denoised (right) image. Noise variance:
  0.02. Optimal parameter $\lambda^*=1770.9$.}
  \label{fig:experiment1}
\end{figure}

As an example, we compute the optimal parameter  $\lambda^*$ in \eqref{eq:gaussian noise problem} for a noisy image distorted by Gaussian noise with zero mean and variance $0.02$ . Results are reported in Figure \ref{fig:experiment1}. The optimisation result has been obtained for the parameter values $ \mu = 1e-12, \gamma=100$ and $h=1/177$.

%We observe from experiments that increasing the variance of the noise changes significantly the value of the optimal parameter. This is intuitively clear, since as the
%image becomes noisier less information that can be directly
%obtained. When that happens, the TV-regularisation plays an
%increasingly important role.

\begin{figure}[h!]
\begin{center}
\includegraphics[height=3.5cm,width=6.cm]{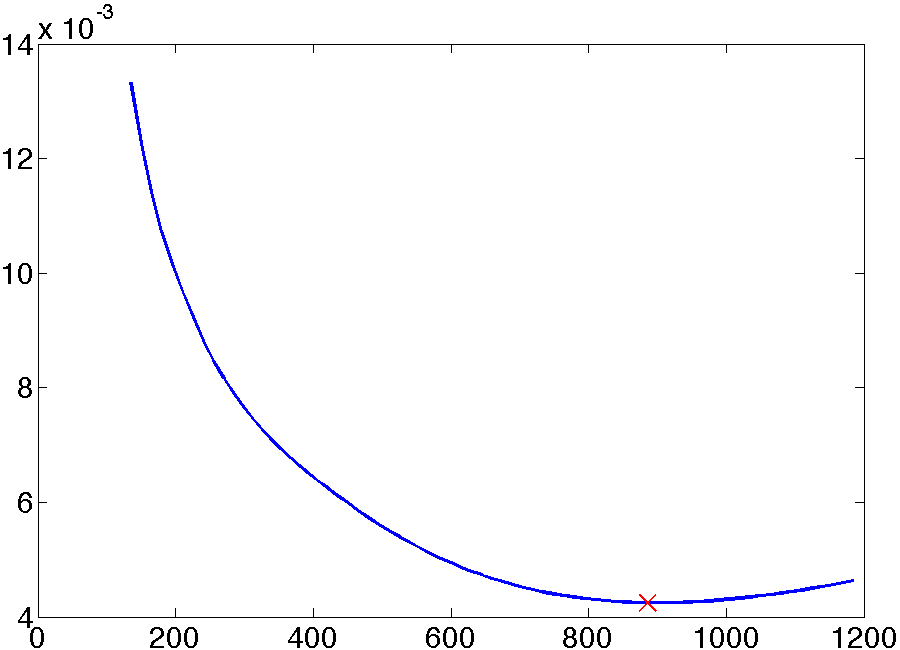} \hfill \includegraphics[height=3.5cm,width=6.cm]{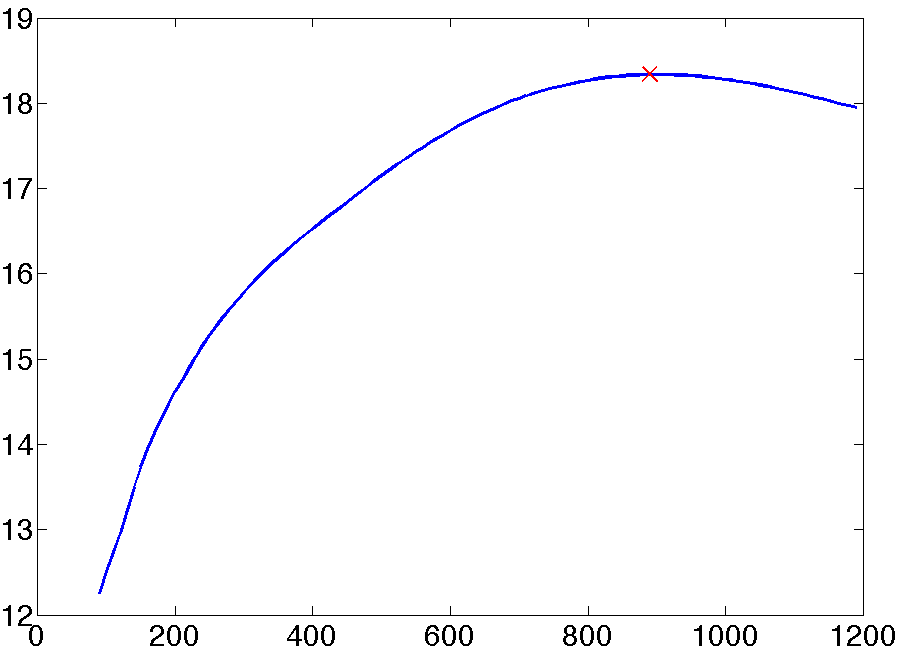}
\end{center}
\caption{Plot of the cost functional value (left) and the SNR (right) vs. the parameter $\lambda$. Parameters: the input is the $80\times 80$ pixel crop of the bottom left corner of the noisy image in Figure \ref{fig:experiment1}, $h=1/80, \, \gamma=100, \, \mu=1e-12.$ The red cross in the plot corresponds to the optimal $\lambda^*=885.5$.}
\label{fig:valuef}
\end{figure}

In order to check the optimality of the computed regularisation parameter $\lambda^*$, we consider the $80\times 80$ pixel bottom left corner of the noisy image in Figure \ref{fig:experiment1}. In Figure \ref{fig:valuef} the values of the cost functional and of the {\bf S}ignal to {\bf N}oise {\bf R}atio $SNR=20 \times \log_{10}\bigg(\frac{\|f_0\|_{L^2}}{\|u-f_0\|_{L^2}}\bigg),$ for parameter values between 150 and 1200, are plotted. Also the cost functional value corresponding to the computed optimal parameter $\lambda^*=885.5$ is shown with a cross. It can be observed that the computed weight actually corresponds to an optimal solution of the bilevel problem. Here we have used $h=1/80$ and the other parameters as above.

\medskip

The problem presented consists in the optimal choice of the TV regularisation parameter, if the original image is known in advance. This is a toy example for proof of concept only. In applications, this image would be replaced by a training set of images.
%Concerning mesh independency, in Table \ref{table: mesh dependence} the number of quasi-Newton iterations for different mesh size steps is shown.  A mesh independent behavior can be observed from the data.
%\begin{table}
%\centering
%\begin{tabular}{|c|c|c|c|c|} \hline
%$\hspace{.2cm} \# \text{mesh points} \hspace{.2cm}$ &
%\hspace{.2cm} 20 \hspace{.2cm} &
%\hspace{.2cm} 40 \hspace{.2cm} & \hspace{.2cm} 80 \hspace{.2cm} & \hspace{.2cm}
%160 \hspace{.2cm}\\
%\hline $\#$ iterations  &9 &9 &8 &10\\
%\hline
%\end{tabular}
%\vspace{0.2cm}\caption{Mesh independency of quasi-Newton method. For different-sized details of the noisy image in Figure \ref{fig:experiment1} with noise variance $0.02$ the number of quasi-Newton iterations is plotted against the different mesh sizes used. Parameters: $\varepsilon=1e-15$, $\gamma=100$, $\beta=1e-10$.} \label{table: mesh dependence}
%\end{table}

\paragraph{Robust estimation with training sets}

Gaussian noise images typically arise with\-in the framework of {\bf M}agnetic
{\bf R}esonance {\bf I}maging (MRI). The challenge in this case consists in
training the TV denoising method such that with one fixed optimally computed $\lambda^*$ clearer images are obtained from noisy acquisitions taken on a single MR tomograph with fixed settings. MR images seem to be a natural choice for our methodology, since a training set of images is often at hand. Let us consider a training database $\left\{(f_0^k,f_k)\right\}_{k=1,\ldots,K}, K\gg 1$ of clean and noisy images. We modify \eqref{eq:gaussian noise problem} as:
\begin{equation}  \label{optprob1gauss}
\min_{\lambda \geq 0}\,\frac{1}{2K}\sum_{k=1}^{K}\norm{f_0^k-u_k}^2_{L^2}
\end{equation}
subject to the set of regularised versions of \eqref{eq:gaussian noise problem lower}, for $k=1,\ldots,K$.

As explained in \cite{calatronidynamic}, dealing with large training sets of images and non-smooth PDE constraints of the form  \eqref{eq:gaussian noise problem lower} may result is very high computational costs as, in principle, each constraint needs to be solved in each iteration of the optimisation loop. On the other hand, in MRI applications, a large database of images is desirable in order to make the optimal noise estimation robust. In order to overcome the computational efforts, we estimate $\lambda^*$ using the Dynamic Sampling Algorithm \ref{alg:dynsamp}.

For the following numerical tests, the parameters are chosen as follows: $\mu=1e-12,$ $\gamma=100$ and $h=1/150$. The noise in the images has distribution $\mathcal{N}(0,0.005)$ and the accuracy parameter $\theta$ of the Algorithm \ref{alg:dynsamp}, is chosen to be $\theta=0.5$.
\smallskip

Table \ref{tablegauss} shows the numerical values of the optimal parameter $\lambda^*$ and $\lambda^*_S$ computed varying $N$ after solving all the PDE constraints and using Dynamic Sampling algorithm, respectively. We measure the efficiency of the algorithms in terms of the number of the PDEs solved during the whole optimisation and we compare the efficiency of solving \eqref{optprob1gauss} subject to the whole set of constraints \eqref{eq:gaussian noise problem lower} with the one where solution is computed by means of the Dynamic Sampling strategy, observing a clear improvement. %This corresponds to an increasing number of BFGS iterations which does not appear to be an issue as BFGS iterations are themselves very fast. For the sake of computational efficiency, what really matters is the number of PDEs that need to be solved in \emph{each} iteration of BFGS. Moreover, thanks to modern parallel computing methods and to the decoupled nature of the constraints in each BFGS iteration, solving such a reduced amount of PDEs makes the computational efforts very reasonable. We note that the size of the sample is generally maintained very small in comparison to $N$ or just slightly increased. 
Computing also the relative error $\norm{\hat{\lambda}_S-\hat{\lambda}}_1/\norm{\lambda_S|}_1$ we note a good level of accuracy: the error remains always below $5\%$.
%\vspace{-0.5cm}
\begin{table}[t]
\centering
\scriptsize    
\begin{tabular}{| c | c | c | c | c | c | c | c | c | c | } 
    \hline
$K$ & $\lambda^*$ & $\lambda^*_S$ & $|S_0|$ & $|S_{end}|$ & eff. & eff. Dyn.S. & BFGS its. & BFGS its. Dyn.S.& diff. \\ \hline
$10$ & $3334.5$ & $3427.7$ &$2$&$3$& $140$ &\bm{$84$}& $7$ &$21$& $2.7\%$ \\ \hline
$20$ & $3437.0$ & $3475.1$ &$4$&$4$& $240$ &\bm{$120$}& $7$ &$15$& $1.1\%$ \\ \hline
$30$ & $3436.5$ & $3478.2$ &$6$&$6$& $420$ &\bm{$180$}& $7$ & $15$ & $1.2\%$ \\ \hline
$40$ & $3431.5$ & $3358.3$ &$8$&$9$& $560$ &\bm{$272$}& $7$ &$16$& $2.1\%$ \\ \hline
$50$ & $3425.8$ & $3306.4$ &$10$ & $10$ & $700$ & \bm{$220$} & $7$ & $11$ & $3.5\%$ \\ \hline
$60$ & $3426.0$ & $3543.4$ & $12$ & $12$ & $840$ & \bm{$264$} & $7$ &$11$ & $3.3\%$ \\ \hline
$70$ & $3419.7$ & $3457.7$ & $14$ & $14$ &  $980$ & \bm{$336$ }& $7$ & $12$ & $1.1\%$ \\ \hline
$80$ & $3418.1$ & $3379.3$  & $16$ & $16$ & $1120$ & \bm{$480$} &$7$ & $15$ & $<1\%$ \\ \hline
$90$ & $3416.6$ & $3353.5$  & $18$  & $18$  & $1260$ & \bm{$648$}  & $7$ & $18$  & $2.3\%$ \\ \hline
$100$ & $3413.6$ & $3479.0$ & $20$ & $20$ & $1400$ & \bm{$520$} & $7$ & $13$ & $1.9\%$ \\ \hline
 \end{tabular}
%\caption{$N$ is the size of the database, $\hat{\lambda}$ is the optimal parameter obtained by solving all the $N$ constraints, whereas $\hat{\lambda}_S$ is the one computed by solving the problem with Algorithm 1. We compare the accuracy of the  approximation in terms of the  difference $\norm{\hat{\lambda}_S-\hat{\lambda}}_1/\norm{\lambda_S|}_1$.}
\caption{Optimal $\lambda^*$ estimation for large training sets: computational costs are reduced via Dynamic Sampling Algorithm \ref{alg:dynsamp}.}
\label{tablegauss}
\end{table}

\begin{figure}[h!]
\begin{center}
\includegraphics[width=2cm,height=2cm]{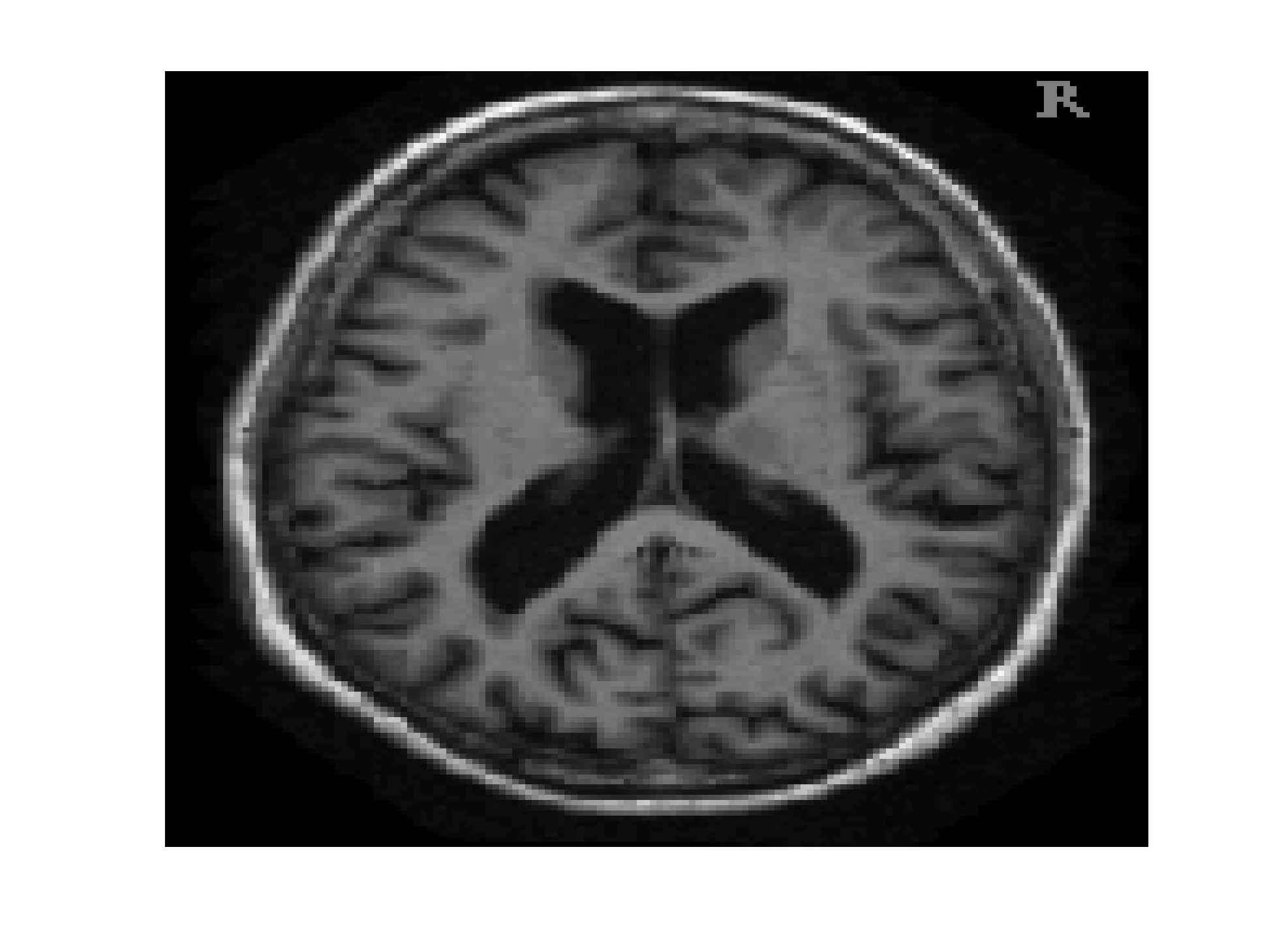}
\hspace{-0.3cm}
\includegraphics[width=2cm,height=2cm]{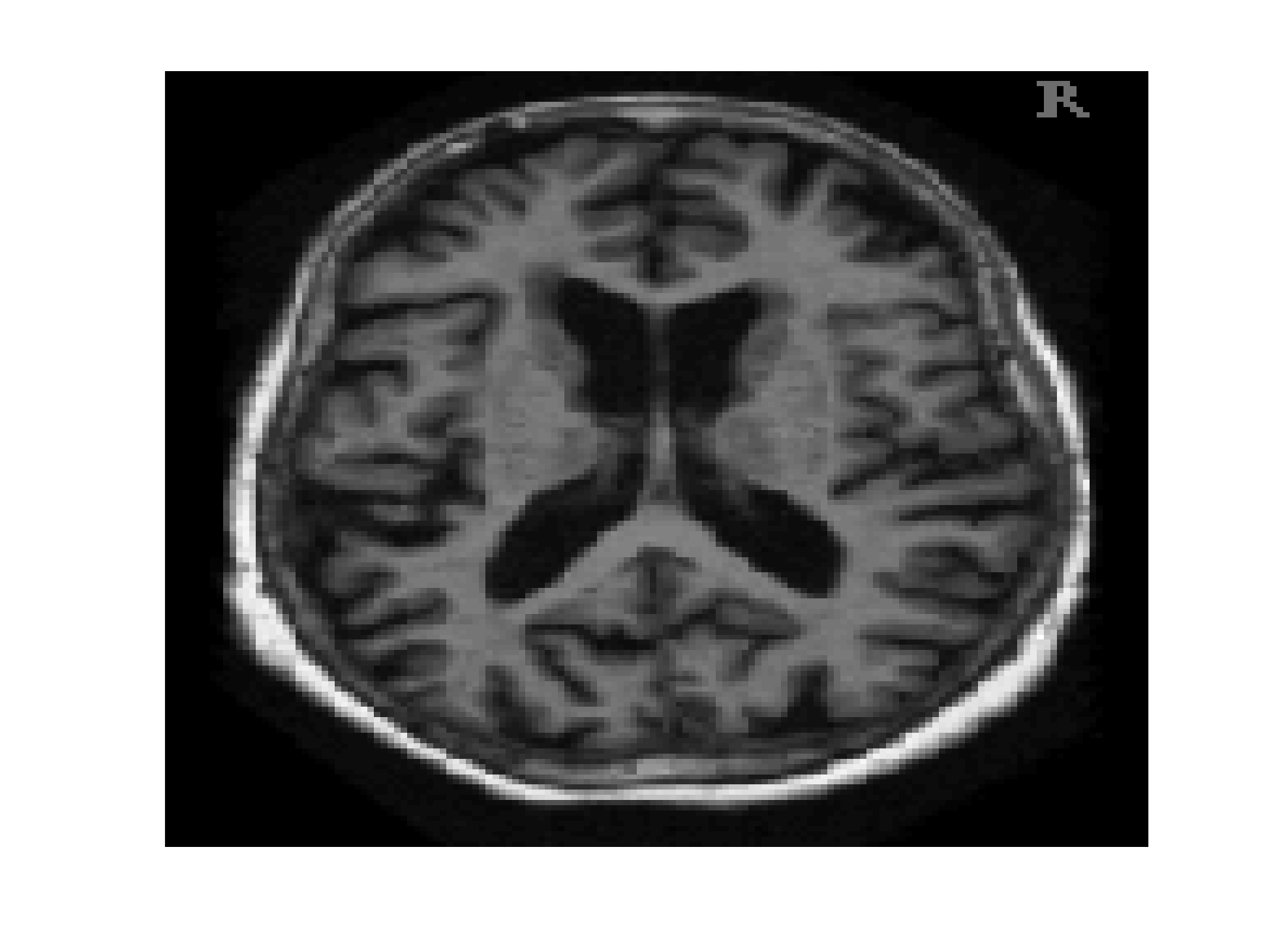}
\hspace{-0.3cm}
 \includegraphics[width=2cm,height=2cm]{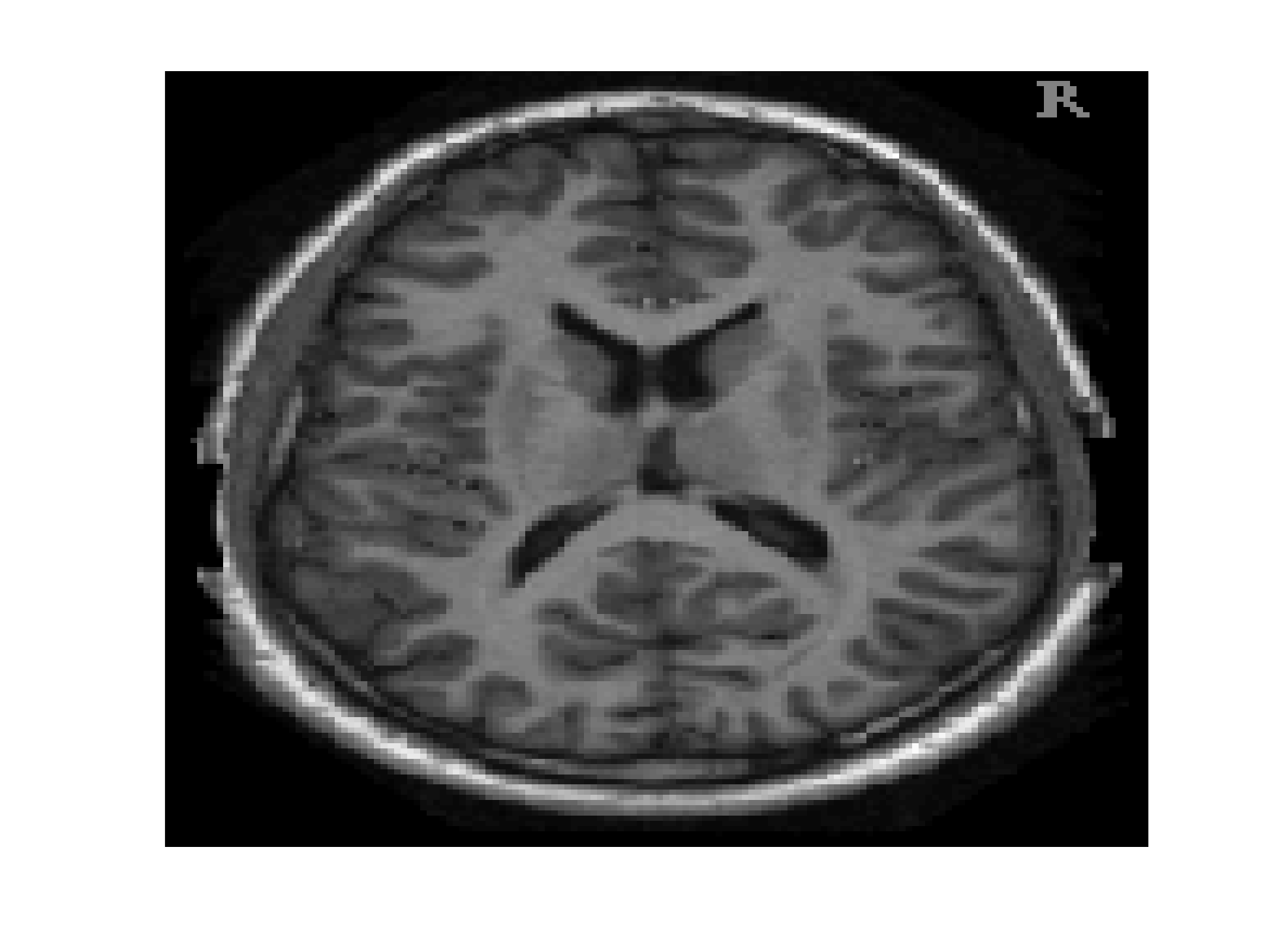}
\hspace{-0.3cm}
 \includegraphics[width=2cm,height=2cm]{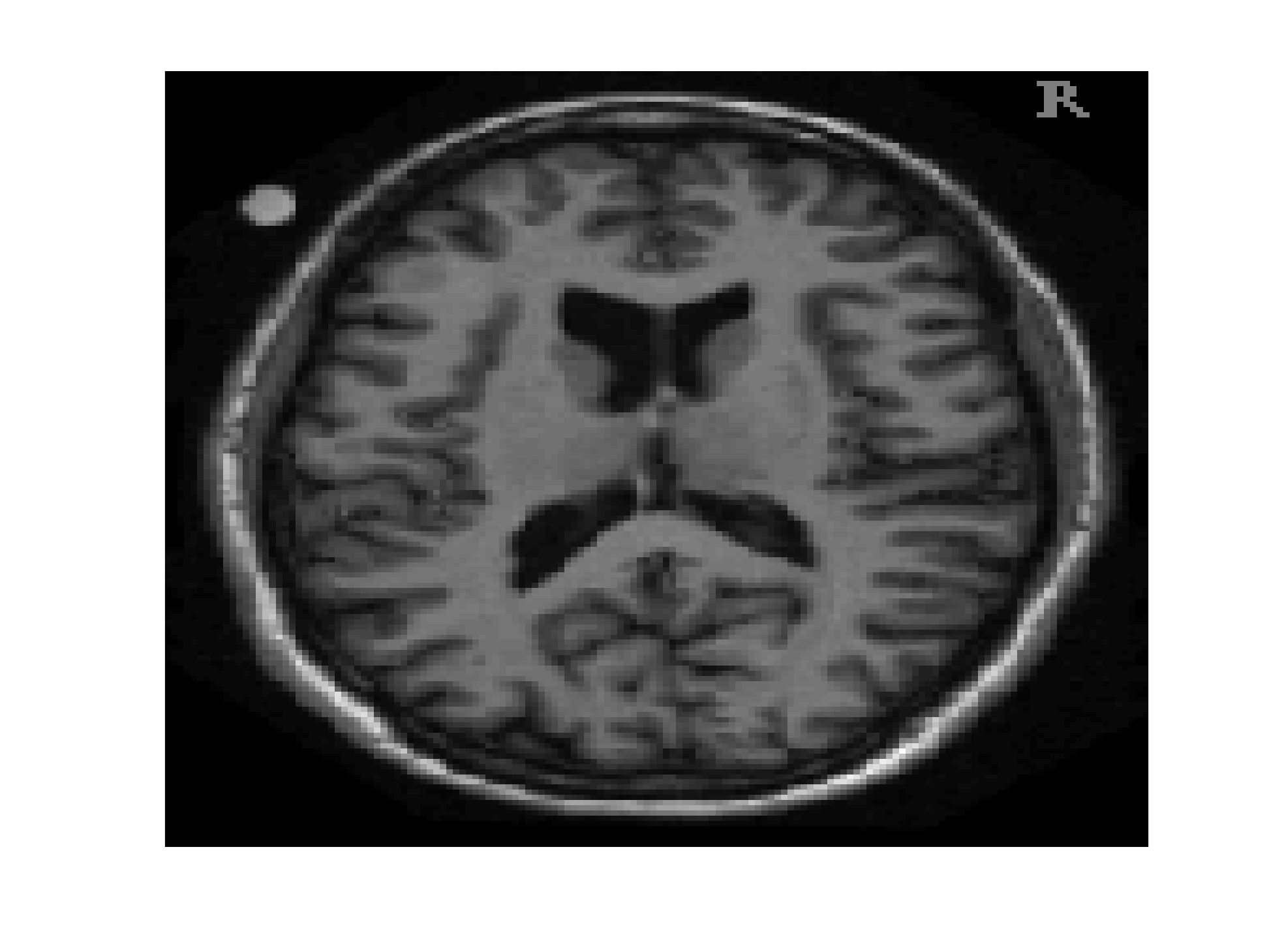}
\hspace{-0.3cm}
\includegraphics[width=2cm,height=2cm]{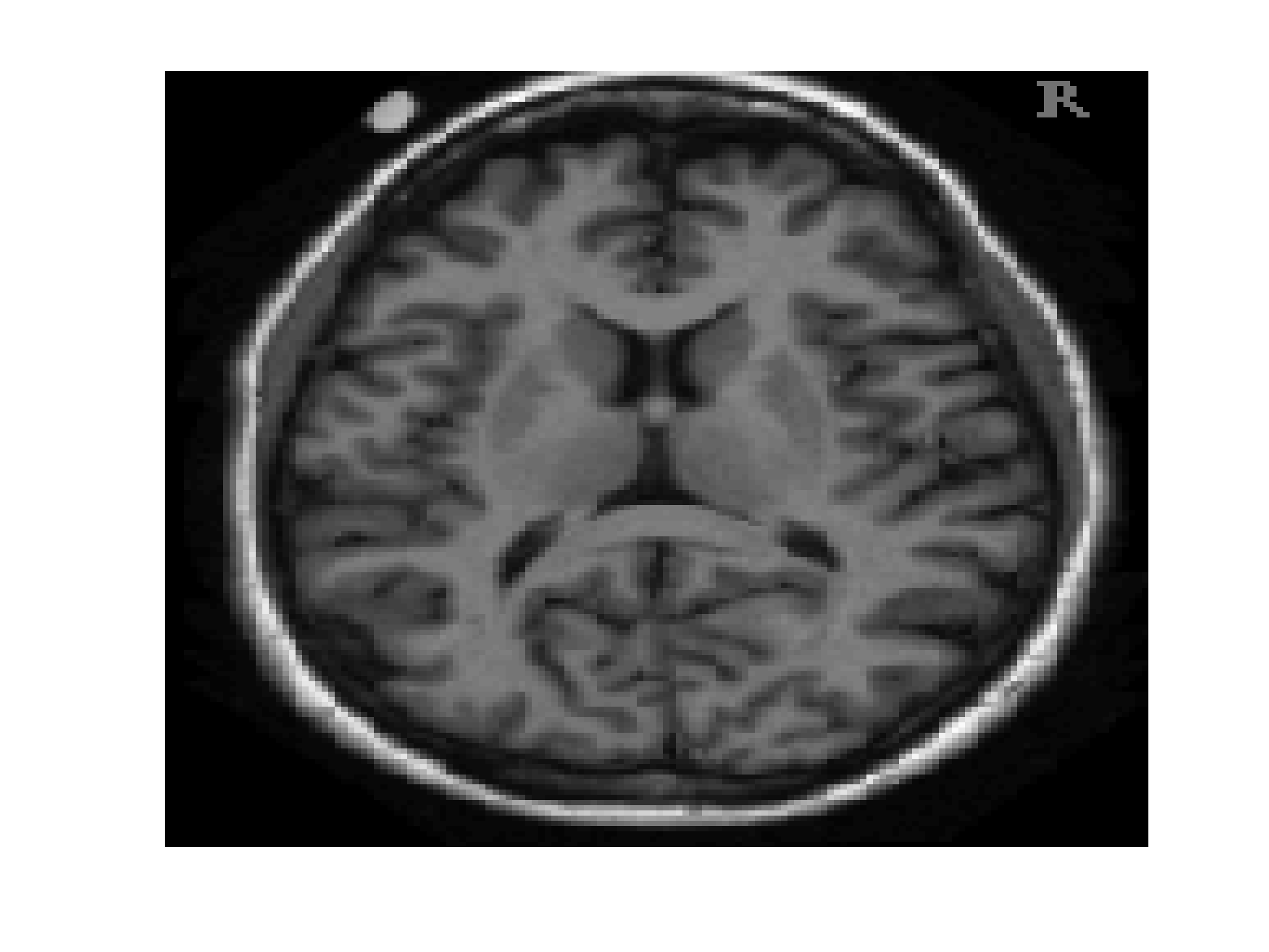}\\
\includegraphics[width=2cm,height=2cm]{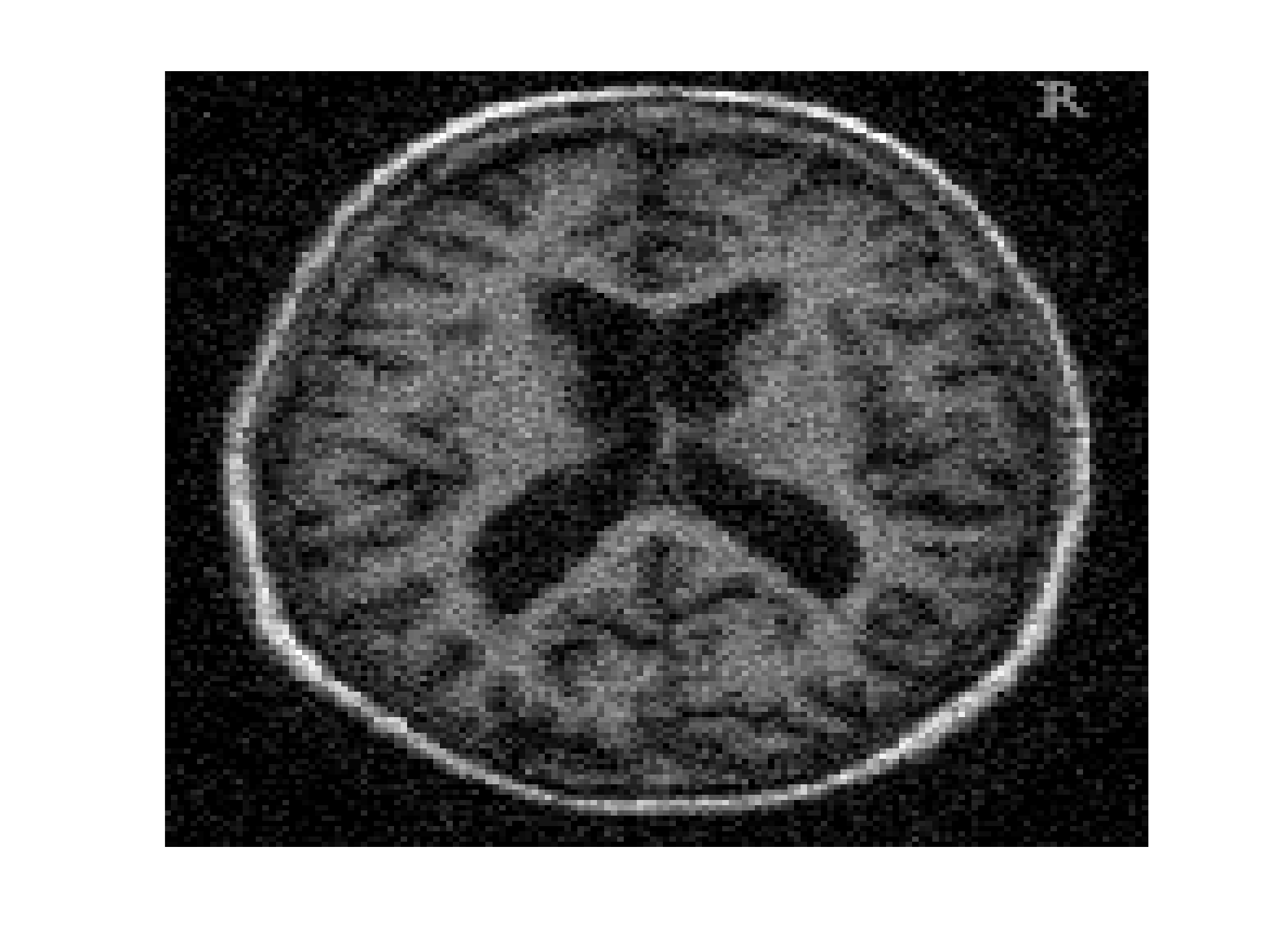}
\hspace{-0.3cm}
\includegraphics[width=2cm,height=2cm]{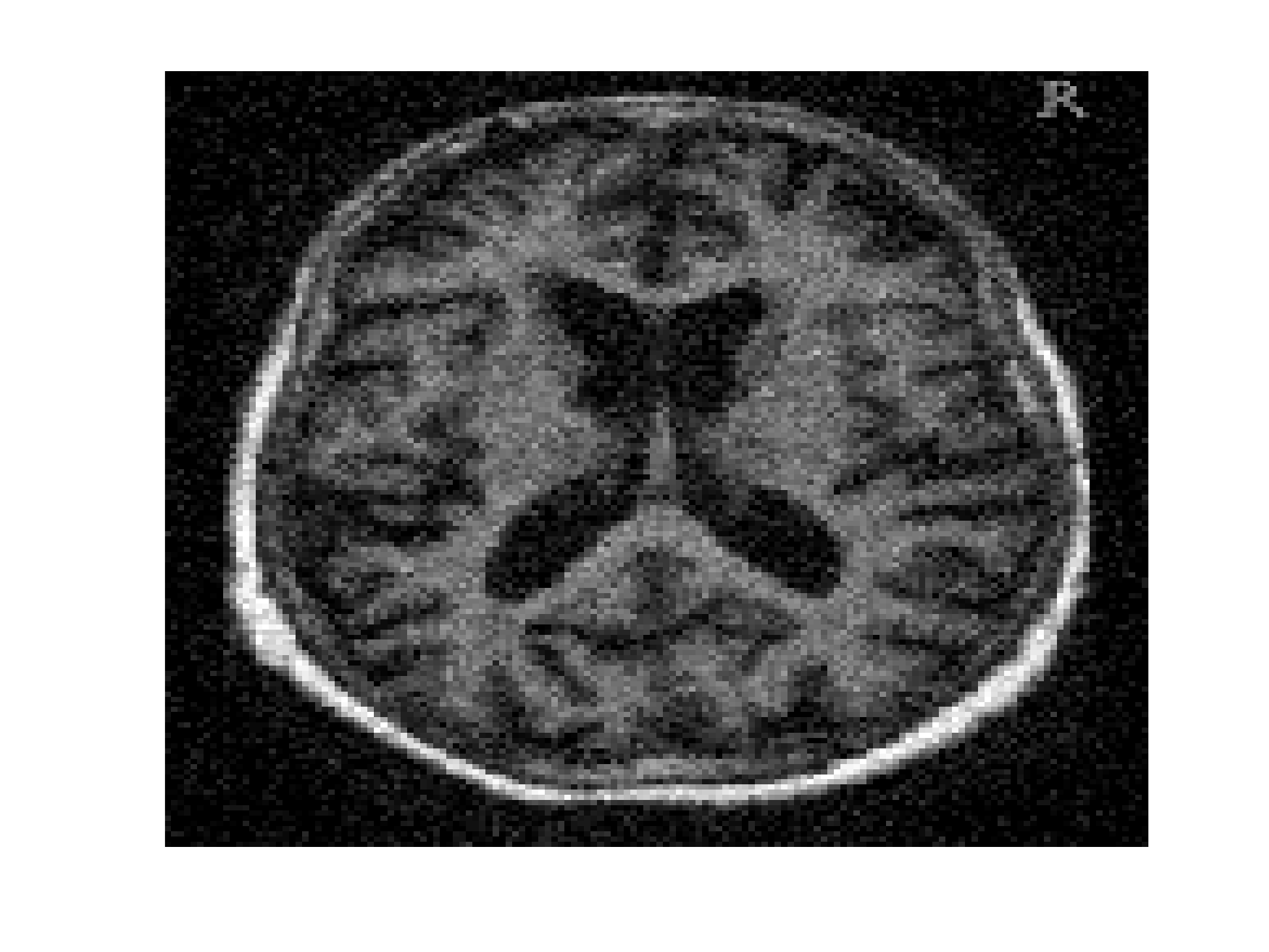}
\hspace{-0.3cm}
 \includegraphics[width=2cm,height=2cm]{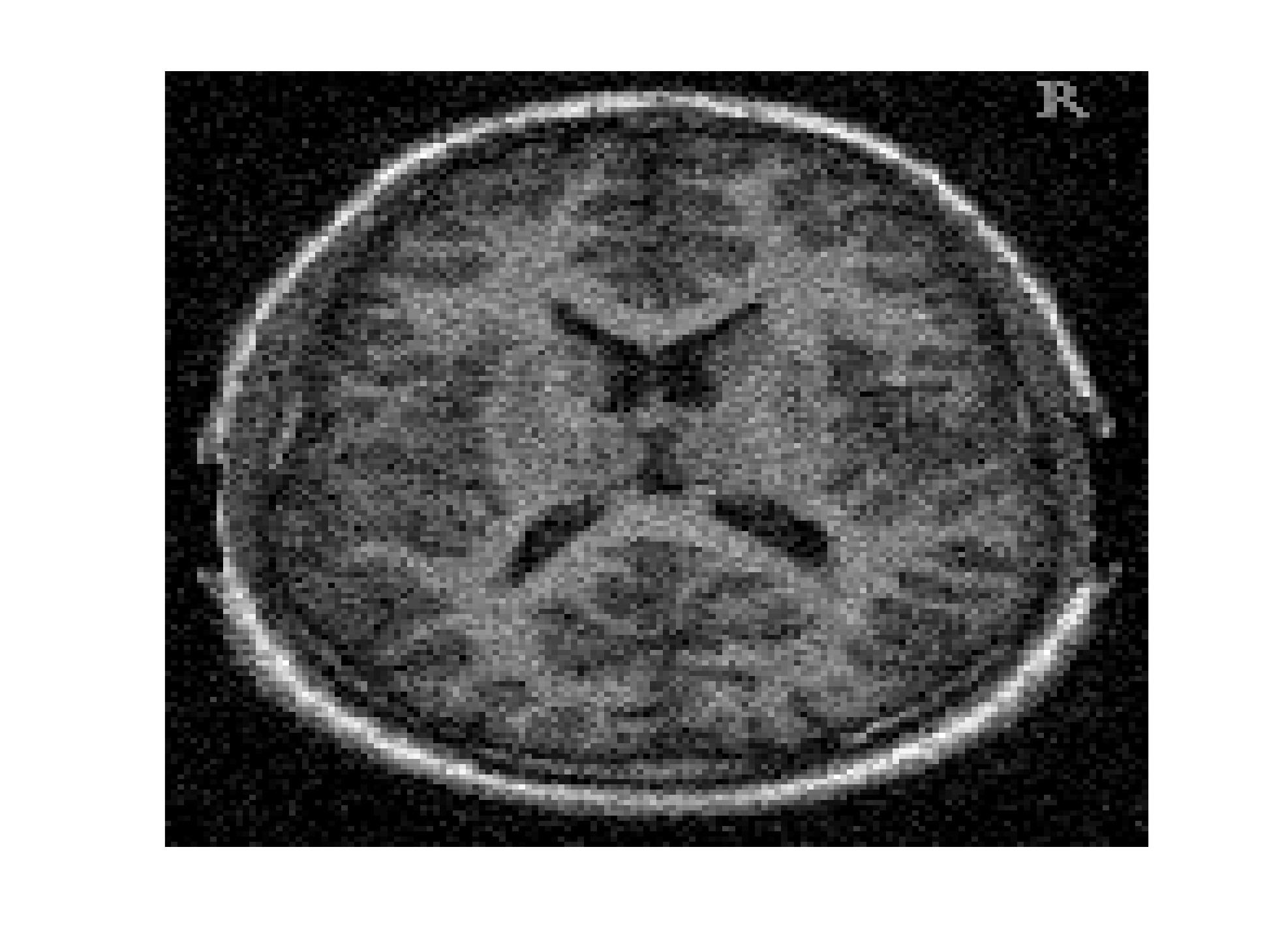}
\hspace{-0.3cm}
 \includegraphics[width=2cm,height=2cm]{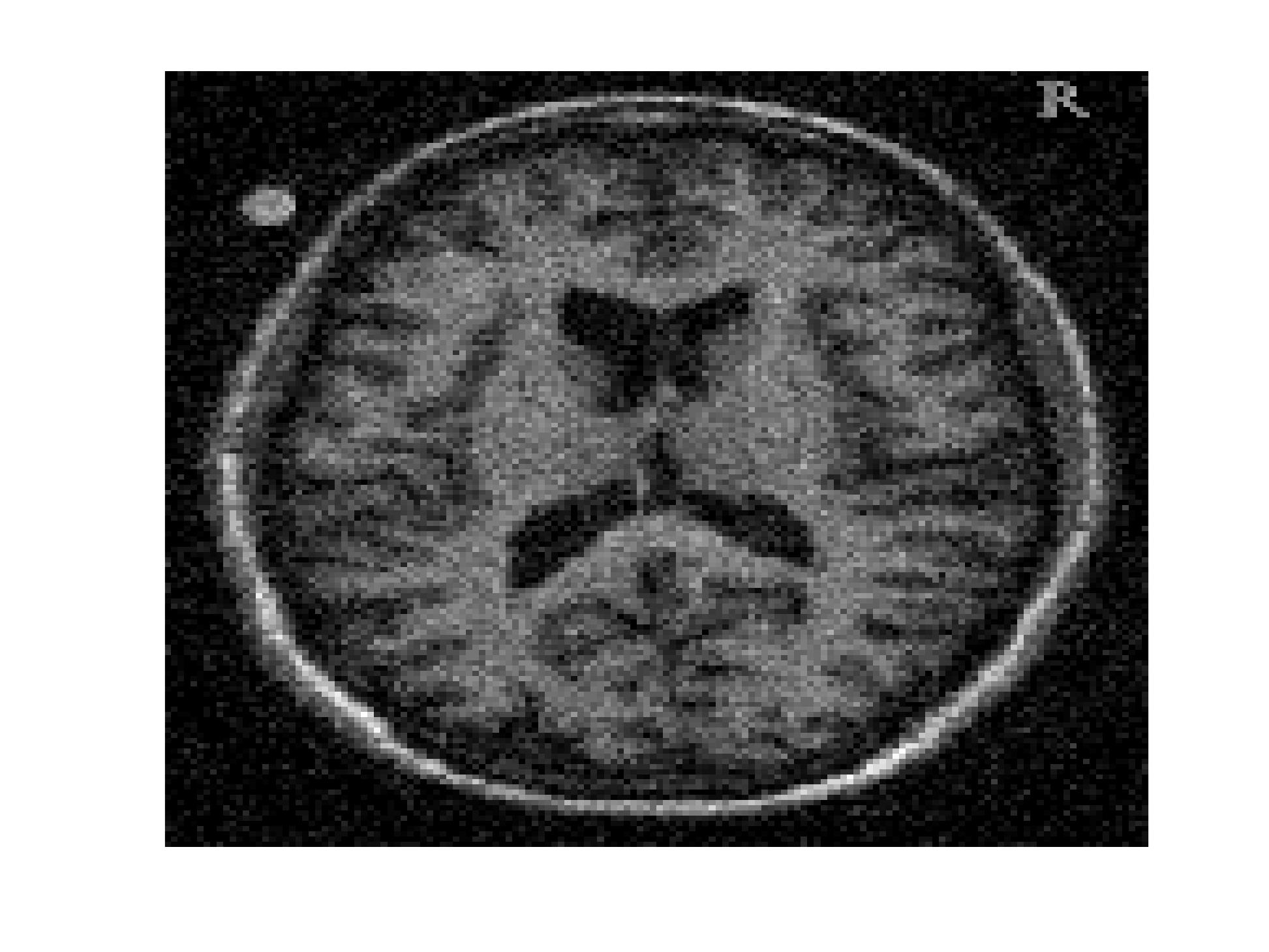}
\hspace{-0.3cm}
\includegraphics[width=2cm,height=2cm]{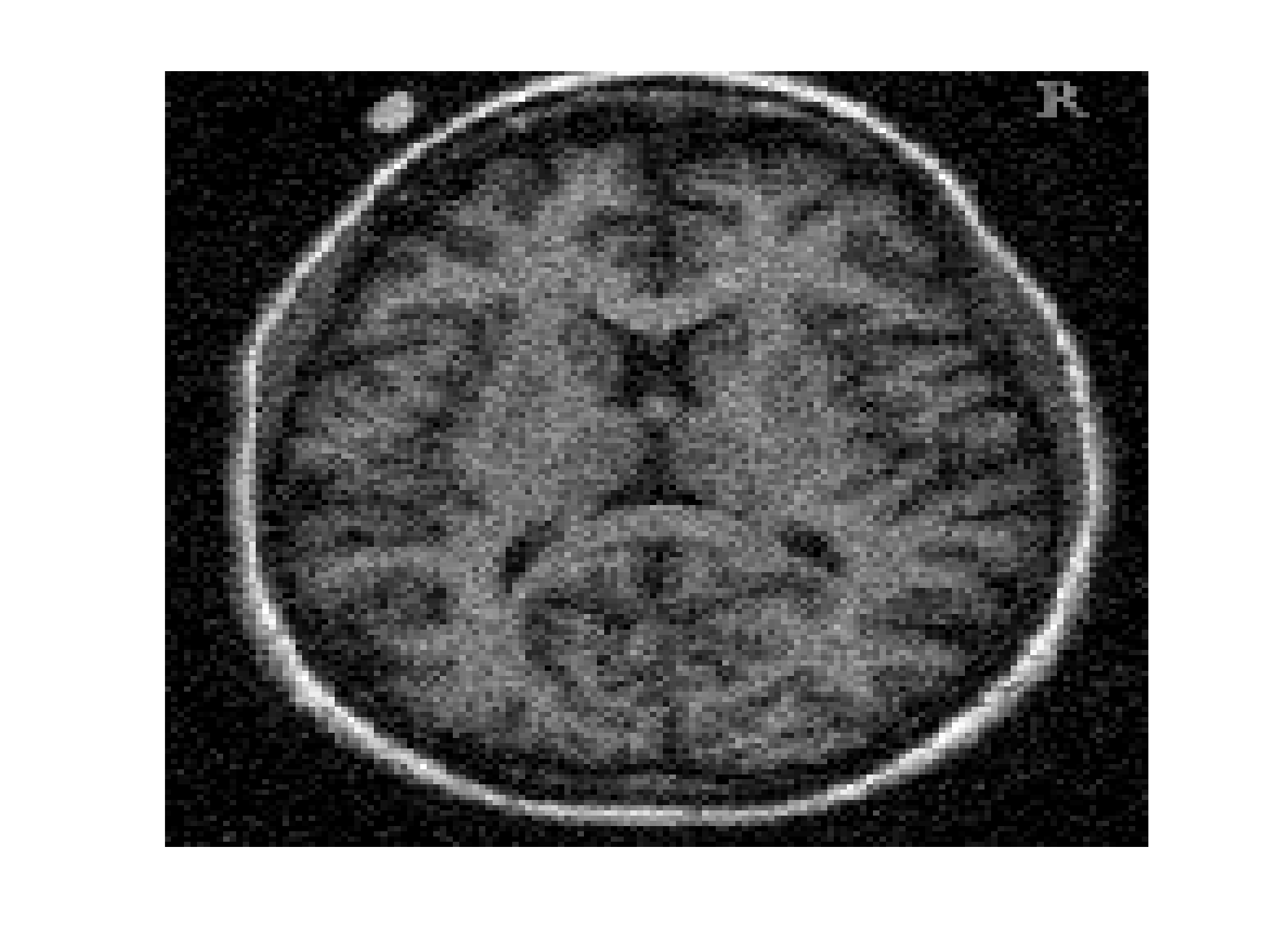}\\
\includegraphics[width=2cm,height=2cm]{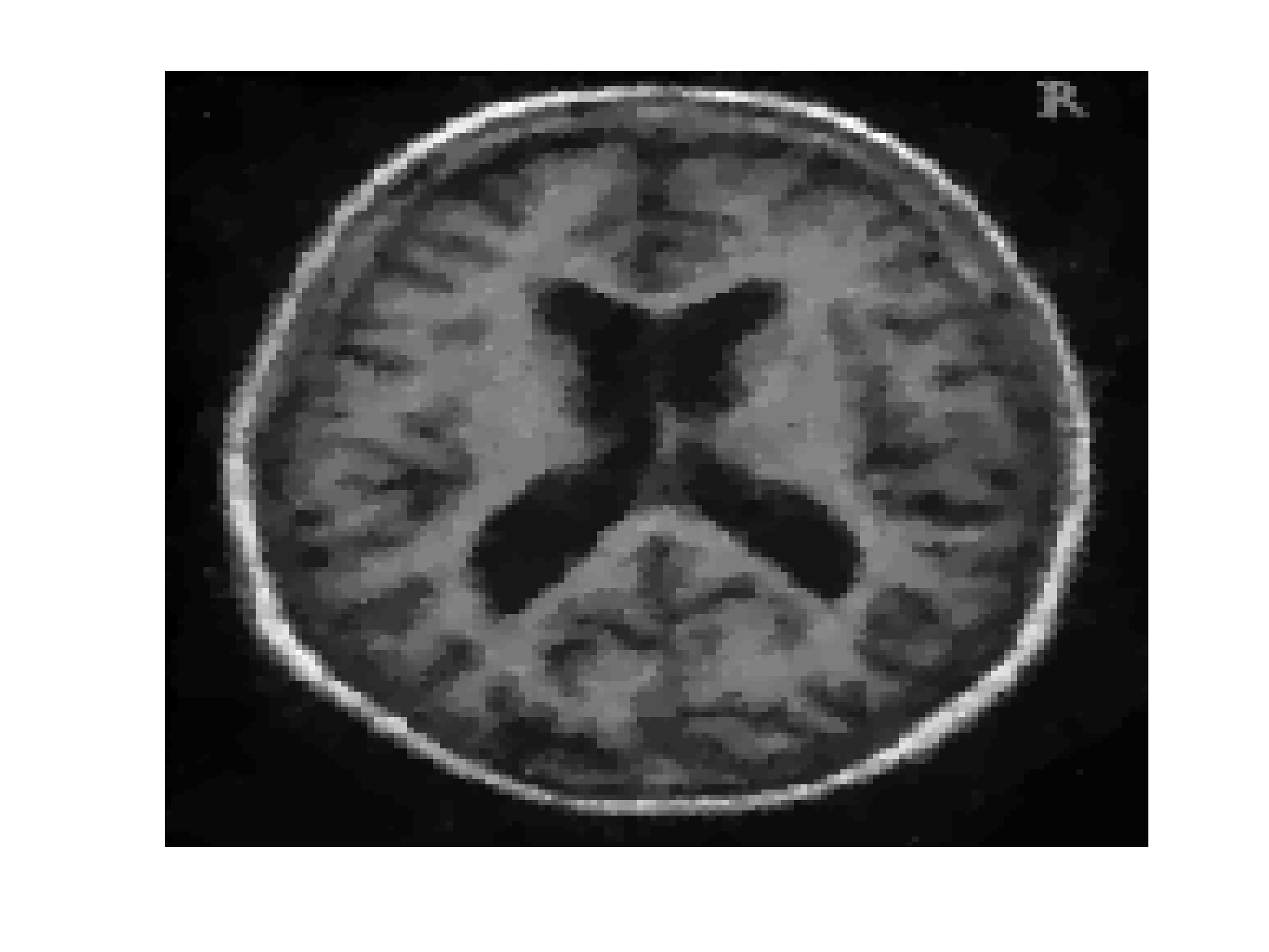}
\hspace{-0.3cm}
\includegraphics[width=2cm,height=2cm]{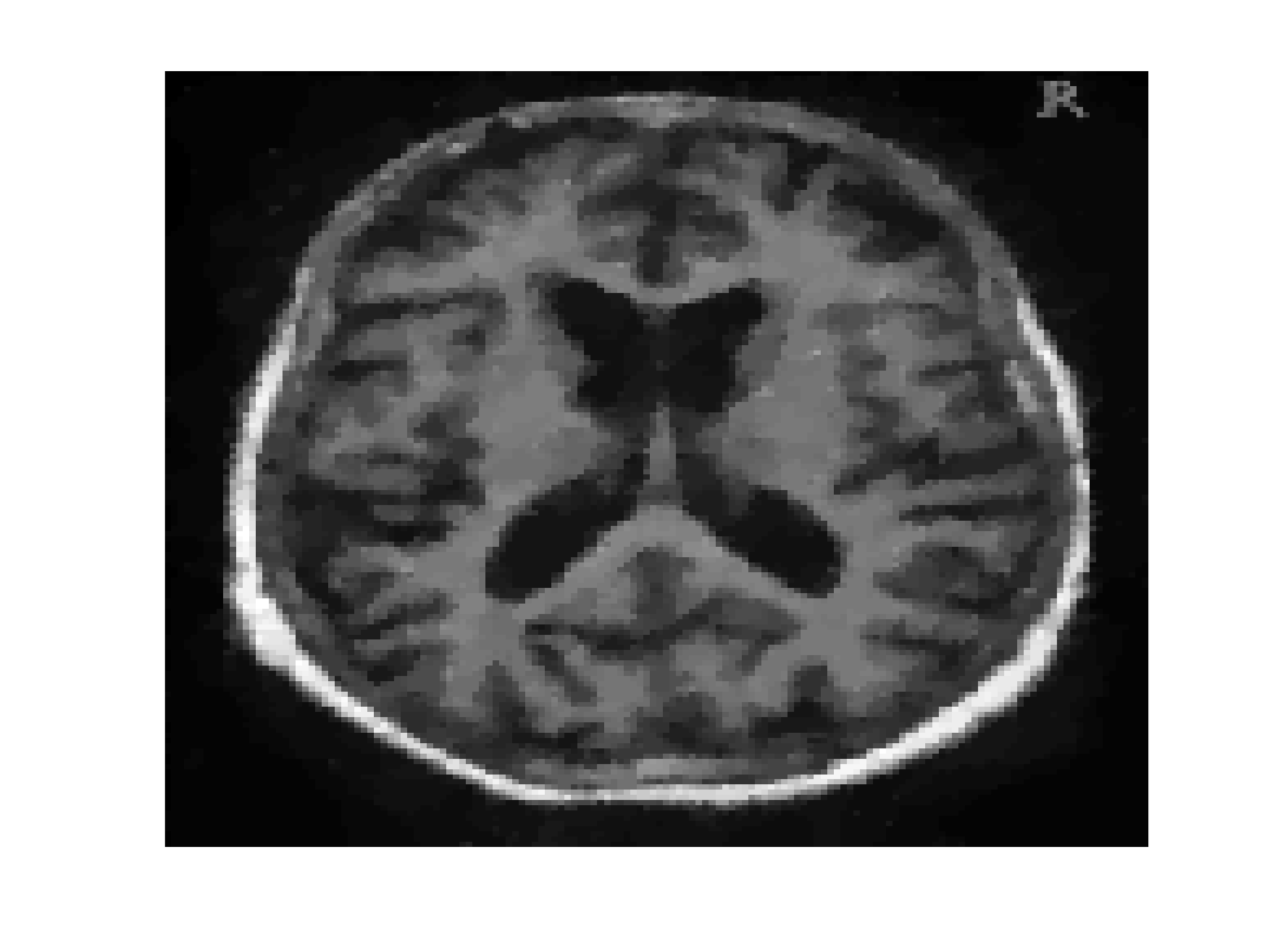}
\hspace{-0.3cm}
 \includegraphics[width=2cm,height=2cm]{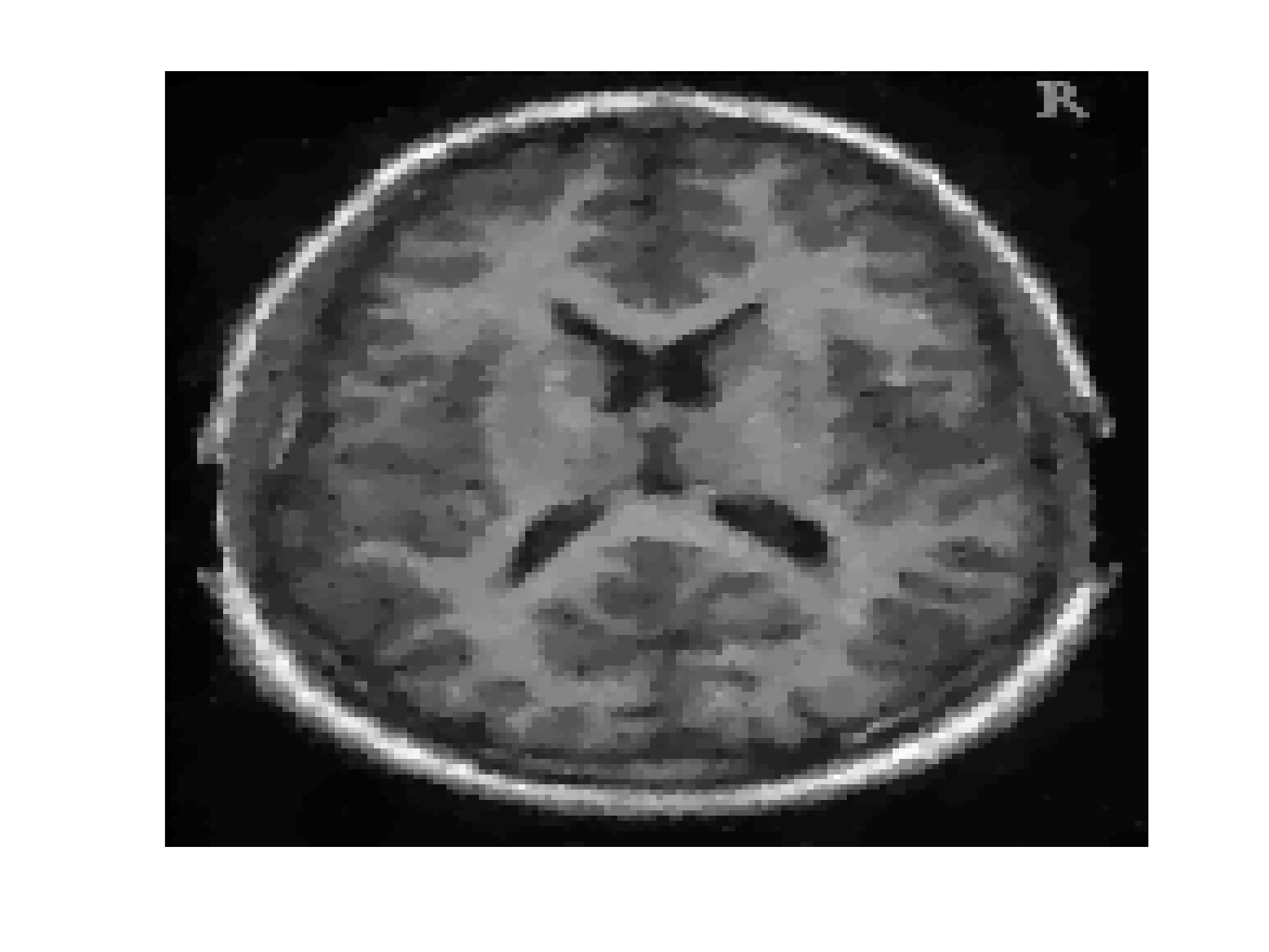}
\hspace{-0.3cm}
 \includegraphics[width=2cm,height=2cm]{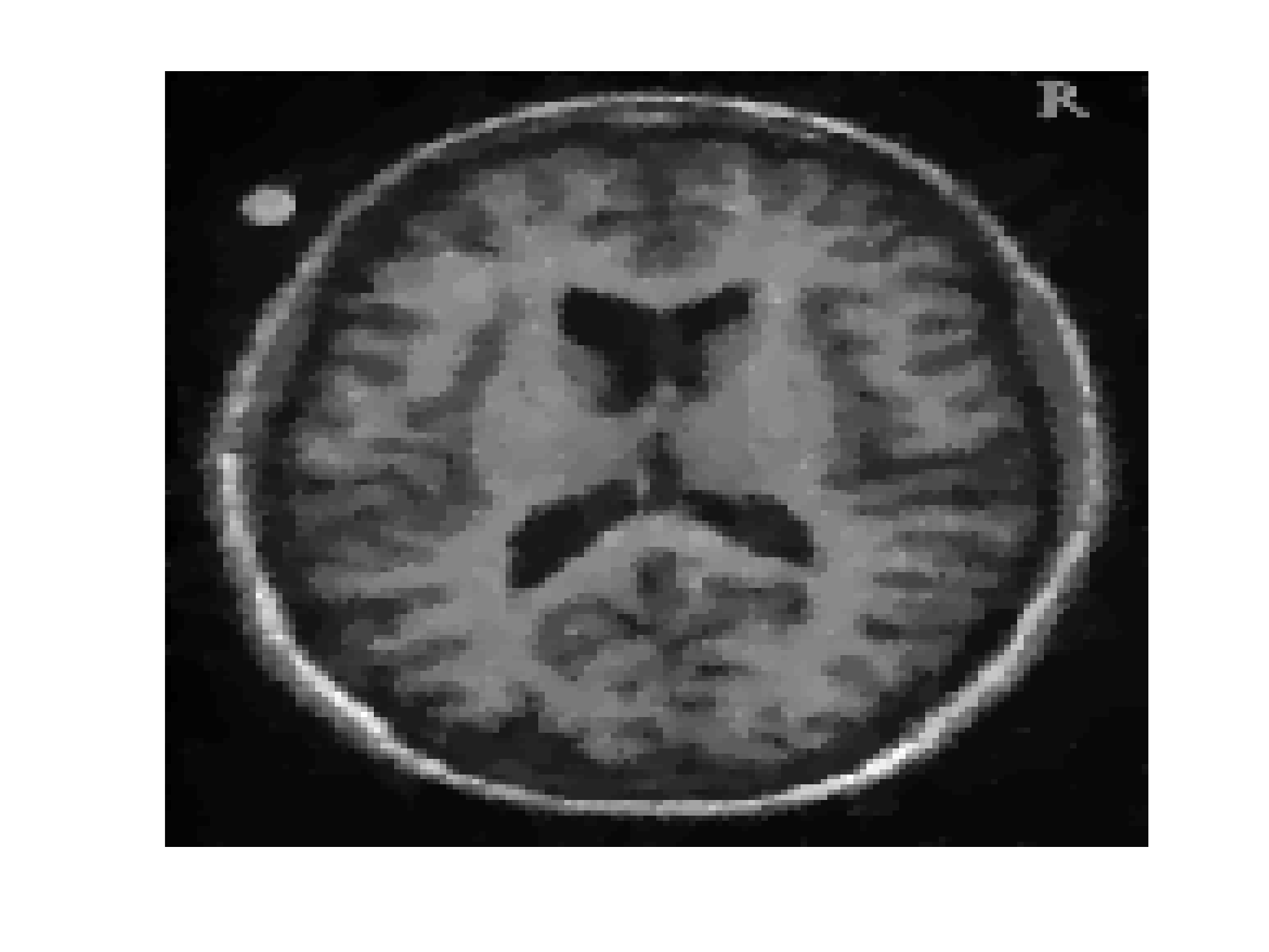}
\hspace{-0.3cm}
\includegraphics[width=2cm,height=2cm]{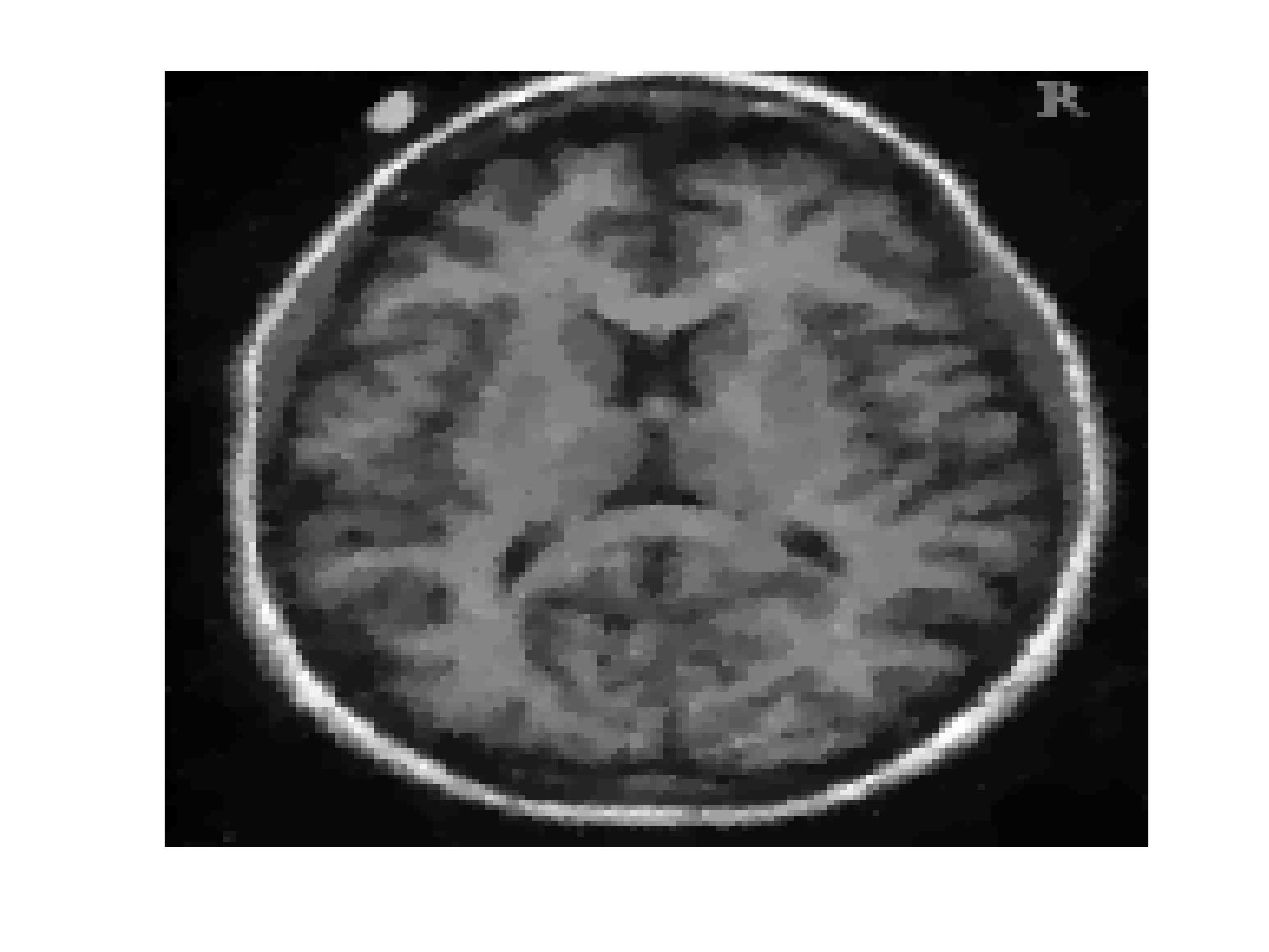}
\end{center}
\caption{Sample of $5$ images of OASIS MRI brain database: original images (upper row), noisy images (middle row) and optimal denoised images (bottom row),  $\hat{\lambda_S}=3280.5$.}
\label{fig:training set}
\end{figure}

Figure \ref{fig:training set} shows an example of database of brain images\footnote{OASIS online database, \texttt{http://www.oasis-brains.org/}.} together with the optimal denoised version obtained by Algorithm \ref{alg:dynsamp} for Gaussian noise estimation.

\begin{figure}[h!]
\begin{center}
\includegraphics[width=1cm,height=1cm]{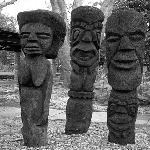}\hspace{0.2cm}
\includegraphics[width=1cm,height=1cm]{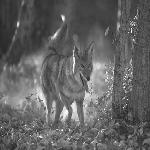}\hspace{0.2cm}
 \includegraphics[width=1cm,height=1cm]{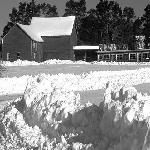}\hspace{0.2cm}
 \includegraphics[width=1cm,height=1cm]{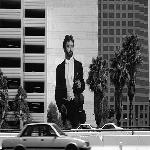}\hspace{0.2cm}
\includegraphics[width=1cm,height=1cm]{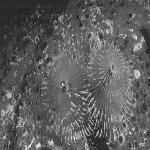}\hspace{0.2cm}
\includegraphics[width=1cm,height=1cm]{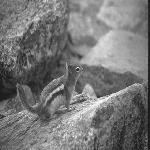}\hspace{0.2cm}
\includegraphics[width=1cm,height=1cm]{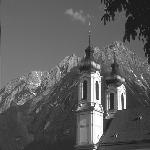}\hspace{0.2cm}
\includegraphics[width=1cm,height=1cm]{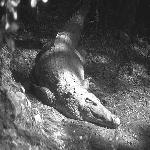}\hspace{0.2cm}
\includegraphics[width=1cm,height=1cm]{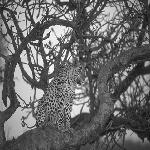}\hspace{0.2cm}
\includegraphics[width=1cm,height=1cm]{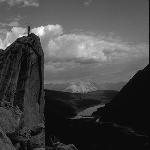} \\ \vspace{0.2cm}
\includegraphics[width=1cm,height=1cm]{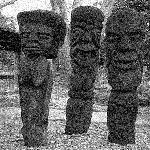}\hspace{0.2cm}
\includegraphics[width=1cm,height=1cm]{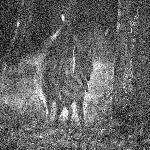}\hspace{0.2cm}
 \includegraphics[width=1cm,height=1cm]{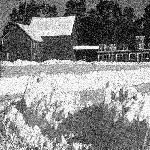}\hspace{0.2cm}
 \includegraphics[width=1cm,height=1cm]{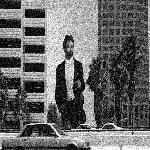}\hspace{0.2cm}
\includegraphics[width=1cm,height=1cm]{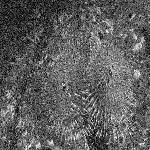}\hspace{0.2cm}
\includegraphics[width=1cm,height=1cm]{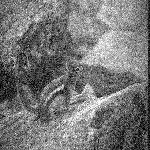}\hspace{0.2cm}
\includegraphics[width=1cm,height=1cm]{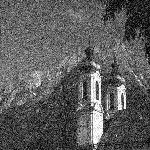}\hspace{0.2cm}
\includegraphics[width=1cm,height=1cm]{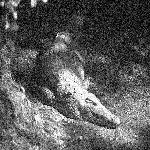}\hspace{0.2cm}
\includegraphics[width=1cm,height=1cm]{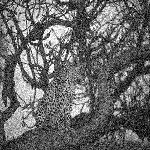}\hspace{0.2cm}
\includegraphics[width=1cm,height=1cm]{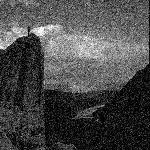} 
\end{center}
\caption{Noise-free and noisy versions of images from the Berkeley database. The Gaussian noise distribution is  $0$ mean and variance $\sigma^2=0.01$.}
\label{fig:berkeley_dat}
\end{figure}

In order to test the adaptability of our method to images which are very diverse between each other, we test our model for a very diversified database \footnote{Berkeley database, available online at: \texttt{http://www.eecs.berkeley.edu/Research/Projects/CS/vision/bsds/BSDS300/html/dataset/images.html}}, see Fig. \ref{fig:berkeley_dat}. From Table \ref{table_div_database} we can observe that increasing the size of the database, the estimation of the optimal parameter $\lambda^*$ may vary significantly, due to the diversity of images considered. This reflects the property of our approach to estimate the parameter $\lambda^*$ which is optimal with respect to the \emph{entire} database, cf. cost functional \eqref{optprob1gauss}.
%We report the results obtained in Figure \ref{fig:berkeley_res}.

\begin{table}[h!]
\centering
%\begin{tabular}{| c | c | } 
%    \hline
%$N$ & $\lambda^*_S$ \\ \hline
%$10$ & $2732.15$ \\ \hline
%$20$ & $2766.32$ \\ \hline
%$30$ & $2170.23$  \\ \hline
%$40$ & $2292.51$  \\ \hline
% \end{tabular}
 \begin{tabular}{| c | c | c | c | c | }
 \hline
 $\bm{K}$ & $10$ & $20$ & $30$ & $40$  \\ \hline
 $\bm{\lambda^*}$ & $2732.15$ & $2766.32$ & $2170.23$ & $2292.51$  \\ \hline
 \end{tabular}
 
\caption{Optimal $\lambda^*$ estimation for heterogeneous database, see Fig. \ref{fig:berkeley_dat}. The numerical value adapts to the diversity of the images considered.}
\label{table_div_database}
\end{table}

\begin{figure}[h!] 
\begin{center}
\includegraphics[height=3cm]{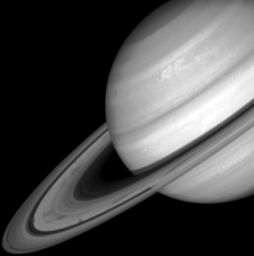} \quad
\includegraphics[height=3cm]{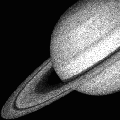} \quad 
\includegraphics[height=3cm]{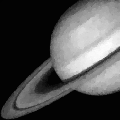}
\end{center}
\caption{Poisson denoising: Original (left), noisy (center) and optimal denoised (right) images. Parameters: $\gamma=1e3, \mu=1e-10, h=1/128, \eta=1e4$. Optimal weight: $\lambda^*=1013.76$.}
\label{fig:poisson}
\end{figure}

\paragraph{Poisson noise} As a second example, we consider the case of images corrupted by Poisson noise. The corresponding data fidelity in this case has been shown in \cite{sawatzky2009total} to be a KL-type fidelity defined as $\Phi(u)=u-f~ \log u$, which requires the additional condition for $u$ to be strictly positive. We enforce this constraint by using a standard penalty method and solve:
\begin{equation*} \label{eq: poiss upper level}
\min_{\lambda \geq 0}  ~\frac{1}{2} \|f_0 - u\|^2_{L^2} 
\end{equation*}
where $u$ is the solution of the minimisation problem:
\begin{equation} \label{eq: poiss lower level}
\min_{v>0} \left\{\frac{\mu}{2} \norm{\grad v}^2_{L^2} + |D v|(\Omega) +
  \lambda \int_\Omega  (v-f~ \log v )~dx + \frac{\eta}{2} \| \min(v,\delta) \|_{L^2}^2 \right\},
\end{equation}
where $\eta\gg 1$ is a penalty parameter enforcing the positivity constraint and $\delta\ll 1$ ensures strict positivity throughout the optimisation. After Huber-regularising the TV term using \eqref{def:huber},  we write the primal-dual form of the corresponding optimality condition for the optimisation problem \eqref{eq: poiss lower level} similarly as in \eqref{eq:primal-dual1}-\eqref{eq:primal-dual2} :
\begin{equation} \label{eq: poiss lower level PDE}
-\mu\Delta u - \mathrm{div}\,q  +\lambda ~(1-\frac{f}{u}) + \eta\chi_{\mathcal T_\gamma }~ u = 0,\quad q=\frac{\gamma\nabla u}{\max(\gamma|\nabla u|, 1)},
\end{equation}
where $\mathcal T_\delta$ is the active set $\mathcal T_\delta := \{ x \in \Omega: u(x)<\delta \}$.  We then design a modified SSN iteration solving \eqref{eq: poiss lower level PDE} similarly as described in Section \ref{sec:adjoint},  see \cite[Section 4]{de2013image} for more details. Figure \ref{fig:poisson} shows the optimal denoising result for the Poisson noise case in correspondence of the value $\lambda^*=1013.76$.

\paragraph{Spatially dependent weight}
We continue with an example where $\lambda$ is spatially-dependent. Specifically, we choose as parameter space $V= \{ v \in H^1(\Omega): \partial_n u=0 \text{ on }\Gamma \}$ in combination with a TV regulariser and a single Gaussian noise model. For this example the noisy image is distorted non-uniformly: A Gaussian noise with zero mean and variance $0.04$ is present on the whole image and an additional noise with variance $0.06$ is added on the area marked by red line.

Since the spatially dependent parameter does not allow to get rid of the positivity constraints in an automatic way, we solved the whole optimality system by means of the semismooth Newton method described in Section 3, combined with a Schwarz domain decomposition method. Specifically, we decomposed the domain first and apply the globalized Newton algorithm in each subdomain afterwards. The detailed numerical performance of this approach is reported in \cite{chungdelosreyes}.

\begin{figure}[!h] \centering
  \begin{tabular}{ccc}
\includegraphics[width=0.3\textwidth]{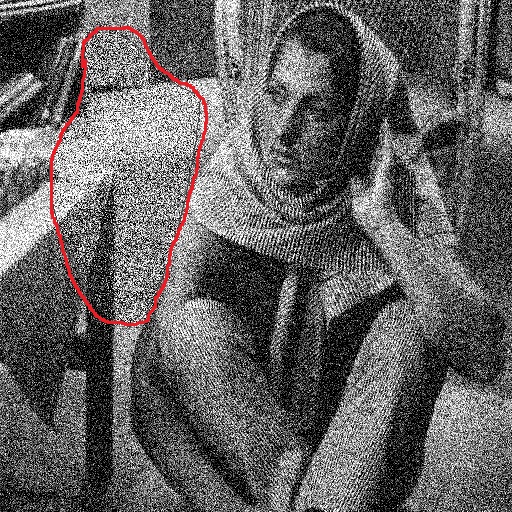} &    \includegraphics[width=0.3\textwidth]{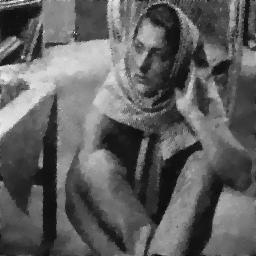} &
    \includegraphics[width=0.3\textwidth]{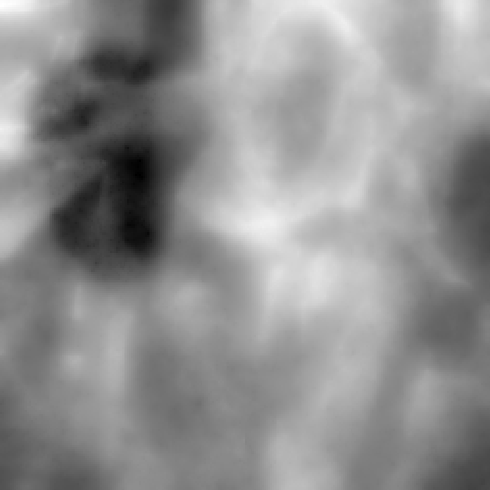}
  \end{tabular}
  \caption{Noisy image (left), denoised image (center) and intensity of $\lambda^*$ (right).}
  \label{fig_4th_rs_whole}
\end{figure}

The results are shown in Figure \ref{fig_4th_rs_whole} for the parameters $\mu = 1e-16$, $\gamma=25$ and $h = 1/500$. The values of $\lambda$ on whole domain are between $100.0$ to $400.0$. From the right image in Figure \ref{fig_4th_rs_whole} we can see the dependence of the optimal parameter $\lambda^*$ on the distribution of noise. As expected, at the high-level noise area in the input image, values of $\lambda^*$ are lower (darker area) than in the rest of the image.

\subsection{Multiple noise estimation}
In many applications, the acquired image may be possibly corrupted by different types of noise, each one corresponding to a different data fidelity term $\Phi_i$ weighted by a non-negative weighting $\lambda_i$. In this multiple noise case, we consider the following optimisation lower level problem:
\begin{equation*} 
\min_{u} \left\{\frac{\mu}{2} \|\grad u\|^2_{L^2} + |Du|(\Omega) +
  \int_\Omega \Psi(\lambda_1,\ldots,\lambda_M,\Phi_1(u),\ldots,\Phi_M(u))~dx \right\},
\end{equation*}
where the modelling function $\Psi$ combines the different fidelity terms $\Phi_i$ and weights $\lambda_i$ in order to deal with the multiple noise case. The case when $\Psi$ is a linear a linear combination of fidelities $\Phi_i$ with coefficients $\lambda_i$ is the one presented in the general model \eqref{eq:learn} and \eqref{eq:learn-numerical-single} and has been considered in \cite{de2013image}. In the following, we present also the case when $\Psi$ is an infimal-convolution operation of fidelities, as considered in \cite{lucainfimal}.

\paragraph{Impulse and Gaussian noise} 

 Motivated by some previous work in literature on the use of the infimal-convolution operation 
 \cite[Chapter 16]{bauschke2011convex} for image decomposition, cf. \cite{chambolle97image,vaggelisinfimal}, we consider in \cite{lucainfimal} the modelling of mixed noise distribution through such operation with the intent of obtaining an optimal denoised image thanks to the decomposition of the noise into its different components. In the case of combined Gaussian and impulse noise, the optimisation model reads:
$$
\min_{\lambda_1,\lambda_2 \geq 0}  ~\frac{1}{2} \|f_0 - u\|^2_{L^2} 
$$
where $u$ is the solution of the optimisation problem:
\begin{equation}\label{impgaussinfconv}
\min_{\substack{ v\in BV \\
 n\in L^2}}  \left\{\frac{\mu}{2}\|\grad v\|^2_{L^2}+|Dv|(\Omega) + 
\lambda_1\|n\|_{L^1}  + \lambda_2\|f-v-n\|_{L^2}^2 \right\},
\end{equation}
where $n$ represents the impulse noise component (and, as such, is treated using the $L^1$ norm) and the optimisation runs over $v$ and $n$. We use once again a single training pair $(f_0,f)$ 
%Similarly as in \eqref{eq:impulseVI}-\eqref{eq: impulse lower level}
and consider a Huber-regularisation depending on a parameter $\gamma$ for both the TV term and the $L^1$ norm in \eqref{impgaussinfconv}. The corresponding Euler-Lagrange equations are:
\begin{align*}
 -\mu\Delta u - \mathrm{div}\,\left(\frac{\gamma\nabla u}{\max(\gamma|\nabla u|, 1)} \right) -\lambda_2 (f-u-n) & = 0,\\
\lambda_1 \frac{\gamma\ n}{\max(\gamma|n|, 1)} -\lambda_2 (f-u-n) & = 0.
\end{align*}
Again, writing the equations above in a primal-dual form, we can write the modified SSN iteration and solve the optimisation problem with BFGS as described in Section \ref{sec:adjoint}.
% as follows:
%\begin{align}
%& -\mu \Delta \delta_u - \mathrm{div}\, \delta_q +\lambda_2 \delta_u + \lambda_2 \delta_{n} = - \left(-\mu \Delta u - \mathrm{div}\, q +\lambda_2 (f-u-n) \right) \label{eq:numerics imp gauss modif ssn iter 1}\\
%& \lambda_1\delta_p+\lambda_2\delta_u+\lambda_2\delta_{n} = - \left(\lambda_1p -\lambda_2 (f-u-n) \right) \label{eq:numerics imp gauss modif ssn iter 2} \\
%&\delta_q - \frac{\gamma \nabla \delta_u}{\max(1,\gamma|\nabla u|)}+ \chi_{\mathcal A_\gamma}  \gamma^2 \frac{\nabla u^T \nabla \delta_u}{\max(1,\gamma|\nabla u|)^2} \frac{q}{\max(1,|q|)} \nonumber \\&\hspace{6cm}=-q+ \frac{\gamma \nabla u}{\max(1,\gamma|\nabla u|)}, \label{eq:numerics gauss posson modif ssn iter 3} \\
%&\delta_p - \frac{\gamma\ \delta_{n}}{\max(1,\gamma|n|)}+ \chi_{\mathcal R_\gamma}  \gamma^2 \frac{n\ \delta_{n}}{\max(1,\gamma|n|)^2} \frac{p}{\max(1,|p|)} \nonumber \\&\hspace{6cm}=-p+ \frac{\gamma\  n}{\max(1,\gamma|n|)}, \label{eq:numerics gauss posson modif ssn iter 4}
%\end{align}
%for the increments $\delta_u, \delta_{n}, \delta_q$ and $\delta_p$ and where the active sets $\mathcal{A}_\gamma$ and $\mathcal{R}_\gamma$ are defined to be: $\mathcal{A}_\gamma=\{ x \in \Omega: \gamma|\nabla u(x)|\geq1 \}$ and $\mathcal{R}_\gamma=\{ x \in \Omega: \gamma|n(x)|\geq1 \}$.

\begin{figure}[h!] 
\begin{center}
\includegraphics[width=2cm]{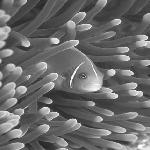} \hfill\includegraphics[width=2cm]{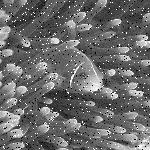} \hfill \includegraphics[width=2cm]{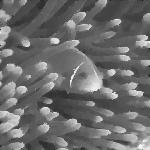}\hfill \includegraphics[width=2cm]{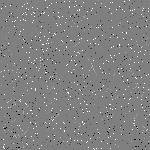}\hfill\includegraphics[width=2cm]{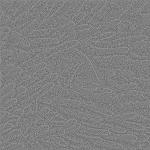}
\end{center}
\caption{Impulse-Gaussian denoising. From left to right: Original image, noisy image corrupted by impulse noise and Gaussian noise with mean zero and variance $0.005$, denoised image, impulse noise residuum and Gaussian noise residuum. Optimal parameters: $\lambda_1^*=351.23$ and $\lambda_2^*=5200.1$.}
\label{fig:infconv impulse gauss}
\end{figure}

In Figure \ref{fig:infconv impulse gauss} we present the results of the model considered. The original image $f_0$ has been corrupted with Gaussian noise of zero mean and variance $0.005$ and then a percentage of $5\%$ of pixels has been corrupted with impulse noise. The parameters have been chosen to be $\gamma =1e3$, $\mu=1e-15$ and the mesh step size $h=1/120$. The computed optimal weights are $\lambda_1^*=351.23$ and $\lambda_2^*=5200.1$. The results show the actual decomposition of the noise into its sparse and Gaussian components.

\paragraph{Gaussian and Poisson noise}  We consider now the optimisation problem with $\Phi_1(u)=|u-f|^2$ for the Gaussian noise component and $\Phi_2(u)=(u-f\log u)$ for the Poisson distributed one. We aim to determine the optimal weighting $(\lambda_1,\lambda_2)$ as follows:
$$
\min_{\lambda_1,\lambda_2 \geq 0}  ~\frac{1}{2} \|f_0 -u \|^2_{L^2} 
$$
subject to $u$ be the solution of:
\begin{equation}\label{gausspoissonorg}
\min_{v > 0} \left\{\frac{\mu}{2} \|\grad v\|^2_{L^2} + |Dv|(\Omega) + \frac{\lambda_1}{2} \|v-f\|^2_{L^2} + \lambda_2 \int_\Omega (v-f\log v)~ dx\right\},
\end{equation}
for one training pair $(f_0,f)$, where $f$ corrupted by Gaussian and Poisson noise. After Huber-regularising the Total Variation term in \eqref{gausspoissonorg}, we derive (formally) the following Euler-Lagrange equation 
\begin{align*}
& -\mu\Delta u - \mathrm{div}\,\left(\frac{\gamma\nabla u}{\max(\gamma|\nabla u|, 1)} \right) +\lambda_1 (u-f) + \lambda_2 (1-\frac{f}{u})  -\alpha = 0\\
& \alpha\cdot u = 0,
\end{align*}
with non-negative Lagrange multiplier $\alpha \in L^2(\Omega)$. As in \cite{sawatzky2009total} we multiply the first equation by $u$ and obtain
$$
u\cdot \left(-\mu\Delta u - \mathrm{div}\left(\frac{\gamma\nabla u}{\max(\gamma|\nabla u|, 1)} \right) +\lambda_1 (u-f) \right) + \lambda_2(u-f) = 0,
$$
where we have used the complementarity condition $\alpha\cdot u =0$. Next, the solution $u$ is computed iteratively by using a semismooth Newton type method combined with the outer BFGS iteration as above.

\begin{figure}[h!]
\begin{center}
\includegraphics[width=0.3\textwidth]{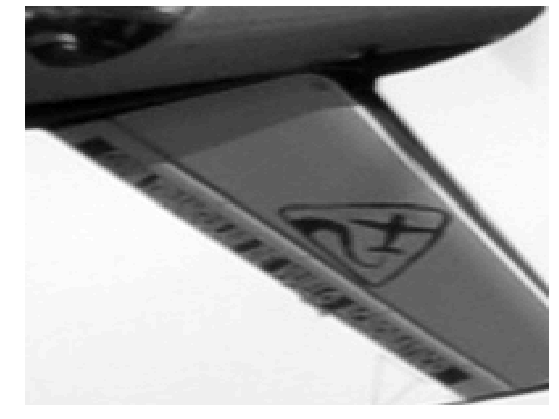} \hfill \includegraphics[width=0.3\textwidth]{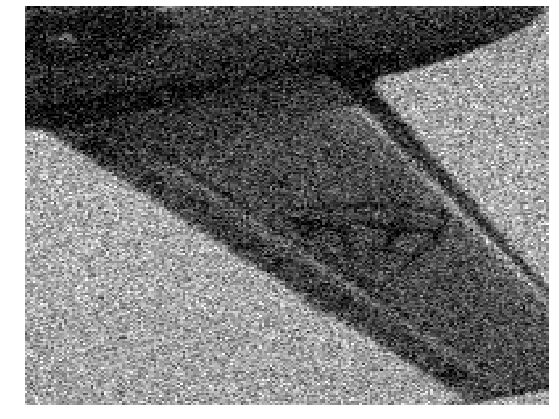} \hfill \includegraphics[width=0.3\textwidth]{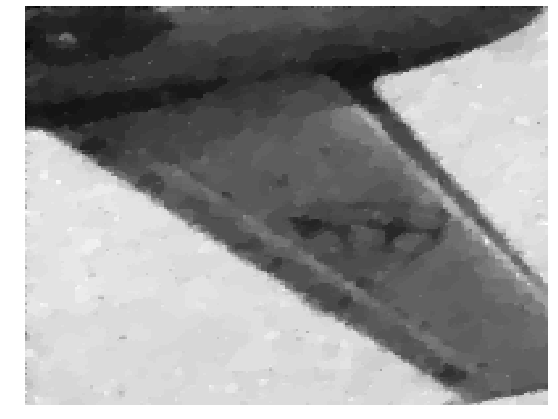}  
\end{center}
\caption{Poisson-Gaussian denoising: Original image (left), noisy image corrupted by Poisson noise and Gaussian noise with mean zero and variance $0.001$ (center) and denoised image (right). Optimal parameters $\lambda_1^*=1847.75$ and $\lambda_2^*=73.45$.}
\label{fig:experimentpoissongauss}
\end{figure}

In Figure \ref{fig:experimentpoissongauss} we show the optimisation result. The original image $f_0$ has been first corrupted by Poisson noise and then Gaussian noise was added, with zero mean and variance $0.001$. Choosing the parameter values to be $\gamma =100$ and $\mu=1e-15$, the optimal weights $\lambda_1^*=1847.75$ and $\lambda_2^*=73.45$ were computed on a grid with mesh size step $h=1/200$. 
%From Figure \ref{fig:experimentpoissongauss} a good match between the original and the denoised images can be observed. The noise appears to be successfully removed with the computed optimal weights.

\section{Conclusion and outlook}\label{sec:concl}
Machine learning approaches in image processing and computer vision have mostly developed in parallel to their mathematical analysis counterparts, which have variational regularisation models at their core. Variational regularisation techniques offer rigorous and intelligible image analysis -- which gives reliable and stable answers that provide us with insight in the constituents of the process and error estimates. This guarantee of giving a meaningful and stable result is crucial in most image processing applications, in biomedical and seismic imaging, in remote sensing and astronomy: provably giving an answer which is correct up to some error bounds is important when diagnosing patients, deciding upon a surgery or when predicting earthquakes. Machine learning methods, on the other hand, are extremely powerful as they learn from examples and are hence able to adapt to a specific task. The recent rise of deep learning gives us a glimpse on what is possible when intelligently using data to learn from. Todays (29 April 2015) search on a Google on `deep learning image' just gave 59,800,000 hits. Deep learning is employed for all kinds of image processing and computer vision tasks, with impressive results! The weak point of machine learning approaches, however, is that they generally cannot offer stability or error bounds, neither provide most of them understanding about the driving factors (e.g. the important features in images) that led to their answer.

In this paper we wanted to give an account to a recent realisation in mathematical image processing that a marriage between machine learning and variational regularisation might be interesting -- an attempt to bring together the Good from both worlds. In particular, we have discussed bilevel optimisation approaches in which optimal image regularisers and data fidelity terms are learned making use of a training set. We discussed the analysis of such a bilevel strategy in the continuum as well as their efficient numerical solution by quasi-Newton methods, and presented numerical examples for computing optimal regularisation parameters for TV, $\TGV^2$ and $ICTV$ denoising, as well as for deriving optimal data fidelity terms for TV image denoising for data corrupted with pure or mixed noise distributions.

Although the techniques presented in this article are mainly focused on denoising problems, the perspectives of using similar approaches in other image reconstruction problems (inpainting, segmentation, etc.) appear to be promising. Also the extension to color images deserves to be further studied.

Finally, there are still several open questions which deserve to be investigated in the future. Here a short list:
\begin{itemize}
\item Is it possible to obtain an optimality system for (P) by performing an asymptotic analysis when $\mu \rightarrow 0$?
 \item How to measure optimality? Are quality measures such as the
signal-to-noise ratio and generalisations thereof \cite{wang2004ssim} enough? Should one try to match characteristic expansions of the image such as Fourier or Wavelet expansions?
\cite{mumford2001stochastic}.
\item How to decide about the presence of a specific noise model? Is it possible to use sparse optimization for the automatic identification of one specific model? Can it be used to identify mixed models?
\end{itemize}

\section{Acknowledgments}

The original research behind this review has been supported by the King Abdullah University of Science and Technology (KAUST) Award No.~KUK-I1-007-43, the EPSRC grants Nr.~EP/J009539/1 ``Sparse \& Higher-order Image Restoration'', and Nr.~EP/M00483X/1 ``Efficient computational tools for inverse imaging problems'', the Escuela Politécnica Nacional de Quito under award PIS 12-14 and the MATHAmSud project SOCDE ``Sparse Optimal Control of Differential Equations''. C. Cao and T.~Valkonen have also been supported by Prometeo scholarships of SENESCYT (Ecuadorian Ministry of Higher Education, Science, Technology and Innovation).

%% If there is an additional footnote on page 1, place ``\makethankshere'' subsequent
%% to that footnote and use the class option ``nothanks''.

\medskip

    \noindent\textbf{A data statement for the EPSRC}
    The data leading to this \emph{review} publication will be made available, as appropriate, as part of the original publications that this work summarises.

 \providecommand{\homesiteprefix}{http://iki.fi/tuomov/mathematics}
  \providecommand{\eprint}[1]{\href{http://arxiv.org/abs/#1}{arXiv:#1}}

%\bibliography{abbrevs2,bib}

\begin{thebibliography}{10}
\providecommand{\url}[1]{\texttt{#1}}
\providecommand{\urlprefix}{URL }
\providecommand{\eprint}[2][]{\url{#2}}

\bibitem{allard2008total}
W.~Allard, Total variation regularization for image denoising, {I}. {G}eometric
  theory. \emph{SIAM J. Math. Anal.} 39 (2008) 1150--1190.

\bibitem{ambcosdal96}
L.~Ambrosio, A.~Coscia and G.~Dal~Maso, Fine Properties of Functions with
  Bounded Deformation. \emph{Arch. Ration. Mech. Anal.} 139 (1997) 201--238.

\bibitem{ambrosio2000fbv}
L.~Ambrosio, N.~Fusco and D.~Pallara, \emph{Functions of Bounded Variation and
  Free Discontinuity Problems}, Oxford University Press (2000).

\bibitem{baus2014fully}
F.~Baus, M.~Nikolova and G.~Steidl, Fully smoothed L1-TV models: Bounds for the
  minimizers and parameter choice. \emph{J. Math. Imaging Vision} 48 (2014)
  295--307.

\bibitem{bauschke2011convex}
H.~H. Bauschke and P.~L. Combettes, \emph{Convex Analysis and Monotone Operator
  Theory in Hilbert Spaces}, CMS Books in Mathematics, Springer (2011).

\bibitem{benning2011higher}
M.~Benning, C.~Brune, M.~Burger and J.~Müller, Higher-Order {TV}
  Methods—Enhancement via {Bregman} Iteration. \emph{J. Sci. Comput.} 54
  (2013) 269--310.

\bibitem{benning2012ground}
M.~Benning and M.~Burger, Ground states and singular vectors of convex
  variational regularization methods. \emph{Methods and Applications of
  Analysis} 20 (2013) 295--334, \eprint{1211.2057}.

\bibitem{biegler2011large}
L.~Biegler, G.~Biros, O.~Ghattas, M.~Heinkenschloss, D.~Keyes, B.~Mallick,
  L.~Tenorio, B.~van Bloemen~Waanders, K.~Willcox and Y.~Marzouk,
  \emph{Large-scale inverse problems and quantification of uncertainty}, volume
  712, John Wiley \& Sons (2011).

\bibitem{bredies2012total}
K.~Bredies and M.~Holler, A total variation-based JPEG decompression model.
  \emph{SIAM J. Imaging Sci.} 5 (2012) 366--393.

\bibitem{bredies2011tgv}
K.~Bredies, K.~Kunisch and T.~Pock, Total Generalized Variation. \emph{SIAM J.
  Imaging Sci.} 3 (2011) 492--526.

\bibitem{bredies2013properties}
K.~Bredies, K.~Kunisch and T.~Valkonen, Properties of $L^1$-$\mbox{TGV}^2$: The
  one-dimensional case. \emph{J. Math. Anal Appl.} 398 (2013) 438--454.

\bibitem{sampta2011tgv}
K.~Bredies and T.~Valkonen, Inverse problems with second-order total
  generalized variation constraints, in: \emph{Proc.~{SampTA} 2011} (2011).

\bibitem{bui2008model}
T.~Bui-Thanh, K.~Willcox and O.~Ghattas, Model reduction for large-scale
  systems with high-dimensional parametric input space. \emph{SIAM J. Sci.
  Comput.} 30 (2008) 3270--3288.

\bibitem{papoutsellis2014}
M.~Burger, J.~M{\"u}ller, E.~Papoutsellis and C.-B. Sch{\"o}nlieb, Total
  Variation Regularisation in Measurement and Image space for PET
  reconstruction. \emph{Inverse Problems} 10 (2014).

\bibitem{vaggelisinfimal}
M.~Burger, K.~Papafitsoros, E.~Papoutsellis and C.-B. Sch{\"o}nlieb, Infimal
  convolution regularisation functionals on BV and $L^p$ spaces. Part I: The
  finite p case (2015), submitted.

\bibitem{byrd2012sample}
R.~H. Byrd, G.~M. Chin, J.~Nocedal and Y.~Wu, Sample size selection in
  optimization methods for machine learning. \emph{Math. Program.} 134 (2012)
  127--155.

\bibitem{impulsegauss2008}
J.-F. Cai, R.~H. Chan and M.~Nikolova, Two-phase approach for deblurring images
  corrupted by impulse plus gaussian noise. \emph{Inverse Probl. Imaging} 2
  (2008) 187--204.

\bibitem{calatronidynamic}
L.~Calatroni, J.~C. De~los Reyes and C.-B. Sch{\"o}nlieb, Dynamic sampling
  schemes for optimal noise learning under multiple nonsmooth constraints, in:
  \emph{System Modeling and Optimization}, 85--95, Springer Verlag (2014).

\bibitem{lucainfimal}
L.~Calatroni, J.~C.~D. los Reyes and C.-B. Sch{\"o}nlieb, Learning the optimal
  Total Variation denoising model for multiple noise distributions, in
  preparation.

\bibitem{caselles2007discontinuity}
V.~Caselles, A.~Chambolle and M.~Novaga, The discontinuity set of solutions of
  the TV denoising problem and some extensions. \emph{Multiscale Model. Simul.}
  6 (2007) 879--894.

\bibitem{chambolle97image}
A.~Chambolle and P.-L. Lions, Image recovery via total variation minimization
  and related problems. \emph{Numer.~Math.} 76 (1997) 167--188.

\bibitem{Chen2012}
Y.~Chen, T.~Pock and H.~Bischof, Learning $\ell_1$-based analysis and synthesis
  sparsity priors using bi-level optimization, in: \emph{Workshop on Analysis
  Operator Learning vs. Dictionary Learning, NIPS 2012} (2012).

\bibitem{chen2014}
Y.~Chen, R.~Ranftl and T.~Pock, Insights into analysis operator learning: From
  patch-based sparse models to higher-order MRFs. \emph{Image Processing, IEEE
  Transactions on}  (2014), to appear.

\bibitem{chen_cvpr2015}
Y.~Chen, W.~Yu and T.~Pock, {On learning optimized reaction diffusion processes
  for effective image restoration}, in: \emph{IEEE Conference on Computer
  Vision and Pattern Recognition} (2015), to appear.

\bibitem{chungdelosreyes}
C.~V. Chung and J.~C. De~los Reyes, Learning optimal spatially-dependent
  regularization parameters in total variation image restoration, in
  preparation.

\bibitem{chung2011designing}
J.~Chung, M.~Chung and D.~P. O'Leary, Designing optimal spectral filters for
  inverse problems. \emph{SIAM J. Sci. Comput.} 33 (2011) 3132--3152.

\bibitem{chung2014optimal}
J.~Chung, M.~I. Espa{\~n}ol and T.~Nguyen, Optimal Regularization Parameters
  for General-Form Tikhonov Regularization. \emph{arXiv preprint
  arXiv:1407.1911}  (2014).

\bibitem{cichocki2002adaptive}
A.~Cichocki, S.-i. Amari et~al., \emph{Adaptive Blind Signal and Image
  Processing}, John Wiley Chichester (2002).

\bibitem{costantini2004virtual}
R.~Costantini and S.~Susstrunk, Virtual sensor design, in: \emph{Electronic
  Imaging 2004}, 408--419, International Society for Optics and Photonics
  (2004).

\bibitem{de2015numerical}
J.~C. De~los Reyes, \emph{Numerical PDE-Constrained Optimization}, Springer
  (2015).

\bibitem{de2013image}
J.~C. De~los Reyes and C.-B. Sch{\"o}nlieb, Image denoising: Learning the noise
  model via Nonsmooth {PDE}-constrained optimization. \emph{Inverse Probl.
  Imaging} 7 (2013).

\bibitem{reyes2015b}
J.~C. De~los Reyes, C.-B. Sch{\"o}nlieb and T.~Valkonen, Optimal parameter
  learning for higher-order total variation regularisation models, in
  preparation.

\bibitem{reyes2015a}
J.~C. De~los Reyes, C.-B. Sch{\"o}nlieb and T.~Valkonen, The structure of
  optimal parameters for image restoration problems, in preparation.

\bibitem{domke2012generic}
J.~Domke, Generic methods for optimization-based modeling, in:
  \emph{International Conference on Artificial Intelligence and Statistics},
  318--326 (2012).

\bibitem{dong2011automated}
Y.~Dong, M.~Hinterm{\"u}ller and M.~M. Rincon-Camacho, Automated regularization
  parameter selection in multi-scale total variation models for image
  restoration. \emph{J. Math. Imaging Vision} 40 (2011) 82--104.

\bibitem{engl1996regularization}
H.~W. Engl, M.~Hanke and A.~Neubauer, \emph{Regularization of Inverse
  Problems}, volume 375, Springer (1996).

\bibitem{evans2002inverse}
S.~N. Evans and P.~B. Stark, Inverse problems as statistics. \emph{Inverse
  Problems} 18 (2002) R55.

\bibitem{Foigausspoiss}
A.~Foi, Clipped noisy images: Heteroskedastic modeling and practical denoising.
  \emph{Signal Processing} 89 (2009) 2609 -- 2629, special Section: Visual
  Information Analysis for Security.

\bibitem{frick2012statistical}
K.~Frick, P.~Marnitz, A.~Munk et~al., Statistical multiresolution Dantzig
  estimation in imaging: Fundamental concepts and algorithmic framework.
  \emph{Electronic Journal of Statistics} 6 (2012) 231--268.

\bibitem{gilboa2014total}
G.~Gilboa, A total variation spectral framework for scale and texture analysis.
  \emph{SIAM J. Imaging Sci.} 7 (2014) 1937--1961.

\bibitem{haber2010numerical}
E.~Haber, L.~Horesh and L.~Tenorio, Numerical methods for the design of
  large-scale nonlinear discrete ill-posed inverse problems. \emph{Inverse
  Problems} 26 (2010) 025002.

\bibitem{haber2003learning}
E.~Haber and L.~Tenorio, Learning regularization functionals -- a supervised
  training approach. \emph{Inverse Problems} 19 (2003) 611.

\bibitem{hintermuller2013subspace}
M.~Hinterm{\"u}ller and A.~Langer, Subspace Correction Methods for a Class of
  Nonsmooth and Nonadditive Convex Variational Problems with Mixed $L^1/L^2$
  Data-Fidelity in Image Processing. \emph{SIAM J. Imaging Sci.} 6 (2013)
  2134--2173.

\bibitem{hintermuller2006infeasible}
M.~Hinterm{\"u}ller and G.~Stadler, An Infeasible Primal-Dual Algorithm for
  Total Bounded Variation--Based Inf-Convolution-Type Image Restoration.
  \emph{SIAM J. Sci. Comput.} 28 (2006) 1--23.

\bibitem{hintermuller2014bilevel}
M.~Hinterm{\"u}ller and T.~Wu, Bilevel Optimization for Calibrating Point
  Spread Functions in Blind Deconvolution (2014), preprint.

\bibitem{huang2012optimal}
H.~Huang, E.~Haber, L.~Horesh and J.~K. Seo, Optimal Estimation Of
  {L1}-regularization Prior From A Regularized Empirical Bayesian Risk
  Standpoint. \emph{Inverse Probl. Imaging} 6 (2012).

\bibitem{idier2013bayesian}
J.~Idier, \emph{Bayesian approach to inverse problems}, John Wiley \& Sons
  (2013).

\bibitem{poissongauss2013}
A.~Jezierska, E.~Chouzenoux, J.-C. Pesquet and H.~Talbot, {A Convex Approach
  for Image Restoration with Exact Poisson-Gaussian Likelihood}, Technical
  report (2013).

\bibitem{kaipio2006statistical}
J.~Kaipio and E.~Somersalo, \emph{Statistical and computational inverse
  problems}, volume 160, Springer Science \& Business Media (2006).

\bibitem{kingsbury2001complex}
N.~Kingsbury, Complex wavelets for shift invariant analysis and filtering of
  signals. \emph{Applied and Computational Harmonic Analysis} 10 (2001)
  234--253.

\bibitem{Klatzer2015}
T.~Klatzer and T.~Pock, {Continuous Hyper-parameter Learning for Support Vector
  Machines}, in: \emph{Computer Vision Winter Workshop (CVWW)} (2015).

\bibitem{knoll2010second}
F.~Knoll, K.~Bredies, T.~Pock and R.~Stollberger, Second order total
  generalized variation {(TGV)} for {MRI}. \emph{Magnetic Resonance in
  Medicine} 65 (2011) 480--491.

\bibitem{kolehmainen2011marginalization}
V.~Kolehmainen, T.~Tarvainen, S.~R. Arridge and J.~P. Kaipio, Marginalization
  of uninteresting distributed parameters in inverse problems—application to
  diffuse optical tomography. \emph{International Journal for Uncertainty
  Quantification} 1 (2011).

\bibitem{kunisch2013bilevel}
K.~Kunisch and T.~Pock, A bilevel optimization approach for parameter learning
  in variational models. \emph{SIAM J. Imaging Sci.} 6 (2013) 938--983.

\bibitem{krtv}
J.~Lellmann, D.~Lorenz, C.-B. Sch{\"o}nlieb and T.~Valkonen, Imaging with
  {K}antorovich-{R}ubinstein discrepancy. \emph{SIAM J. Imaging Sci.} 7 (2014)
  2833--2859, \eprint{1407.0221}.

\bibitem{mairal2009online}
J.~Mairal, F.~Bach, J.~Ponce and G.~Sapiro, Online dictionary learning for
  sparse coding, in: \emph{Proceedings of the 26th Annual International
  Conference on Machine Learning}, 689--696, ACM (2009).

\bibitem{MBPSZ08}
J.~Mairal, B.~F, J.~Ponce, G.~Sapiro and A.~Zisserman, Discriminative learned
  dictionaries for local image analysis. \emph{CVPR}  (2008).

\bibitem{meyer2002oscillating}
Y.~Meyer, \emph{Oscillating patterns in image processing and nonlinear
  evolution equations}, AMS (2001).

\bibitem{mumford2001stochastic}
D.~Mumford and B.~Gidas, Stochastic models for generic images. \emph{Quarterly
  of Applied Mathematics} 59 (2001) 85--112.

\bibitem{natterer2001mathematical}
F.~Natterer and F.~W{\"u}bbeling, \emph{Mathematical Methods in Image
  Reconstruction}, Monographs on Mathematical Modeling and Computation Vol 5,
  Philadelphia, PA: SIAM) (2001).

\bibitem{nikolova2004variational}
M.~Nikolova, A variational approach to remove outliers and impulse noise.
  \emph{J. Math. Imaging Vision} 20 (2004) 99--120.

\bibitem{nocedal2006numerical}
J.~Nocedal and S.~Wright, \emph{Numerical Optimization}, Springer Series in
  Operations Research and Financial Engineering, Springer (2006).

\bibitem{ochs_ssvm2015}
P.~Ochs, R.~Ranftl, T.~Brox and T.~Pock, {Bilevel Optimization with Nonsmooth
  Lower Level Problems}, in: \emph{International Conference on Scale Space and
  Variational Methods in Computer Vision (SSVM)} (2015), to appear.

\bibitem{OF96}
B.~Olshausen and D.~Field, Emergence of simple-cell receptive field properties
  by learning a sparse code for natural images. \emph{Nature} 381 (1996)
  607--609.

\bibitem{papafitsoros2013study}
K.~Papafitsoros and K.~Bredies, A study of the one dimensional total
  generalised variation regularisation problem. \emph{arXiv preprint
  arXiv:1309.5900}  (2013).

\bibitem{peyre2011learning}
G.~Peyr{\'e} and J.~M. Fadili, Learning analysis sparsity priors.
  \emph{Sampta'11}  (2011).

\bibitem{Ranftl_GCPR2014}
R.~Ranftl and T.~Pock, A Deep Variational Model for Image Segmentation, in:
  \emph{36th German Conference on Pattern Recognition (GCPR)} (2014).

\bibitem{ring2000structural}
W.~Ring, Structural Properties of Solutions to Total Variation Regularization
  Problems. \emph{ESAIM: Math. Model. Numer. Anal.} 34 (2000) 799--810.

\bibitem{rockafellar-wets-va}
R.~T. Rockafellar and R.~J.-B. Wets, \emph{Variational Analysis}, Springer
  (1998).

\bibitem{roth2005fields}
S.~Roth and M.~J. Black, Fields of experts: A framework for learning image
  priors, in: \emph{Computer Vision and Pattern Recognition, 2005. CVPR 2005.
  IEEE Computer Society Conference on}, volume~2, 860--867, IEEE (2005).

\bibitem{rudin1992nonlinear}
L.~I. Rudin, S.~Osher and E.~Fatemi, Nonlinear total variation based noise
  removal algorithms. \emph{Physica D: Nonlinear Phenomena} 60 (1992) 259--268.

\bibitem{sawatzky2009total}
A.~Sawatzky, C.~Brune, J.~Müller and M.~Burger, Total Variation Processing of
  Images with Poisson Statistics, in: \emph{Computer Analysis of Images and
  Patterns}, \emph{Lecture Notes in Computer Science}, volume 5702, Edited by
  X.~Jiang and N.~Petkov, 533--540, Springer Berlin Heidelberg (2009).

\bibitem{scharf1991statistical}
L.~L. Scharf, \emph{Statistical Signal Processing}, volume~98, Addison-Wesley
  Reading, MA (1991).

\bibitem{schmidt2014shrinkage}
U.~Schmidt and S.~Roth, Shrinkage fields for effective image restoration, in:
  \emph{Computer Vision and Pattern Recognition (CVPR), 2014 IEEE Conference
  on}, 2774--2781, IEEE (2014).

\bibitem{starckmurtagh1994}
J.-L. Starck, F.~D. Murtagh and A.~Bijaoui, Image restoration with noise
  suppression using a wavelet transform and a multiresolution support
  constraint (1994).

\bibitem{tadmor2004multiscale}
E.~Tadmor, S.~Nezzar and L.~Vese, A multiscale image representation using
  hierarchical (BV, L 2) decompositions. \emph{Multiscale Model. Simul.} 2
  (2004) 554--579.

\bibitem{tappen2007utilizing}
M.~F. Tappen, Utilizing variational optimization to learn {M}arkov random
  fields, in: \emph{Computer Vision and Pattern Recognition, 2007. CVPR'07.
  IEEE Conference on}, 1--8, IEEE (2007).

\bibitem{temam1985mpp}
R.~Temam, \emph{Mathematical problems in plasticity}, Gauthier-Villars (1985).

\bibitem{unser1995texture}
M.~Unser, Texture classification and segmentation using wavelet frames.
  \emph{Image Processing, IEEE Transactions on} 4 (1995) 1549--1560.

\bibitem{unser2013unifying}
M.~Unser and N.~Chenouard, A unifying parametric framework for 2D steerable
  wavelet transforms. \emph{SIAM J. Imaging Sci.} 6 (2013) 102--135.

\bibitem{tuomov-jumpset}
T.~Valkonen, The jump set under geometric regularisation. {Part 1}: Basic
  technique and first-order denoising. \emph{SIAM J. Math. Anal.}  (2015),
  accepted, \eprint{1407.1531}.

\bibitem{vardi1985statistical}
Y.~Vardi, L.~Shepp and L.~Kaufman, A statistical model for positron emission
  tomography. \emph{Journal of the American Statistical Association} 80 (1985)
  8--20.

\bibitem{viola2012unifying}
F.~Viola, A.~Fitzgibbon and R.~Cipolla, A unifying resolution-independent
  formulation for early vision, in: \emph{Computer Vision and Pattern
  Recognition (CVPR), 2012 IEEE Conference on}, 494--501, IEEE (2012).

\bibitem{wang2004ssim}
Z.~Wang, A.~C. Bovik, H.~R. Sheikh and E.~P. Simoncelli, Image quality
  assessment: From error visibility to structural similarity. \emph{IEEE Trans.
  Image Processing} 13 (2004) 600--612.

\bibitem{YSM10}
G.~Yu, G.~Sapiro and S.~Mallat, Image modeling and enhancement via structured
  sparse model selection. \emph{Proc. IEEE Int. Conf. Image Processing}
  (2010).

\end{thebibliography}
% Alternatively, you may use the standard bibliography environment:       %
%-------------------------------------------------------------------------%

%\begin{thebibliography}{9}
%  \bibitem{Lange:2007}
%  T.~Lange and I.~E. Shparlinski,
%  Distribution of some sequences of points on elliptic curves,
%  \emph{J. Math. Cryptol.} \textbf{1} (2007), 1--11.
%\end{thebibliography}

\end{document}